\documentclass[a4paper]{amsart}
\usepackage[utf8]{inputenc}

\usepackage{bbm}
\usepackage{amsmath}
\usepackage{amsfonts}
\usepackage{amssymb}
\usepackage{amsthm}
\usepackage[foot]{amsaddr}
\usepackage{enumerate}
\usepackage{xcolor}
\usepackage{mathrsfs}
\usepackage{stmaryrd}
\usepackage{graphicx}
\usepackage{caption}
\usepackage{subcaption}
\usepackage[a4paper,margin=1in,headsep=25pt,footskip=35pt]{geometry}
\usepackage{fancyhdr}
\usepackage{comment}
\usepackage{hyperref}

\SetSymbolFont{stmry}{bold}{U}{stmry}{m}{n}

\newtheorem{thm}{Theorem}
\newtheorem{cor}[thm]{Corollary}
\newtheorem{lem}[thm]{Lemma}
\newtheorem{prop}[thm]{Proposition}

\newtheorem{defn}[thm]{Definition}
\newtheorem*{defn*}{Definition}
\newtheorem{assump}[thm]{Assumption}

\newtheorem{rem}[thm]{Remark}

\newcommand{\deq}{\mathrel{\mathop:}=}
\newcommand{\e}[1]{\mathrm{e}^{#1}}
\newcommand{\R} {\mathbb{R}}
\newcommand{\C} {\mathbb{C}}

\newcommand{\N} {\mathbb{N}}

\newcommand{\E} {\mathbb{E}}

\newcommand{\adj}{^{*}} 
\newcommand{\tp}{^{\intercal}}

\newcommand{\dist} {\mathrm{dist}}

\DeclareMathOperator{\tr}{tr}
\DeclareMathOperator{\Tr}{Tr}

\DeclareMathOperator{\supp}{supp}

\DeclareMathOperator{\re}{\mathrm{Re}}
\DeclareMathOperator{\im}{\mathrm{Im}}

\newcommand{\caA}{{\mathcal A}}

\newcommand{\caC}{{\mathcal C}}
\newcommand{\caD}{{\mathcal D}}
\newcommand{\caE}{{\mathcal E}}

\newcommand{\caH}{{\mathcal H}}

\newcommand{\caM}{{\mathcal M}}

\newcommand{\caS}{{\mathcal S}}

\newcommand{\caX}{{\mathcal X}}
\newcommand{\caY}{{\mathcal Y}}

\newcommand{\bbC}{{\mathbb C}}

\newcommand{\bbF}{{\mathbb F}}

\newcommand{\bbR}{{\mathbb R}}
\newcommand{\bbS}{{\mathbb S}}

\newcommand{\frc}{\mathfrak{c}}

\newcommand{\fre}{{\mathfrak e}}

\newcommand{\frs}{{\mathfrak s}}

\newcommand{\frX}{{\mathfrak X}}
\newcommand{\frY}{{\mathfrak Y}}

\newcommand{\bsa}{{\boldsymbol a}}
\newcommand{\bsb}{{\boldsymbol b}}
\newcommand{\bsc}{{\boldsymbol c}}

\newcommand{\bsh}{{\boldsymbol h}}

\newcommand{\bsm}{{\boldsymbol m}}
\newcommand{\bsn}{{\boldsymbol n}}

\newcommand{\bss}{{\boldsymbol s}}

\newcommand{\bsu}{{\boldsymbol u}}

\newcommand{\bsw}{{\boldsymbol w}}
\newcommand{\bsx}{{\boldsymbol x}}
\newcommand{\bsy}{{\boldsymbol y}}
\newcommand{\bsz}{{\boldsymbol z}}

\newcommand{\bsF}{{\boldsymbol F}}

\newcommand{\bsL}{{\boldsymbol L}}

\newcommand{\bsQ}{{\boldsymbol Q}}

\newcommand{\rmd}{\mathrm{d}}

\newcommand{\rmf}{\mathrm{f}}

\newcommand{\rmi}{\mathrm{i}}

\newcommand{\wt}{\widetilde}
\newcommand{\ol}{\overline}

\newcommand{\mr}{\mathring}
\newcommand{\wh}{\widehat}
\newcommand{\beq}{ \begin{equation} }
	\newcommand{\eeq}{ \end{equation} }
\newcommand{\beqs}{\begin{equation*}}
	\newcommand{\eeqs}{\end{equation*}}

\newcommand{\lone}{\mathbbm{1}} 

\newcommand{\dd}{\mathrm{d}}
\newcommand{\ii}{\mathrm{i}}

\renewcommand{\P}{\mathbb{P}}
\newcommand{\SC}{\mathrm{sc}}

\newcommand{\bsiota}{\boldsymbol{\iota}}

\newcommand{\llbra}{\llbracket}
\newcommand{\rrbra}{\rrbracket}

\newcommand\norm[1]{\Vert#1\Vert}
\newcommand\Norm[1]{\left\Vert#1\right\Vert}

\newcommand\Absv[1]{\left\vert#1\right\vert}
\newcommand\absv[1]{\vert#1\vert}

\newcommand\brkt[1]{\langle#1\rangle}
\newcommand\Brkt[1]{\left\langle#1\right\rangle}
\newcommand\bbrktt[1]{\llbra #1\rrbra}

\numberwithin{equation}{section} 
\numberwithin{thm}{section}

\newcommand{\nc}{\normalcolor}

\title{On the spectral edge of non-{H}ermitian random matrices}
\author[]{Andrew Campbell$^{\dagger}$}
\email{andrew.campbell@ist.ac.at}

\author[]{Giorgio Cipolloni$^{\ddagger}$}
\email{gcipolloni@arizona.edu}

\author[]{L\'{a}szl\'{o} Erd\H{o}s$^{\dagger}$}
\email{lerdos@ist.ac.at}

\author[]{Hong Chang Ji$^{*}$}
\email{hji56@wisc.edu}

\address{$^{\dagger}$Institute of Science and Technology Austria}
\address{$^{\ddagger}$ Department of Mathematics, University of Arizona }
\address{$^{*}$ Department of Mathematics, University of Wisconsin-Madison}
\thanks{$^\dagger$ Supported by ERC Advanced Grant ``RMTBeyond" No. 101020331}

\begin{document}
	\maketitle
	
\begin{abstract}
For general large non-Hermitian random matrices $X$ and deterministic deformation matrices $A$, we prove that the local eigenvalue statistics of $A+X$ close to the typical edge points of its spectrum are universal. Furthermore, we show that under natural assumptions on $A$ the spectrum of $A+X$ does not have outliers at a distance larger than the natural fluctuation scale of the eigenvalues. As a consequence, the number of eigenvalues in each component of $\mathrm{Spec}(A+X)$ is deterministic.
\end{abstract}
	
	\tableofcontents
	
\section{Introduction}

	\subsection{Main results}\label{sec:main results intro} Consider a large $N\times N$ deterministic matrix $A$ (called a \emph{deformation matrix}) and a large random matrix $X$ with independent and identically distributed (i.i.d.) centered  entries, which we refer to as an \textit{IID} matrix. It was established in \cite{Erdos-Ji2023circ} that the eigenvalue density of the \textit{deformed} non-Hermitian random matrix $A+X$ at the edge of its spectrum exhibits either a sharp cutoff behavior, which we refer to as the \textit{sharp} edge, or there may be a small subset of the spectral edge with quadratic decay, which we refer to as the \textit{quadratic} edge. Motivated by the well known analogous principle in Hermitian random matrix theory, our first main result, Theorem \ref{thm:main correlations}, establishes that the local eigenvalue statistics at sharp edge points are completely governed only by
	  the symmetry class of the matrix, i.e.\ real or complex. Specifically, under very mild assumptions on $A$, we prove full universality for the scaling limit of eigenvalue correlation functions localized around a sharp edge point. These limiting statistics are explicitly given by the well known special case when $A=0$ and $X$ has i.i.d.\ real or complex Gaussian entries with variance $\frac{1}{N}$, known as the real or complex \emph{Ginibre ensembles}, respectively. We additionally show that Theorem \ref{thm:main correlations} implies universality for matrices with a non-trivial variance profile or even with certain correlation structure.
	
	 In the pure IID case, i.e.\ when $A=0$, this universality principle (both in the bulk and at the edge of the spectrum) of local eigenvalue correlation functions was previously established for matrices whose entries have the same first four moments as the entries of a Ginibre matrix (commonly referred to as four moment matching) by Tao and Vu \cite{Tao-Vu2015}. Full universality at the spectral edge,   i.e.  with the four moment matching condition of \cite{Tao-Vu2015} removed, is due to \cite{Cipolloni-Erdos-Schroder2021}.  The  special case of complex Ginibre $X$ and normal $A$ was established only very recently by Liu and Zhang at the quadratic edge in \cite{Liu-Zhang2023} and the sharp edge in \cite{Liu-Zhang2024}.
	
	In addition to our first main result determining the local eigenvalue fluctuations, we consider the question of whether there can be any eigenvalues of $A+X$ outside the bulk spectrum beyond the scale of these fluctuations. Answering this question in the negative, 
	our second main result Theorem \ref{theo:noput} is a \emph{no-outlier} theorem: We establish that, with very high probability, $A+X$ has no eigenvalues beyond the natural fluctuation scale at the edge of the spectrum, under a natural assumption on $A$ that excludes outliers of $A$ itself. Previously Bordenave and Capitaine \cite[Theorem 1.4]{Bordenave-Capitaine2016} showed that deformed non-Hermitian random matrices have no outlying eigenvalues of distance on a macroscopic scale from the bulk spectrum, whereas we  achieve the same result on the optimal scale.
	
	\begin{figure}[h]
		\centering
		\begin{subfigure}{0.48\textwidth}
			\includegraphics[width=1.1\linewidth]{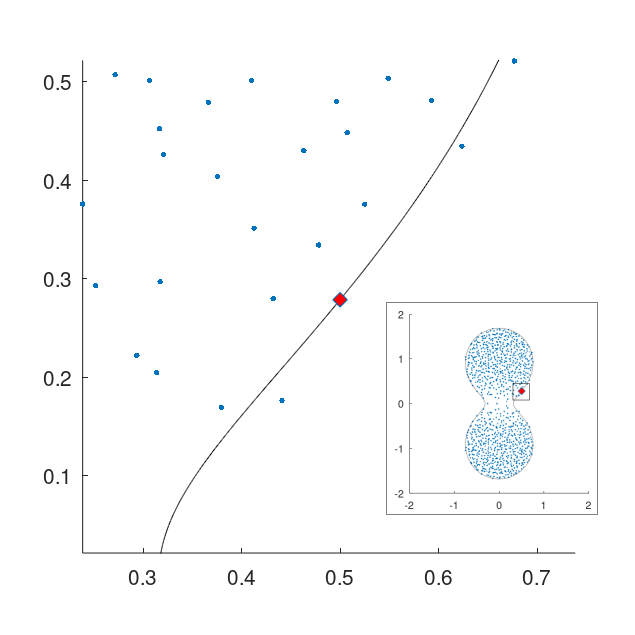}
			\caption{Eigenvalues (blue) around a base point (red) at the edge (black)\label{fig:univ}}
		\end{subfigure}
		\\
		\begin{subfigure}{0.48\textwidth}
			\includegraphics[width=1.1\linewidth]{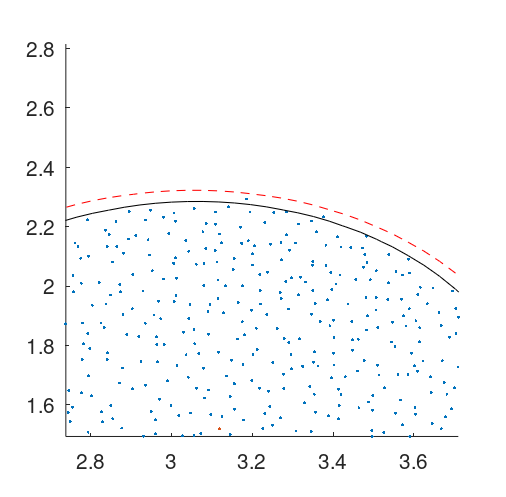}
			\caption{No outlier beyond fluctuation scale (red); see (c) for full spectrum.\label{fig:f_sc}}
		\end{subfigure}
		\,
		\begin{subfigure}{0.48\textwidth}
			\includegraphics[width=1.1\linewidth]{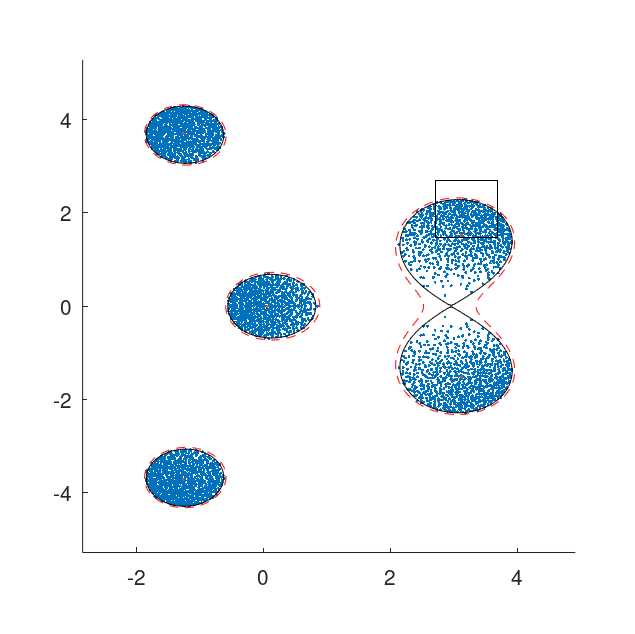}
			\caption{Clusters of eigenvalues (blue) and fluctuation scale (red)\label{fig:clust}}
		\end{subfigure}
	\caption{Illustrations of the main results}\label{fig:main}
	\end{figure}
	Our final main result, Corollary \ref{cor:clustrig}, concerns the components of $\mathrm{Spec}(A+X)$. For general $A$, $\mathrm{Spec}(A+X)$ can concentrate in \textit{clusters} on multiple disconnected domains in $\C$. Corollary \ref{cor:clustrig} asserts that each cluster of $\mathrm{Spec}(A+X)$ has a deterministic number of eigenvalues with very high probability. We refer to this phenomenon as \textit{cluster-rigidity}, providing the two-dimensional analogue of the \emph{band rigidity} of Hermitian random matrices proved in \cite{Alt-Erdos-Kruger-Schroder2020}. See Figure \ref{fig:main} for an illustration of the main results.

	\subsection{History and related results}\label{sec:intro history} 	For Hermitian random matrices the study of local eigenvalue statistics was initiated by the pioneering work of Wigner, and the universality of these statistics regardless of the exact distribution of the matrix entries is now well understood in many general ensembles. These local statistics of the eigenvalues are universal in the sense that they are determined only by the local density of the eigenvalues and the symmetry class of the matrix, i.e.\ real symmetric, complex Hermitian, etc. For mean field Hermitian random matrix ensembles it was established in \cite{Alt-Erdos-Kruger2020}  that there are only three possible local shapes (bulk/edge/cusp) of the eigenvalue density. In all three cases very general universality results have been established over the last decade or so: the bulk \cite{Ajanki-Erdos-Kruger2017,Johansson2001,Lee-Schnelli2018,Lee-Schnelli-Stetler-Yau2016,Sodin2010,Bourgade-Yau-Yin2020,Erdos-Kruger-Schroder2019,Erdos-Peche-Ramirez-Schlein-Yau2010,Erdos-Schlein-Yau2011,Erdos-Schnelli2017,Landon-Yau2017Conv,Pastur-Scherbina2008,Tao-Vu2011,Huang-Landon-Yau2015,Aggarwal-Lopatto-Yau2021} (i.e.\ locally flat spectral density), the edge \cite{Alt-Erdos-Kruger-Schroder2020,Bourgade-Erdos-Yau2014,Bourgade-Erdos-Yau-Yin2017,Landon-Yau2017,Lee-Schnelli2015,Soshnikov1999,Tao-Vu2010CMP,Pastur-Shcherbina2003,Ji-Park2021} (i.e.\ square-root decay of the spectral density), and the cusp \cite{Cipolloni-Erdos-Kruger-Schroder2019,Erdos-Kruger-Schroder2020,Hachem-Hardy-Najim2016} (i.e.\ cube root decay of the spectral density). 
	
	For non-Hermitian random matrices, studying the spectra is significantly more difficult, and much less is known about the universality of local eigenvalue statistics. Following the intuition from the Hermitian world, the first step to understanding local statistics is to understand the local shape of the eigenvalue densities of such matrices. For an IID matrix $X$ the limiting eigenvalue density is uniform on the unit disk, known as the \emph{ circular law}, and hence a point in the spectrum is either in the bulk, or the density has a sharp cutoff at that point. As discussed in Section \ref{sec:main results intro}, \cite{Erdos-Ji2023circ} established the non-Hermitian analogue of \cite{Alt-Erdos-Kruger2020} classifying the potential eigenvalue densities at the edge. In particular \cite[Remark 2.10]{Erdos-Ji2023circ} shows that the set of quadratic points is an analytic sub-manifold of the entire edge, 
	 i.e. typically consists of finitely many points. In this paper we thus focus on universality 
	near the ubiquitous  sharp edge points  for brevity, the case of quadratic edges is deferred to future work.
	
	The spectra of non-Hermitian operators are unstable with respect to small perturbations, and hence the size of the resolvent $\left(X-z \right)^{-1}$ of a non-Hermitian matrix is affected not just by the spectrum of $X$ but also its pseudospectrum. Another obstacle is that for non-normal operators the absence of the spectral theorem makes defining any spectral density more difficult. Fortunately, Brown \cite{Brown1986} defined what is now known as the \textit{Brown measure} of an operator in a von Neumann algebra, and this measure serves as our non-normal analogue of the spectral density. By the ``edge" of $A+X$ we mean truly the boundary of the support of the {Brown} measure of $A+\bsx$, where $\bsx$ is a \textit{circular} element in a von Neumann algebra equipped with a faithful normal trace $\brkt{\cdot}$.  This edge was identified in the  realm 
	 of free probability by \cite{Belinschi-Yin-Zhong2024} and \cite{Zhong2021} as the set of solutions to $\Brkt{\absv{A-z}^{-2} }=1$. For the purposes of interpreting our current work all one needs to know is that $\bsx$ is the appropriate infinite dimensional limit of $X$, and we refer the interested reader to, for example, \cite{Mingo-Speicher2017} for a thorough introduction to free probability.
	
	For Ginibre matrices the local eigenvalue statistics at the spectral edge were explicitly identified for the complex case in \cite{Ginibre1965} and for the more complicated real case in \cite{Borodin-Sinclair2009,Forrester-Nagao2007}; 
	 see also  \cite{Byun-Forrester2022GinUE,Byun-Forrester2023GinoE} and the references therein for a more complete history on eigenvalue statistics of Ginibre ensembles.   Beyond Ginibre ensembles, 
	 there are two prominent generalizations; the \textit{invariant} ensembles and matrices with independent, but not necessarily Gaussian entries like IID matrices. 
	
	The invariant ensembles are random normal matrices  with a density proportional to $\e{-2N\tr Q(M)}$ over the space of $N\times N$ complex normal matrices $M$ for some real potential function $Q$, whose eigenvalue distributions are also known as 2D Coulomb gases. Hedenmalm and Wennman \cite{Hedenmalm-Wennman2021} proved that at the spectral edge the eigenvalues exhibit complex Ginibre edge statistics for a family of real analytic potentials $Q$, see also the related earlier work of \cite{Ameur-Kang-Makarov2019}.  Local eigenvalue statistics at the edge of other non-Hermitian random matrices with Gaussian entries or from invariant ensembles have also been considered in various works including \cite{Akemann-Byun-Kang2022,Bothner-Little2022,Byun-Noda2024,Liu-Zhang2022,Ameur-Kang-Makarov2019}. However, assuming that the matrix has Gaussian entries or unitarily invariant distribution allows for explicit  algebraic formulas on the joint density function and techniques that do not generalize to other entry distributions. 
	
	 In this paper we consider the generalization of independent entries, for which the lack of algebraic formulas and techniques developed for invariant ensembles makes the situation significantly harder.  Beyond  \cite{Tao-Vu2015} and \cite{Cipolloni-Erdos-Schroder2021} little is known about edge universality of non-Hermitian random matrices. Very recently, the problem of bulk universality for matrices with i.i.d. entries was resolved in the breakthrough papers; for the complex case in \cite{Maltsev-Osman2023} by Maltsev and Osman, and for the real case in \cite{Dubova-Yang2024} by Dubova and Yang and in \cite{Osman2024} by Osman. All these works use very different methods compared to this manuscript. More precisely, using the multi--resolvent local laws from \cite{Cipolloni-Erdos-Schroder2024,Cipolloni-Erdos-Henheik-Schroder2023arXiv} as an input, they rely on the partial Schur decomposition to access the $k$--point correlation functions of the eigenvalues. The derivation of the correlation functions  uses supersymmetric (SUSY) techniques that heavily depend on the symmetry class.  The real case is much more complicated than 
	 the complex case  as a consequence of the symmetry of the spectrum with respect to the real axis. Neither these multi-resolvent local laws nor the derivation of the $k$--point correlation functions in the edge regime are known at the moment. We thus follow a simpler strategy which solely relies on the single resolvent local law. Moreover, our proof allows for the real and complex cases to be treated in parallel  with the same effort. 
	
	Extending the edge universality from standard IID matrices \cite{Cipolloni-Erdos-Schroder2021} to arbitrarily deformed ones, $A+X$, is not only a natural generalization, but also allows us  to establish edge universality for other natural non--Hermitian  matrices with nontrivial variance profile and even correlations. Thus our present result is the non--Hermitian counterpart of the Tracy-Widom edge universality established  for deformed Wigner matrices first in \cite{Lee-Schnelli2015} and then in \cite{Alt-Erdos-Kruger-Schroder2020} in full generality. Here we stress only that the ability to cover such  a wide array of cases with Theorem \ref{thm:main correlations} is thanks to our very mild assumptions on $A$. 
 There is a very general idea behind this approach:  the deformed model $A+X$ can be used to study very broad
classes of  random matrix ensembles, including certain correlations, as long as each matrix element has a sizeable variance.  In fact, to prove local laws, central limit theorems for linear eigenvalue statistics or local spectral universality for general matrix ensembles
 it is basically sufficient to prove the same result for $A+X^{\rm Gin}$ and then use some version of the \emph{Green function
 comparison} (GFT). 
	See Section~\ref{sec:examples} for further details  of this general idea for the current edge universality question.

	As a side note we mention that the extremal eigenvalue statistics also exhibit a universality phenomenon; for Ginibre ensembles they were proven to be Gumbel by \cite{Rider2003,Rider-Sinclair2014} and extended to full universality for IID matrices $X$ in the series of papers \cite{Cipolloni-Erdos-Schroder-Xu2022,Cipolloni-Erdos-Schroder-Xu2023,Cipolloni-Erdos-Xu2022,Cipolloni-Erdos-Xu2023}.

	 As far as outliers are concerned,  in Section \ref{sec:main results intro} we already mentioned 
	\cite{Bordenave-Capitaine2016} which 
	 proved that if $A$ does not have any outliers in its spectrum, then $A+X$ also does not have any outliers on a macroscopic scale. This result was an extension of a result of Tao \cite[Theorem 1.4]{Tao2013} which gave the same no outlier result for the special case of $A=0$. Our result Theorem \ref{theo:noput} improves \cite[Theorem 1.4]{Bordenave-Capitaine2016} from an order one result to the optimal scale given by the natural  local  fluctuations of the eigenvalues.
	With \cite{Bordenave-Capitaine2016} already describing the existence of outliers inherited from outliers of $A$, our result Theorem \ref{theo:noput}  provides the no-outlier counterpart on the optimal scale and resolves  the question of when $A+X$ does or does not have outliers in the spectrum beyond the natural scale. In an independent work posted on arXiv on the same day as the first version of this manuscript, Alt and Kr\"uger \cite{Alt-Kruger2024} showed that $A+X$ does not have outliers on a macroscopic scale when $X$ has a variance profile and $A$ is diagonal. A much stronger result, showing the non existence of outliers on the optimal fluctuation scale, follows from the techniques developed in this paper
	if $X$ is \emph{flat}  (see Remark~\ref{rem:compak} for a more detailed discussion). We point out that  \cite{Alt-Kruger2024} also gives an explicit description of the limiting eigenvalue density by solving the Matrix Dyson Equation (MDE) associated to these models. We also mention that there is a very large literature on  outliers (when they exist) 
	  that we do not discuss here; we  refer the interested readers to \cite{Baik-BenArous-Peche2005, Belinschi-Bercovici-Capitaine2017,Belinschi-Capitaine2017,Bigot-Male2021,Belinschi-Bordenave-Capitaine-Cebron2021,Benaych-Rochet2016,Basak-Zeitouni2020,Coston-ORourke-Wood2020,Nonnenmacher-Vogel2021} and references therein.

	\subsection{About the proofs}\label{sec:intro proof} In this section we provide a summary of the proofs and highlight the methodological novelties.  
	As is often the case in the study of non-Hermitian random matrices (with  a few  notable exceptions such as \cite{Maltsev-Osman2023,Osman2024,Bordenave-Chafai-Gracia-Zelada2022}) we use Girko's Hermitization \cite{Girko1984} formula to express the spectral statistics of $A+X$ in terms of traces of the resolvent matrix $G$ of the \textit{Hermitization} of $A+X$, which is defined by 
	\begin{equation}\label{eq:herm}
	H^{z}\deq \begin{pmatrix}
		0 & A+X-z\\
		(A+X-z)^* & 0
	\end{pmatrix},\quad z\in\C.
\end{equation} 
 The advantage of considering the Hermitization of a matrix is that it can be related to linear statistics of the non--Hermitian matrix via Girko's celebrated Hermitization formula:
\begin{equation}
	\label{eq:girkof}
	\int_{\C}F(z)\dd\rho_{A+X}(z)=-\frac{1}{4\pi N}\int_{\C}\Delta F(z) \log\Absv{\det H^{z}}\dd^2 z,
\end{equation}
where $\rho_{A+X}$ is the normalized empirical eigenvalue measure of $A+X$ and $F$ is a nice test function. We point out that the equality in \eqref{eq:girkof} can easily be derived by performing two integration by parts and using the two simple identities\footnote{Here we used that $\log(\cdot)$ is the Green's function of the Laplacian on $\C$, i.e. that $\Delta\log|z_0-\cdot|=-2\pi\delta_{z_0}$.}
\[
\frac{1}{2}\log\Absv{\det H^{z}}= \log\Absv{\det (X+A-z)}, \qquad\quad -\frac{1}{2\pi N}\Delta \log\Absv{\det (X+A-\cdot)}=\rho_{A+X}(\cdot).
\]
From \eqref{eq:girkof}, simple functional calculus allows one to express $\int F\dd\rho_{A+X}$ in terms of the resolvent $G^{z}=(H^{z}-\ii\eta)^{-1}$, and this observation by Girko has become the backbone of non-Hermitian random matrix theory. We then use this formula in conjunction with a \textit{local law}, i.e.\ an approximation of $G^{z}$ by the deterministic solution to a \textit{matrix Dyson equation} (MDE). With Girko's formula and this deterministic approximation available we can then compare the eigenvalue statistics of $A+X$ to the known statistics of a Ginibre matrix $X^{\mathrm{Gin}}$ through a Green function comparison (GFT). 
	
	To prove the optimal local law for the Hermitization $H^{z}$ we apply a bootstrapping argument, which we refer to as the \textit{iterated} \textit{Zig-Zag} method. This is a repeated combination of the characteristic flow method (\textit{Zig}) followed by a GFT (\textit{Zag}). This Zig-Zag argument has previously been used to prove multi-resolvent local laws in \cite{Cipolloni-Erdos-Schroder2024,Cipolloni-Erdos-Xu2023, Cipolloni-Erdos-Henheik2023,Cipolloni-Erdos-Henheik2024,Erdos-Riabov2024}. Crucially these previous works always had the single resolvent local law available as an a priori input, hence they could complete the process in only one or two Zig-Zag steps. One of the novelties of our paper is  the first self-consistent proof of a single resolvent local law 
	from scratch  using the Zig-Zag approach.  However,  lacking the  a priori input, we need to use a large (but finite) number of Zig-Zags in the process. We remark the Zig step has been previously used for a single resolvent local law in \cite{Adhikari-Landon2023,Bourgade2022,Adhikari-Huang2020,Huang-Landon2019}. Our local law covers both possible spectral edges of $A+X$, and hence so do our no-outlier and cluster-rigidity results. 
		 We mention that a clever self-consistent way of proving the GFT part of the isotropic
	  local law has already  appeared in \cite{Knowles-Yin2016}, where it was used to 
	  handle  the bulk and edge regimes for deformed Wigner and sample covariance matrices. However, 
	  as an input for the GFT argument, 
	  this paper still  relied on an isotropic 
	  local law for the analogous Gaussian model exploiting its unitary invariance property.
	  We use the convenient characteristic flow method, that was not available at the time of \cite{Knowles-Yin2016}, 
	  to prove local laws directly, using only the trivially available global laws as the sole input.
	  Moreover, the Zig-Zag method is very robust; here we apply it in the most critical cusp
	  regime without any additional complication.
	     
	We next discuss our proof of universality. The comparison between $A+X$ and $X^{\mathrm{Gin}}$ can be naturally broken into two steps by inserting the intermediate ensemble $A+X^{\mathrm{Gin}}$, both of which are carried out via GFT. We compare $A+X$ and $A+X^{\mathrm{Gin}}$ by evolving $X$ along an Ornstein--Uhlenbeck flow, and the comparison between $A+X^{\mathrm{Gin}}$ and $X^{\mathrm{Gin}}$ follows a continuous flow of the deformation. We focus below on the latter step since the former step is relatively standard. As this step is sensitive to the symmetry class, we first explain the complex case and then handle the additional difficulties in the real case.
	
\begin{figure}
	\centering
	\begin{subfigure}{0.32\textwidth}
		\includegraphics[width=\linewidth]{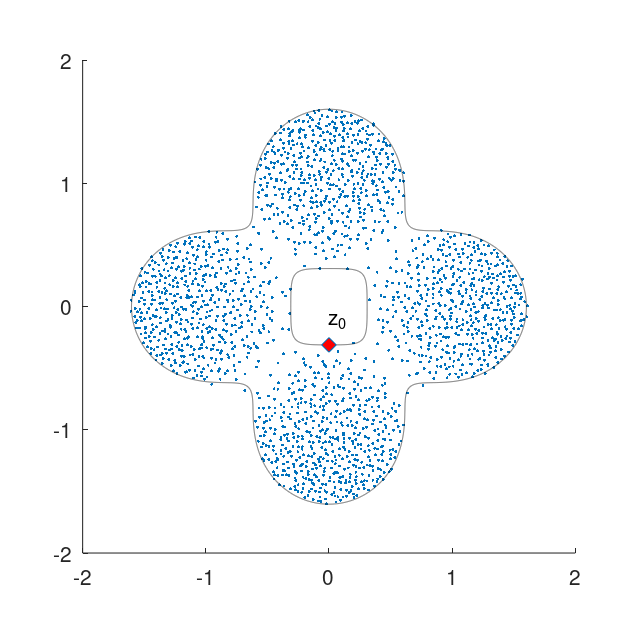}
	\end{subfigure}
\begin{subfigure}{0.32\textwidth}
	\includegraphics[width=\linewidth]{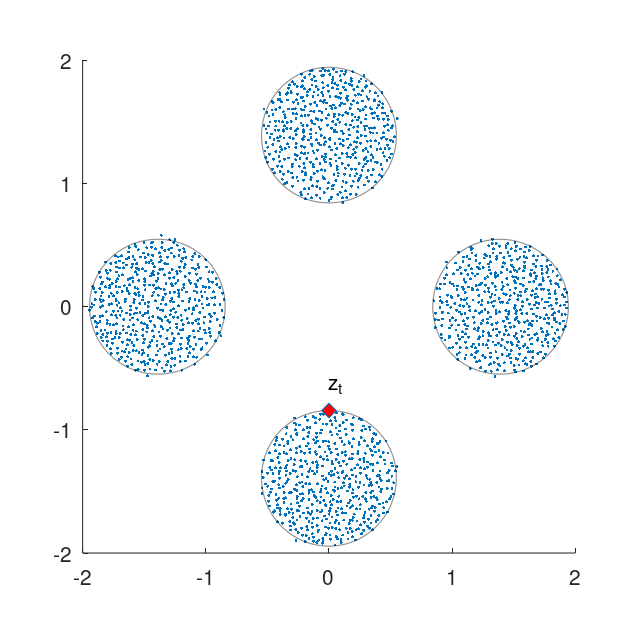}
\end{subfigure}
\begin{subfigure}{0.32\textwidth}
	\includegraphics[width=\linewidth]{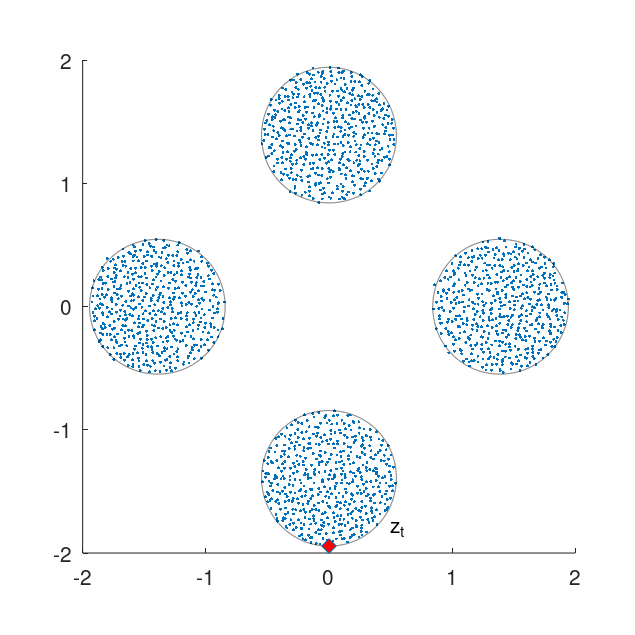}
\end{subfigure}	\\
\begin{subfigure}{0.32\textwidth}
	\includegraphics[width=\linewidth]{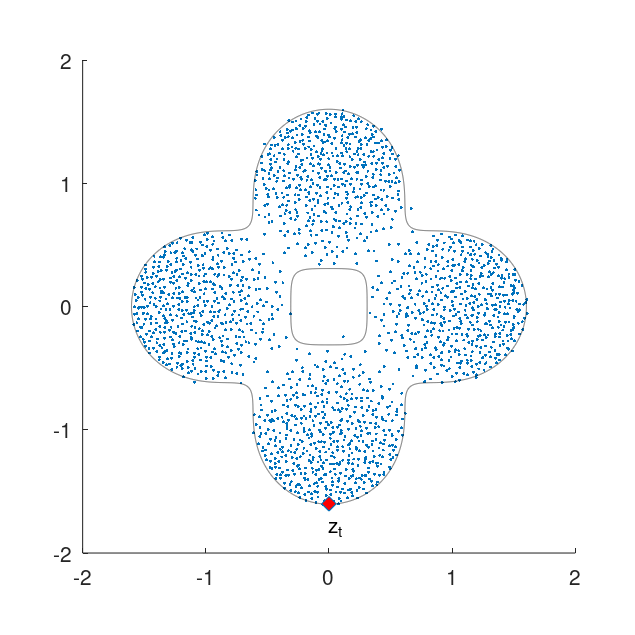}
\end{subfigure}
\begin{subfigure}{0.32\textwidth}
	\includegraphics[width=\linewidth]{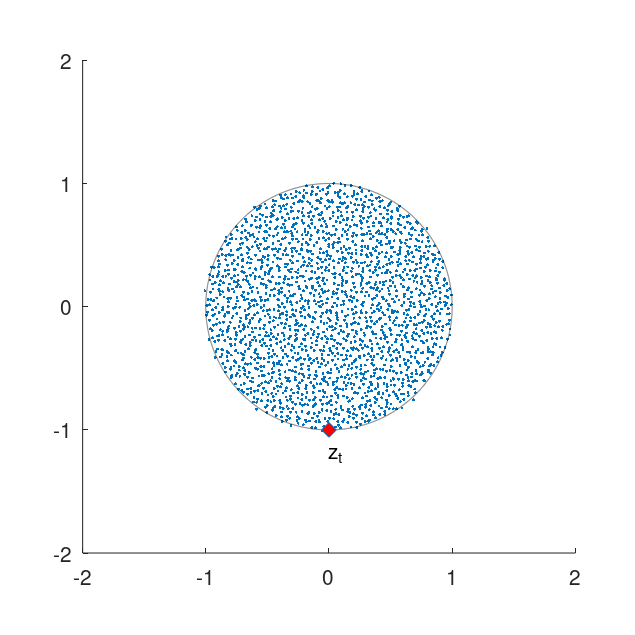}
\end{subfigure}
	\caption{Eigenvalues of the path $\caA_{t}+X$ and the evolution $z_t$ of the base point}\label{fig:comp_path}
\end{figure}

	In the complex case, to compare $A+X^{\mathrm{Gin}}$ and $X^{\mathrm{Gin}}$, we evolve $A-z$ continuously along a carefully chosen curve $A_t-z_t$ for $t\in[0,1]$ such that this $z_t$ remains a sharp edge of $A_t+X^{\mathrm{Gin}}$ and that $A_{0}=A$ and $A_{1}=0$. Along this curve, we show that the time derivatives of the eigenvalue and singular value statistics are sufficiently small.
	The key technical challenge in this comparison is that the naive estimates of the time derivatives available from the local law are not sufficient, and we need to identify non-trivial cancellations within the derivatives.
	These cancellations are revealed by introducing  and properly matching  \emph{four} time-dependent real parameters; one  in order to keep $z_{t}$ at  a sharp edge; another one   to keep  the (1D) local density of singular values unchanged; and two to keep the (2D) local density of eigenvalues constant. 
	Our approach is reminiscent of Lee and Schnelli's \cite{Lee-Schnelli2015} work on deformed Wigner matrices, where they tracked cancellations through \emph{two} analogous parameters from the evolutions of spectral edge and (1D) eigenvalue density.
	However, due the non-Hermitian nature of the problem, we have to introduce two extra parameters to get such strong cancellations.
	In fact, what we show with careful adjustments of these four parameters is that, along any $C^1$ curve of deformations starting from $A-z_{0}$ for which the origin remains a sharp edge, the local eigenvalue statistics at the origin remain unchanged. 
	Exploiting the freedom to compare eigenvalue statistics along a fully general curve is a methodological novelty, which allows us to handle the more delicate ``inner" edges of the spectrum which do not appear in Ginibre matrices. See Figure \ref{fig:comp_path} in Section \ref{sec:comp_Gini} for an illustration of how we can use this flexibility.
	 
	The real case has an additional, fundamental difficulty due to the symmetry of the spectrum about the real axis. For a mean-field real random matrix, this symmetry results in a one-parameter family of universality classes around a point $z$, depending on the parameter $\sqrt{N}\absv{\im z}$. Our proof must accommodate this complication so that the comparison between $A+X^{\mathrm{Gin}}$ and $X^{\mathrm{Gin}}$ must be performed along a carefully designed flow $A_{t}-z_{t}$ which keeps the universality class. More precisely, the real case requires to handle a few additional terms that are not present in the complex case. These terms stem from correlations between eigenvalues around $z$ and $\ol{z}$. The key difficulty is to show that these quantities do not affect the evolution of the local statistics along the flow. We divide the proof into two cases: $\absv{\im z}\lesssim N^{-1/2}$ and $\absv{\im z}\gg N^{-1/2}$. In the former case we show that these additional terms have their own cancellations using the construction that $A_{t}$ keeps the exact universality class. In the latter case the extra terms are small by themselves, as one can expect from the decay of correlation of eigenvalues further apart than the natural scale $N^{-1/2}$. In particular we show that the local statistics for the real ensemble at a distance much larger than $N^{-1/2}$ away from the real axis coincide with the statistics for the complex ensemble.
	
	In both the real and complex cases these eigenvalue comparisons are done in parallel with a comparison on the lower bound of the least singular value of $A+X$. Control of the least singular value of non-Hermitian matrices is both essential to applying Girko's Hermitization and an interesting question in its own right. We refer the interested reader to the recent survey \cite{Tikhomirov2023survey} and the references therein. We focus only on works which are most relevant to eigenvalue universality. For Ginibre matrices, \cite{Cipolloni-Erdos-Schroder2020,Cipolloni-Erdos-Schroder2022SIAM} gave an optimal lower tail on the least singular value of $X^{\mathrm{Gin}}-z$ for $z$ at the spectral edge, i.e.\ $\absv{z}=1$. The optimal tail bound in \cite{Cipolloni-Erdos-Schroder2020} for shifted Ginibre matrices was used in the proof of full edge universality for IID matrices in \cite{Cipolloni-Erdos-Schroder2021}. Recently, Shcherbina and Shcherbina \cite{Shcherbina-Shcherbina2022} considered the least singular value tail bound for deformed complex Ginibre matrices $A+X^{\mathrm{Gin}}-z$, however because their result is only for complex matrices and not uniform in the optimal scale of the least singular value we instead use a GFT to the known results of \cite{Cipolloni-Erdos-Schroder2020}.

	For our no-outlier and cluster-rigidity results we first show an improvement of the local law   for $z$ outside $\mathrm{Spec}(A+\bsx)$. Then, for the proof of cluster-rigidity we consider the interpolation $A+\sqrt{t}X$, for which we demonstrate that the number of eigenvalues of $A+\sqrt{t}X$ in each cluster of $A+X$ is fixed over $t$ with very high probability.

\subsection{Notation} In this section we introduce some common notations used throughout the paper. We write $\C_{+}$ for the upper half plane $\C_{+}:=\{z\in\C:\im z>0\}$ and $\bbS^{1}$ for the unit circle $\bbS^{1}\deq\{z\in\C:\absv{z}=1\}$. For a positive integer $k\in\N$ we write $[k]:=\{1,\dots,k\}$. 
We denote the Lebesgue measure on $\C\cong\R^2$ by $\rmd^2 z$, and that on $\C^{k}\cong\R^{2k}$ by $\dd^{2k}\bsz$.
We define the $2N\times 2N$ block constant matrices
\beq\label{eq:def_E1}
	E_{1}\deq\begin{pmatrix} 1 & 0 \\ 0 & 0 \end{pmatrix} ,\qquad E_{2}\deq\begin{pmatrix} 0 & 0 \\ 0 & 1 \end{pmatrix}.
\eeq
For a generic complex matrix $B$, we write its real, imaginary, Hermitian, and skew-Hermitian parts as
\beq\label{eq:def_REIM}
	\Re[B]\deq \frac{B+\ol{B}}{2},\quad 
	\Im[B]\deq\frac{B-\ol{B}}{2\ii},\quad
	\re[B]\deq \frac{B+B\adj}{2}, \quad
	\im[B]\deq \frac{B-B\adj}{2\ii}, 
\eeq
respectively. For a $d\times d$ (in most cases $d=N$ or $2N$) matrix $A$ we use $\brkt{A}=d^{-1}\Tr A$ to denote the normalized trace. We use $\bsu\adj\bsy$ to denote the standard Euclidean inner product of two vectors $\bsu,\bsy\in\C^d$. 

We write $c$ and $C$ for generic positive constants (independent of $N$) whose precise values may vary by lines. For positive quantities $f,g$ we write $f\lesssim g$ and $f\sim g$ if $f\le Cg$ and $cg\le f \le Cg$, respectively. We say that an $N$-dependent event holds with \emph{very high probability} if for any $D>0$, the event holds with probability at least $1-N^{-D}$ for any $N\ge N_0$, where  the threshold $N_0$ depends on $D$. If this event depends on additional parameters,
we say that the event holds with \emph{very high probability uniformly} 
over a given set of parameters if a common threshold $N_0$ can be chosen uniformly for all parameters in this set.
Moreover, for two random variables $X$ and $Y$ we say $X$ is stochastically dominated by $Y$, which we denote by $X\prec Y$, if for every $\xi>0$ the event $\{\absv{X}\leq N^{\xi}\absv{Y} \}$ holds with very high probability.

\subsection*{Acknowledgement}  The authors would like to thank the anonymous referee for providing helpful comments and suggestions. We also thank Joscha Henheik and Volodymyr Riabov for pointing out a gap in an earlier version of the proof of Eq. \eqref{eq:imprllaw}.
	 
\section{Model and Main Results}

We consider $N\times N$ non--Hermitian random matrices with independent, identically distributed (IID) entries satisfying the following assumption. 

\begin{defn}\label{def:matrix assumption}
	An $N$-dependent sequence of complex (resp. real) random matrices $X\equiv X^{(N)}\in\C^{N\times N}$ (resp. $\R^{N\times N}$) is called an \emph{IID random matrix} if it satisfies the following:
	\begin{itemize}
		\item The entries of $X$ are independent.
		\item $\E[X]=0$, $\E[\absv{X_{ij}}^{2}]=1/N$. Additionally, $\E[X_{ij}^{2}]=0$ in the complex case.
		\item The entries of $\sqrt{N}X$ have finite moments, i.e. for each $p\in\N$ there exists a constant $c_{p}>0$ with
		\beq
		\label{eq:momentass}
		\sup_{N\in\N}\max_{i,j}\E\absv{\sqrt{N}X_{ij}}^{p}\leq c_{p}.
		\eeq
	\end{itemize}
	A complex (resp, real) IID matrix with Gaussian entries is referred to as \emph{complex} (resp. \emph{real}) \emph{Ginibre matrix} and denoted by $X^{\mathrm{Gin}(\bbC)}$ (resp. $X^{\mathrm{Gin}(\bbR)}$). 
\end{defn}
\begin{defn}\label{defn:A}
	Let $\bbF$ be either $\R$ or $\C$, $X\in\bbF^{N\times N}$ be an IID matrix, and $A\equiv A^{(N)}\in\bbF^{N\times N}$ be a general sequence of deterministic matrices such that
	\beq\label{eq:Abound}
	\sup_{N\in\N}\norm{A}\leq C_{0}
	\eeq 
	for a constant $C_{0}>0$. We define the $(N\times N)$ \emph{deformed IID matrix} as
	\beq
	A+X\in\bbF^{N\times N}.
	\eeq
\end{defn}

\subsection{The limiting operator $A+\bsx$}
It is expected that the empirical eigenvalue distribution of $A+X$ concentrates around a deterministic ($N$-dependent) measure denoted by $\rho_{A+\bsx}$ for large $N$ (see e.g. \cite{Sniady2002} for $X$ with Gaussian entries).  More precisely, $\rho_{A+\bsx}$ is the \emph{Brown measure} of an infinite-dimensional operator $A+\bsx$. We now introduce a few definitions to properly define the operator $A+\bsx$ and its Brown measure $\rho_{A+\bsx}$. 

\begin{defn}[Brown measure, ]\label{defn:BM}
	Let $(\caM,\brkt{\cdot}_{\caM})$ be a $W\adj$-probability space, i.e., let $\caM$ be a von-Neumann algebra with a faithful, normal, tracial state $\brkt{\cdot}_{\caM}$. For an element $\bsa\in\caM$, we define its \emph{Brown measure} to be the distributional Laplacian
	\beq\label{eq:def_BM}
	\rho_{\bsa}\deq \frac{1}{2\pi}\Delta\brkt{\log\absv{\bsa-\cdot}}_{\caM}.
	\eeq
	Consistently, for a matrix $B\in\bbC^{N\times N}$, we write $\rho_{B}$ for its empirical eigenvalue distribution.
\end{defn}
We refer to \cite[Chapter 11]{Mingo-Speicher2017} for a comprehensive background on Brown measures, in particular
Theorem 5 therein for the existence of the distributional Laplacian in \eqref{eq:def_BM}. As mentioned above, we will approximate the eigenvalue distribution of the large random matrix $A+X\in\C^{N\times N}$ with that of a deterministic, infinite-dimensional operator $A+\bsx$ where $A$ and $\bsx$ live in a common $W\adj$-probability space. Here, the operator $\bsx$ shall be a \emph{circular element} that is \emph{$*$-free} from $A$:
\begin{defn}
	Let $(\caM,\brkt{\cdot}_{\caM})$ be a $W\adj$-probability space.
	\begin{itemize}
		\item A collection $(\caM_{i})_{i\in I}$ of unital subalgebras of $\caM$ is called free if it satisfies
		\begin{equation}
			\brkt{\bsb_{i_{1}}\cdots \bsb_{i_{n}}}_{\caM}=0
		\end{equation}
		whenever $n\in\N$, $i_{j}\in I$, $\bsb_{i_{j}}\in\caM_{i_{j}},\brkt{\bsb_{i_{j}}}_{\caM}=0$, and $i_{j}\neq i_{j+1}$ holds for $j=1,\ldots,n$ (with the convention $i_{n+1}=i_{1}$). A collection of operators $(\bsa_{i})_{i\in I}$ is called $*$-free if the $*$-subalgebras generated by $\{1,\bsa_{i}\}$ are free.
		\item A circular element in $\caM$ is the sum $\bsc=(\bss_{1}+\ii\bss_{2})/\sqrt{2}$ where $\bss_{1},\bss_{2}\in\caM$ are $*$-free semi-circular elements\footnote{A semi-circular element $\bss\in\caM$ is a self-adjoint operator whose spectral distribution under $\brkt{\cdot}_{\caM}$ is the semicircle law $\dd\rho_{\SC}(x)=\frac{1}{2\pi}\sqrt{(4-x^{2})_{+}}\dd x$, i.e. $\brkt{\bss^{k}}_{\caM}=\int_{\R}x^{k}\dd\rho_{\SC}(x)$ for all $k\in\N$.}.
	\end{itemize}	
\end{defn}
We need one more step to properly define $A+\bsx$, that is, to embed the $(N\times N)$ matrix $A$ into the same von-Neumann algebra as $\bsx$. This is achieved by the so-called \emph{meta-argument}, where we lift $\bsx$ to an operator-valued $N\times N$ matrix:
\normalcolor
\begin{defn}\label{defn:vN}
	Let $(\bsx_{ij})_{i,j\in\bbrktt{N}}\in\caM$ be a collection of $*$-free circular elements,
	and define $\bsx=(N^{-1/2} \normalcolor \bsx_{ij})\in\caM^{N\times N}$. We canonically embed $\C^{N\times N}$ in $\caM^{N\times N}$, so that $A+\bsx$ becomes a sum of two elements of $\caM^{N\times N}$. Finally, we consistently extend $\brkt{\cdot}_{\caM}:\caM^{N\times N}\to\C^{N\times N}$ to be the partial trace, i.e.
	\beq\label{eq:vN_pt}
	\brkt{\bsy}_{\caM}=(\brkt{\bsy_{ij}}_{\caM})_{1\leq i,j\leq N}\in\C^{N\times N},\qquad \bsy=(\bsy_{ij})_{1\leq i,j\leq N}\in\caM^{N\times N}.
	\eeq
	Notice that $(\caM^{N\times N}, \brkt{\brkt{\cdot}_{\caM}})$ is a $W\adj$-probability space, where the outer bracket $\brkt{\cdot}$ is the normalized trace on $\C^{N\times N}$.
\end{defn}
\begin{rem}
	With a slight abuse of notations, when we work with (square) matrices of dimension $2N$ instead of $N$, we naturally also denote by $\brkt{\cdot}_{\caM}$ the partial trace $\caM^{2N\times 2N}\to\C^{2N\times 2N}$ applied entrywise as in \eqref{eq:vN_pt}. Likewise, we often canonically embed $\C^{2N\times 2N}$ in $\caM^{2N\times 2N}$. However, we emphasize that we never embed a random matrix in $\caM^{d\times d}$, e.g. we never consider $X$ or $A+X$ as an element of $\caM^{N\times N}$.
\end{rem}
{
	\begin{lem}\label{lem:meta}
		The two operators $A$ and $\bsx$ are $*$-free in $(\caM^{N\times N},\brkt{\brkt{\cdot}_{\caM}})$, and $\bsx$ is a circular element.
	\end{lem}
	With Definitions \ref{defn:BM} and \ref{defn:vN}, we define $\rho_{A+\bsx}$ to be the Brown measure of $A+\bsx$ in $(\caM^{N\times N},\brkt{\brkt{\cdot}_{\caM}})$. We prove Lemma \ref{lem:meta} in Appendix A along with other properties of $\rho_{A+\bsx}$. With these definitions, the Brown measures $\rho_{A+\bsx}$ and $\rho_{A+X}$ are $1/N$
	close to each other in a weak sense, which is a natural 
	extension of the standard circular law for $A=0$ (see e.g. Theorem~\ref{theo:llaw} below with $a=0$).}
\normalcolor

We now recall the classification of the edge behaviors of $\rho_{A+\bsx}$ from \cite{Erdos-Ji2023circ}. 
We say that a point $z_{0}$ on the boundary of the support of $\rho_{A+\bsx}$ belongs to the \emph{regular edge} of $\rho_{A+\bsx}$ if $\lVert (A-z_0)^{-1}\rVert\le C$ for some (large) constant $C>0$.
Regular edge points $z_0\in \C$ are characterized by the 
condition that $ \brkt{|A-z_0|^{-2}}=1$ (see Lemma~\ref{lem:supp} for a proof).
Then \cite[Theorem 2.10]{Erdos-Ji2023circ} (see also Lemma~\ref{lem:sharp}) asserts that $\rho_{A+\bsx}$ has a density\footnote{With a slight abuse of notation, we write $\rho_{A+\bsx}$ for both the Brown measure and its density.} only admitting one of the two possible asymptotic 
behaviors around any regular edge $z_{0}$:
\begin{itemize}
	\item[(i)] (Sharp edge) As $z\to z_0$ we have
	\begin{equation*}
		\rho_{A+\bsx}(z)\sim \mathbbm{1}\big(z\in\supp(\rho_{A+\bsx})\big).
	\end{equation*}
	
	\item [(ii)] (Critical edge) As $z\to z_0$ we have, for some (typically nondegenerate) quadratic form $Q$,
	\[
	\rho_{A+\bsx}(z)\sim Q[z-z_0]\mathbbm{1}\big(z\in\supp(\rho_{A+\bsx})\big).
	\]
\end{itemize}

In the next two subsections we  present our main results.
The first main result concerns the universality of the local statistics of the eigenvalues of $A+X$ (see Section~\ref{sec:univsharp} below) at a sharp edge point. Notice that $\dd\rho_{\bsx}(z)=\pi^{-1}\lone(\absv{z}\leq 1)\dd^{2}z$ by the circular law, hence the sharp edge prevails in $\rho_{\bsx}$. The universality around critical edges is deferred to future work. One of our fundamental inputs to prove universality is an a priori control of the spectrum down to any mesoscopic scale, i.e. on any scale that is slightly 
bigger than the local fluctuations scale of the eigenvalues. 
As a byproduct of this control on mesoscopic scales,  the second main result of this paper proves that with very high probability there are no outliers in the spectrum of $A+X$ at a distance from the edge of the support of $\rho_{A+\bsx}$ which is bigger than the natural local fluctuation scale of the eigenvalues (see Section~\ref{sec:noout} below). 

\subsection{Universality at the sharp edges}
\label{sec:univsharp}

In \cite{Erdos-Ji2023circ}, it is proved that (the density of) $\rho_{A+\bsx}$ has a sharp jump discontinuity around a point $z_{0}\in\C$ (called the \emph{base-point}) under the following assumption:
\begin{assump}\label{assump:edge}
	We assume that $z_{0}\in\C$ satisfies the following:
	\begin{itemize}
		\item[(i)] $\brkt{\absv{A-z_{0}}^{-2}}=1$.
		\item[(ii)] for a constant $C_{1}>0$, $z_{0}\in \caD_{C_{1}}$ where 
		\beq\label{eq:regular}
		\caD_{C}\deq\left\{z\in\C:\Norm{\frac{1}{A-z}}< C,\quad 
		(1+\absv{z})^{3}\Absv{\Brkt{\frac{1}{\absv{A-z}^{4}}(A-z)}}>\frac{1}{C}\right\}.
		\eeq 
	\end{itemize}
	The last assumption depends on the symmetry class. Let $C_{2}>C_{1}$ be another constant. When $\bbF=\bbC$, we assume that
	\begin{itemize}
		\item[(iii$\C$)] $z_{0}$ is in the unique\footnote{ The fact that there exists a unique unbounded component in $\mathcal{D}_{C_2}$ follows from the fact $\mathcal{D}_{C_2}^c$ is contained in a large enough disk, as a consequence of $\lVert A\rVert\le C_0$ uniformly in $N$. The same reasoning applies to the real case.} unbounded component of $\caD_{C_{2}}$.
	\end{itemize}
	When $\bbF=\bbR$, we fix a small constant $\frc>0$ and assume the following:
	\begin{itemize}
		\item[(iii$\R$)] If $|\im z_{0}|>N^{-1/2+\frc}$, $z_{0}$ is in one of the two unbounded components of $\caD_{C_{2}}\setminus\R$.
		\item[(iii$\R'$)] If $\absv{\im z_{0}}<N^{-1/2+\frc}$, $\re z_{0}$ is in the unique unbounded component of $\caD_{C_{2}}\cap\R$.
	\end{itemize}
\end{assump}
The equality $\brkt{\absv{A-z_{0}}^{-2}}=1$ guarantees that the base point $z_{0}$ is at the spectral edge of $A+X$. The first inequality in \eqref{eq:regular} ensures that the associated Dyson equation (see Lemma \ref{lem:MDE} below) is stable around $z_{0}$, and the second excludes the critical edges. Assumptions (iii) guarantee that there exists a $C^1$ path from $z_0$ to $\infty$ that purely consists of sharp edge points of $A+\sqrt{t}X$ for some $t>0$. In the real case the path has an additional constraint from the symmetry class, hence we divide (iii$\R$) and (iii$\R'$).

Our first main theorem asserts that the local eigenvalue statistics of $A+X$ around a base point $z_0$ at a sharp edge
is universal in the $N\to\infty$ limit (after proper rescaling). In particular, we show that the corresponding $k$--point correlation functions coincide with those of Ginibre ensembles in the same symmetry class. 
To be precise, for a given $z\in\bbS^{1}$, our reference functions $p_{z,k}^{\mathrm{Gin}(\bbF)}$ are defined by the limiting $k$--point correlation functions of eigenvalues of $\sqrt{N}(z^{-1}X^{\mathrm{Gin}(\bbF)}-1)$ near the origin.  These are well-known functions; in the complex case they are independent of $z\in\bbS^{1}$ and given by (see \cite{Ginibre1965})
\begin{equation}
	\label{eq:gincorrf}
	\begin{split}
		p^{\mathrm{Gin}(\C)}_{z,k}(w_1,\dots,w_k)&=\det\left[K^{\mathrm{Gin}(\C)}(w_i,w_j) \right]_{1\leq i,j\leq k},\\
		K^{\mathrm{Gin}(\C)}(w_1,w_2)&=\frac{1}{2\pi}\left[1+\mathrm{erf}\left(-\sqrt{2}(\overline{{w_2}}+w_1)\right) \right]\e{-\frac{\absv{w_1}^2}{2}-\frac{\absv{w_2}^2}{2}+w_1\overline{{w_2}}}.
	\end{split}
\end{equation}
Explicit formulas for $p^{\mathrm{Gin}(\R)}_{z,k}$ are known 
also in the real case \cite{Borodin-Sinclair2009}, but they
are much more involved as a consequence of the special role played by the real axis.

The proper rescaling is governed by a single scaling parameter defined as
\beq\label{eq:def_gamma0}
\gamma_{0}\deq -I_{4}^{-1/2}I_{3},
\eeq
where we set
\beq\label{eq:def_I3I4_0}
I_{3}\deq\Brkt{\frac{1}{\absv{A-z_{0}}^{4}}(A-z_{0})\adj},\qquad\quad I_{4}\deq\Brkt{\frac{1}{\absv{A-z_{0}}^{4}}}.
\eeq
The parameter $\gamma_0\in\C$ is chosen so that the density of the Brown measure of $\gamma_{0}(A+\bsx-z_{0})$ matches that of $(\bsx-1)$ near their sharp edge point at the origin
(see Lemma~\ref{lem:edge} below for more details).  In order to state our main universality result for any $k\in\N$ we implicitly define the \emph{$k$-point correlation functions} $p_k^{(N)}$ of the eigenvalues $\{\sigma_i\}_{i\in [N]}$ of $A+X$ by
\begin{equation}
	\label{eq:corrf}
	\left(\begin{matrix}
		N \\
		k
	\end{matrix}\right)^{-1}\E\sum_{i_1\ne \dots \ne i_k} F(\sigma_{i_1},\dots, \sigma_{i_k})=\int_{\C^k} F(\bsw) p_k^{(N)}(\bsw)\, \dd^{2k} \bsw.
\end{equation}
for any smooth compactly supported test function $F:\C^k\to \C$. The sum in the left--hand side of \eqref{eq:corrf} is over indices all distinct from each other.

\begin{thm}[Sharp edge universality]
	\label{thm:main correlations}
	Let $A+X\in\bbF^{N\times N}$ be the deformed IID matrix defined in Definition \ref{defn:A}, and let the base point
	$z_0\in \C$ satisfy 
	Assumption~\ref{assump:edge}. Fix any $k\in\N$, and let $p_{k}^{(N)}$ be the $k$-point correlation function of the eigenvalues of $\sqrt{N}\gamma_{0}(A+X-z_0)$ near the origin. Then, for any compactly supported smooth function $F:\C^{k}\rightarrow\C$ we have
	\begin{equation}\label{eq:main_corr}
		\int_{\C^{k}} F(\bsw)\left[p_{k}^{(N)}(\bsw)-p^{\mathrm{Gin}(\mathbb{F})}_{z_{1},k}(\bsw) \right]\dd^{2k}\bsw=O\left(N^{-c}\right),
	\end{equation}
	for some small $c>0$, where $z_{1}\in\bbS^{1}$ is chosen\footnote{This choice of $z_{1}$ is relevant only when $\bbF=\R$, since $p_{z,k}^{\mathrm{Gin}(\C)}$ is independent of $z\in\bbS^{1}$.  Furthermore, in the real case, we can choose one of the two $z_1$ satisfying $\im z_{1}=(\absv{\gamma_{0}}\wedge 1/\absv{\im z_{0}})\im z_{0}$, since the local statistics around these two points are the same.} such that $\im z_{1}=(\absv{\gamma_{0}}\wedge 1/\absv{\im z_{0}})\im z_{0}$.
	The implicit constant and $c>0$ in \eqref{eq:main_corr} depend only on $k,F$, and the constants from Definitions \ref{def:matrix assumption} -- \ref{defn:A} and Assumption \ref{assump:edge}.
\end{thm}
The same edge universality for $A+X$ was proven in \cite{Liu-Zhang2024} for the very special case where $X$ is a Ginibre matrix
and $A$ is a normal matrix with finitely supported empirical spectral measure and for $A=0$ and IID $X$ in \cite{Cipolloni-Erdos-Schroder2020}.

While Theorem \ref{thm:main correlations} is formulated for deterministic $A$, it is clear that we may assume that $A$ is a random matrix, independent of $X$, properly chosen so that it satisfies Assumption~\ref{assump:edge} with sufficiently high probability. This simple idea allows us to conclude
edge universality for many other natural random matrix ensembles as a direct consequence of 
Theorem~\ref{thm:main correlations}; we give some examples below. We point out that in all these examples 
it is fundamental that Theorem~\ref{thm:main correlations} allows for very general $A$'s, which in particular may not be normal matrices. We remark that the normality of $A$ was  a fundamental assumption in the previous works for special $A+X^{\mathrm{Gin}}$-type matrices in \cite{Liu-Zhang2023, Liu-Zhang2024}.

\subsubsection{Examples}
\label{sec:examples}

We now present a few examples of relevant random matrix ensembles for which 
edge universality can be proven  by a simple application of our general Theorem~\ref{thm:main correlations} and our GFT arguments. 

To illustrate the idea, we first 
consider an $N\times N$ random matrix $X$ with independent centered entries $X_{ij}$  with a general \emph{flat}
variance profile such that 
\begin{equation}
	\label{eq:defvarprof}
	\E |X_{ij}|^2 =:  \frac{S_{ij}}{N}, \qquad\quad \frac{c}{N}\le S_{ij}\le \frac{C}{N},
\end{equation}
for two $N$--independent constants $c,C>0$. By \cite[Theorem 2.5]{Alt-Erdos-Kruger2018} it follows that the limiting eigenvalue density of such $X$ is supported on a disk of radius $\sqrt{\rho(S)}$, with $\rho(S)$ denoting the spectral radius of the variance matrix $S$. In the first step,  by a GFT with two moments matching we can prove 
that the edge statistics of $X$ match with that of
\begin{equation}
	\label{eq:decvarpro}
	\widetilde{X}+\epsilon X^{\mathrm{Gin}}, \qquad \epsilon\sim 1,
\end{equation}
where $\widetilde{X}$ is Gaussian matrix\footnote{The Gaussianity of $\widetilde{X}$ is not important, but the easiest
	way to construct  \eqref{eq:decvarpro} is first to match $X$ with a Gaussian matrix with the same
	variances $S_{ij}$ and then decompose this Gaussian matrix
	following $S=[S-(\epsilon^2/N)\lone\lone\tp]+ (\epsilon^2/N)\lone\lone\tp$.}
independent of $X^{\mathrm{Gin}}$ and has independent entries with the flat variance
profile $S-(\epsilon^2/N)\lone\lone\tp$. We point out that GFT with 
two moments matching has two versions depending on whether 
the single resolvent local law is available or not. Traditionally, GFT arguments use local laws as inputs, 
in particular this was the case in the first edge universality proof for iid matrices in~\cite{Cipolloni-Erdos-Schroder2021}
and we follow the same path in Appendix~C, to show that the local statistics close to the edges of 
$X$ and  $\widetilde{X}+\epsilon X^{\mathrm{Gin}}$ are asymptotically the same. This proof directly works 
for the  situation of \eqref{eq:decvarpro} with a variance profile as well. Note that for this ensemble
the local law could be borrowed directly from \cite{Alt-Erdos-Kruger2021}. 
However, our Zig-Zag method also proves a stronger version of GFT without {\it a priori} knowledge of the local law.
In fact, a trivial extension of 
our  proof presented in Section~\ref{sec:Proof of the local law} (see e.g. the proof of Theorem~\ref{theo:llaw})
also yields  a self-contained proof of a local law for \eqref{eq:decvarpro}. This observation is important
in situations where local law {\it a priori} is not yet known, see e.g. the correlated case discussed below.

In the second step we prove edge universality for~\eqref{eq:decvarpro}.
We will condition on $\widetilde{X}$ and, after a trivial rescaling, we can
use our Theorem~\ref{thm:main correlations} 
with $A=\epsilon^{-1}\widetilde{X}$ and $X=X^{\mathrm{Gin}}$.  To check the  conditions
of Theorem~\ref{thm:main correlations}  we now argue that
whenever $\wh{z}_0$
is on the spectral edge of $\widetilde{X}+\epsilon X^{\mathrm{Gin}}$, i.e. 
$|\wh{z}_{0}|=\sqrt{\rho(S)}$, the pair $(A,z_{0})\deq(\epsilon^{-1}\widetilde{X},\epsilon^{-1}\wh{z}_{0})$ satisfies  Assumption~\ref{assump:edge} with very high probability. To see this, we only need to recall two facts: 
\begin{itemize}
	\item[(i)] By \cite[Theorem 2.1]{Erdos-Kruger-Schroder2019}, for any $z\in\C$, on a macroscopic scale
	the empirical eigenvalue density $\rho_{\absv{\wt{X}-z}}$ concentrates around a deterministic probability 
	measure $\mu_{z}$ on $[0,\infty)$ with rate $O(N^{-1})$ and by \cite[Corollary 2.3]{Erdos-Kruger-Schroder2019}
	there is no outlying eigenvalue of $\absv{\wt{X}-z}$ beyond a distance $N^{-\delta}$, $\delta>0$, from the support of 
	$\mu_{z}$. 
	\item[(ii)] By \cite[Proposition 3.2(iii)]{Alt-Erdos-Kruger2018} and from the easy fact that 
	$\rho(\wt{S})<\rho(S)-c\epsilon^2$, where $\wt {S}$ is the variance
	profile of $\wt{X}$, we conclude that the solution of the MDE for  the Hermitization of $\widetilde{X}-\wh{z}_0$
	vanishes at zero, moreover it is of order $O(\eta)$ along imaginary axis. Using the implication $(vi) \Longrightarrow (v)$
	of \cite[Lemma D.1]{Alt-Erdos-Kruger2020}, we thus conclude that
	$\mu_{\wh{z}_{0}}$ is supported away from zero,  i.e.,
	$\dist(0,\supp\mu_{\wh{z}_{0}})\sim1$. 
\end{itemize}
These two facts together imply that $\absv{\brkt{\absv{\epsilon^{-1}(\wt{X}-\wh{z}_{0})}^{-2}}-1}\prec N^{-1}$ and $\norm{(\wt{X}-\wh{z}_{0})^{-1}}=O(1)$ with very high probability. Likewise, applying \cite[Section 3]{Alt-Erdos-Kruger2018} to the MDE associated to $\wt{X}$ proves that $\{\absv{z}\geq\epsilon^{-1}\sqrt{\rho(S)}\}$ is contained in the domain $\caD_{C}$ from Assumption \ref{assump:edge} with very high probability, confirming\footnote{Assumption \ref{assump:edge} (i) holds with a random error of size $\prec N^{-1}$, which does not affect the local statistics.} Assumption~\ref{assump:edge} for $(A,z_{0})=(\epsilon^{-1}\wt{X},\epsilon^{-1}\wh{z}_{0})$ with very high probability.
Thus, conditioning on $\wt{X}$ and working in the probability space of  $X^{\mathrm{Gin}}$, we can directly apply Theorem~\ref{thm:main correlations} to show that the local statistics of $\widetilde{X}+\epsilon X^{\mathrm{Gin}}$
around $\wh{z}_{0}$ are universal and that they coincide with the corresponding real/complex Ginibre statistics.
Combining it with the first step above we proved edge universality for the ensemble~\eqref{eq:defvarprof}.

In fact, modulo some technicalities (see footnotes), a similar argument 
is expected to work for a much larger class of ensembles, in particular for matrices 
with correlated entries\footnote{\label{fn:corr_extra} When $X$ has correlated entries we expect that 
	the stronger self-contained version of GFT with two moments matching still works showing that $X$ with general non-Gaussian correlated entries has the same edge statistics as a Gaussian matrix with the same correlations.
	The latter then can be decomposed as in \eqref{eq:decvarpro} and  the second step of the argument above would complete the proof.}
(e.g. for the non--Hermitian counterpart of \cite[Assumptions (A)--(F)]{Erdos-Kruger-Schroder2020}, as considered in 
\cite{Alt-Kruger2021}), and even more generally for matrices of the form $A+X$, with $X$ correlated 
and $A$ being a deterministic 
matrix\footnote{\label{fn:A+flat_extra}In this case, after the decomposition 
	$A+\wt{X} +\epsilon X^{\mathrm{Gin}}$ as in \eqref{eq:decvarpro}, an extra complication is that the MDE for $A+\wt{X}$ is not fully explored yet, 
	i.e. we cannot rely on \cite{Alt-Erdos-Kruger2018} as we did for $A=0$ above. The MDE for such an $A+\wt{X}$, in the 
	special case of $\wt{X}$ as in \eqref{eq:defvarprof}, was studied only recently in \cite{Alt-Kruger2024}, restricted to diagonal $A$; see also Remark~\ref{rem:compak} for a related discussion.}
satisfying an assumption similar to Assumption \ref{assump:edge}.
We omit the details for brevity. We conclude this section pointing out that in particular this would imply edge universality for the so called elliptic ensemble, i.e. when $\E X_{ij}X_{ji}=\tau/N$, with $\tau\in [0,1)$ (see \cite{Lee-Riser2016} for the special case when the $X_{ij}$ are Gaussians).

\subsection{No outliers and cluster-rigidity for $A+X$}
\label{sec:noout}
Consider a deterministic matrix $A$ satisfying the following assumption.
\begin{assump}
	\label{ass:newassout}
	There exists a constant $C_{3}>0$ such that for all $z\in\partial\supp\rho_{A+\bsx}$ we have
	\begin{equation}
		\label{eq:opnormb}
		\left\lVert \frac{1}{A-z}\right\rVert\le C_{3}.
	\end{equation}
	We also assume a non-degeneracy condition\footnote{ The non-degeneracy condition is introduced to simplify the proof; in the general case we would need higher order expansions depending on the singularity types of the degenerate edges in \cite{Alt-Kruger2024}. }, that for all $z\in\partial\supp\rho_{A+\bsx}$
	\begin{equation}\label{eq:nondegenerate} 
		\absv{\nabla_{z}\brkt{\absv{A-z}^{-2}}}+\absv{\det \nabla_{z}^{2}\brkt{\absv{A-z}^{-2}}}>\frac{1}{C_{3}},
	\end{equation}
	where we consider the map $z\mapsto \brkt{\absv{A-z}^{-2}}$ as that on $\C\cong\R^{2}$ and $\nabla_{z}$ is the gradient on $\R^{2}$.
\end{assump}

\begin{rem}
	\label{rem:nooutlocal}
	Unlike in Assumption~\ref{assump:edge}, we formulated a global assumption on $A$ rather than only at a single point, as we want to exclude the presence of outliers in the whole complement of the support of $\rho_{A+\bsx}$; note that now we allow critical edges. However, we point out that we assumed the bounds \eqref{eq:opnormb}--\eqref{eq:nondegenerate} for all points on the edge of the support of $\rho_{A+\bsx}$ only because in Theorem~\ref{theo:noput} below we exclude outliers globally in the complement of the spectrum. If we are interested in excluding outliers only locally around an edge point $z_0$, then it is enough to assume that \eqref{eq:opnormb}--\eqref{eq:nondegenerate} holds for $z=z_0$.
\end{rem}

We now define the typical fluctuation scale of the eigenvalues of $A+X$:
\begin{defn}[Fluctuation scale]
	\label{def:flucscale}
	Fix a point $z\in \C$, we say that $\sigma_{\mathrm{f}}(z)=\sigma_{\mathrm{f}}(z,N)$ is the fluctuation scale of the eigenvalues of $A+X$ around the point $z$ if
	\begin{equation}
		\int_{|w-z|\le \sigma_{\mathrm{f}}(z)} \rho_{A+\bsx}(w)\, \rmd^2 w=\frac{1}{N}.
	\end{equation}
\end{defn}
Note that for $z$'s in the bulk of the spectrum or close to sharp edges we have $\sigma_{\mathrm{f}}(z,N)\sim N^{-1/2}$, while for $z$'s close to a critical edge $\sigma_{\mathrm{f}}(z,N)\sim N^{-1/4}$;  see Lemma \ref{lem:rmf} for the size of $\sigma_{\rmf}(z,N)$ for a general edge point $z$.  Fix any $\epsilon>0$, and define the deterministic set
\begin{equation}
	\label{eq:kappneigh}
	\mathrm{Spec}_\epsilon(A+\bsx):=\bigcup_{z\in \supp(\rho_{A+\bsx})} D\left(z,N^\epsilon\sigma_{\mathrm{f}}(z)\right).
\end{equation}
Here $D(z,r)$ denotes the open disk of radius $r$ centered at $z$.

\begin{thm}[No outliers]
	\label{theo:noput}
	Let $A+X$ be a deformed IID matrix as in Definition~\ref{defn:A}, and suppose that $A$ satisfies Assumption~\ref{ass:newassout}. Then, for any small $\epsilon>0$, we have
	\begin{equation}
		\label{eq:noout}
		\mathrm{Spec}(A+X)\subset  \mathrm{Spec}_\epsilon(A+\bsx),
	\end{equation}
	with very high probability uniformly in $A$ and $X$  satisfying ~\eqref{eq:momentass}, \eqref{eq:Abound} and \eqref{eq:opnormb}.
\end{thm}

We also prove that the number of eigenvalues in each component (which we often also refer to as ``cluster") of the support of $\rho_{A+\bsx}$ is deterministic. 
\begin{cor}[Cluster--rigidity]
	\label{cor:clustrig}
	Under the assumptions of Theorem~\ref{theo:noput}, the following holds. Fix any small $\epsilon,\delta>0$, and let $\mathcal{C}$ be a union of (one or several) components of $\mathrm{Spec}_\epsilon(A+\bsx)$ such that $\mathrm{dist}(\mathcal{C}, \mathrm{Spec}_\epsilon(A+\bsx)\setminus \mathcal{C})\ge\delta$. Then, counting the eigenvalues with multiplicity, we have
	\begin{equation}
		\big|\mathrm{Spec}(A+X)\cap \mathcal{C}\big|=\big|\mathrm{Spec}(A)\cap \mathcal{C}\big|,
	\end{equation}
	with very high probability uniformly over $A$ and $X$ satisfying ~\eqref{eq:momentass}, \eqref{eq:Abound} and \eqref{eq:opnormb}.
	In particular, the number of eigenvalues in $\mathcal{C}$ is a deterministic integer with very high probability.
\end{cor}

Both of these results are  consequences of the following very general local law close to the edge of the spectrum of $A+X$. Fix $a\in [0,1/2]$, consider a smooth compactly supported function $g:\C \to \C$, and define
\[
g_{z_0,a}(z):=g\big(N^a(z-z_0)\big).
\]
Denote the eigenvalues of $A+X$ by $\sigma_1,\dots,\sigma_N$.
\begin{thm}[Local law at regular edges]
	\label{theo:llaw}
	Let $A+X$ be a deformed IID matrix as in Definition~\ref{defn:A}. Fix any $C_{4}>0$, and fix $z_0$ on the edge of the support of $\rho_{A+\bsx}$ such that $\lVert (A-z_0)^{-1}\rVert\le C_{4}$. Then, for any $\xi>0$ and $a\in[0,1/2)$, we have
	\begin{equation}
		\label{eq:llaw}
		\left|\frac{1}{N}\sum_{i=1}^N g_{z_0,a}(\sigma_i)-\int_\C g_{z_0,a}(z)\rho_{A+\bsx}(z)\,\rmd^2 z\right|\le \frac{N^\xi}{N^{1-2a}},
	\end{equation}
	with very high probability uniformly in $A,X,$ and $z_0$ satisfying ~\eqref{eq:momentass}, \eqref{eq:Abound}, and $\norm{(A-z_{0})^{-1}}\leq C_{4}$.
\end{thm}
Note that the error term in \eqref{eq:llaw} is smaller than leading deterministic term only on mesoscopic scales around the base point $z_0$, i.e. when $N^{-a}\gg \sigma_{\mathrm{f}}(z_0)$. In particular, even when $z_{0}$ is a critical edge, \eqref{eq:llaw} remains effective for all mesoscopic scales $a\in[0,1/4)$. We also point out that to prove \eqref{eq:llaw} we do not need Assumption~\ref{ass:newassout} for every point at the edge of the support of $\rho_{A+\bsx}$, but it is enough to assume that \eqref{eq:opnormb} holds for the base point $z_0$. The same result for $z_0$ in the bulk of the spectrum of $A+X$ follows by the averaged local law in \cite[Theorem 2.6]{Cipolloni-Erdos-Henheik-Schroder2023arXiv}. If $z_0$ is in an intermediate regime (between bulk and edge) or is outside of the spectrum we can prove a similar local law, but we do not present it here for brevity.

\begin{rem}
	Using a strong self-contained two moments matching GFT as described below \eqref{eq:decvarpro} we can conclude that Theorem~\ref{theo:noput}, Corollary~\ref{cor:clustrig}, and \eqref{eq:llaw} hold for matrices $X$ with
	independent entries and  variance profile as in \eqref{eq:defvarprof}.
	In particular, 
	\begin{equation}
		\label{eq:specradas}
		\big|\rho(X)-\sqrt{\rho(S)}\big|\le \frac{N^\xi}{\sqrt{N}},
	\end{equation}
	with very high probability where $\rho(X)$ denotes the spectral radius. We point out that the bound \eqref{eq:specradas} was already proven in \cite[Theorem 2.1]{Alt-Erdos-Kruger2021}. 
	More generally, when $X$ is centered and has some correlations as in \cite{Alt-Kruger2021},\footnote{In these models, by \cite[Theorem 2.5]{Alt-Kruger2021}, the self-consistent spectrum consists of one connected component, the limiting density is rotationally symmetric, and there are no critical edges.} we expect that a strong self-contained GFT with two moments matching 
	similar to the one described below \eqref{eq:decvarpro} can prove the bound \eqref{eq:specradas}, but with $\rho(S)$ being replaced with $\rho(\mathscr{S})$, where $\mathscr{S}$ is the covariance operator of $X$.
\end{rem}

\begin{rem}
	\label{rem:compak}
	Simultaneously with the current paper, 
	an independent work by Alt and Kr\"uger \cite{Alt-Kruger2024} considered the $\epsilon$-pseudospectrum 
	of $A+X$ for
	some  $N$-independent $\epsilon\sim 1$, where $X$ has a flat variance profile\footnote{After we have made them 
		aware of the current  work, the authors of \cite{Alt-Kruger2024} informed us that they can relax the flatness assumption.}
	as in \eqref{eq:defvarprof} and $A$ is a diagonal matrix. Their main result \cite[Theorem 2.4]{Alt-Kruger2024} roughly reads: for each fixed $\epsilon,\delta>0$, it holds
	in the large $N$ limit 
	that
	\begin{equation}
		\label{eq:jtres}
		\Psi\mathrm{Spec}_{\epsilon}(A+X)\subset \Psi\mathrm{Spec}_{\epsilon}(A+\bsx)\subset\Psi\mathrm{Spec}_{\epsilon+\delta}(A+X),
	\end{equation}
	where $\bsx$ is a suitable deterministic counterpart of $X$ and $\Psi\mathrm{Spec}_{\epsilon}$ stands for the $\epsilon$-pseudospectrum\footnote{The $\epsilon$--pseudospectrum of any matrix  $X$ is defined by $\Psi\mathrm{Spec}_{\epsilon}(X)\deq \{z\in\C:\norm{(X-z)^{-1}}\geq\epsilon^{-1}\}$.}. In particular, their result implies that $A+X$ has no outliers on macroscopic scale. 
	Notice that the subscript $\epsilon$ in $\mathrm{Spec}_{\epsilon}$ from \eqref{eq:kappneigh} and in 
	$\Psi\mathrm{Spec}_{\epsilon}$ from \eqref{eq:jtres} plays very different roles.
	
	In the case of constant variance profile, i.e. for $X_{\mathrm{iid}}+A$, when $X_{\mathrm{iid}}$ is an IID matrix and $A$ is a general deformation satisfying Assumption~\ref{ass:newassout}, 
	we immediately get a much stronger result: we obtain both inclusions in \eqref{eq:jtres} with the order one $\epsilon,\delta$ replaced with the optimal \emph{local fluctuation scale for the smallest singular value},
	see~\eqref{etaf} (up to an $N^\epsilon$ factor).
	This follows from our optimal local laws both inside the spectrum~\eqref{eq:local law goal} and 
	outside~\eqref{eq:imprllaw}. Moreover a similar two-sided inclusion with optimal precision also 
	holds for the spectrum  instead of the pseudospectrum in the following form:
	$$
	\mathrm{Spec}(A+X)\subset  \mathrm{Spec}_\epsilon(A+\bsx) \subset  \bigcup_{z\in \mathrm{Spec}(A+X)} 
	D\left(z,N^{\epsilon+\delta}\sigma_{\mathrm{f}}(z)\right).
	$$
	More precisely, the first inclusion is just \eqref{eq:noout}, while the second inclusion follows from the fact that there many eigenvalues close to the edge of the support of $\rho_{A+\bsx}$, as a consequence of the local law \eqref{eq:llaw}.
	
	In fact, our method can prove the optimal no-outlier theorem for $A+X$ with $X$ having a flat
	variance profile using the same decomposition as in Section \ref{sec:examples}. More precisely, we consider 
	$A+X$ as a deformed IID matrix, i.e. work with $A+\wt{X}+\epsilon X^{\mathrm{Gin}}$ instead\footnote{This can again be achieved using a self-contained two moments matching GFT as described below \eqref{eq:decvarpro}.}, and then
	condition on $\wt{X}$. Then we need to guarantee that Assumption \ref{ass:newassout}  holds
	with $A$ replaced by $A+\wt{X}$ with very high probability, given the same pair of inputs as in Section \ref{sec:examples} for  $A+\wt{X}$ instead of $\wt{X}$. Namely, we need
	to guarantee  (i) the macroscopic convergence  and no-outlier theorem for $\rho_{|A+\wt{X}-z|}$ and (ii) $\inf\supp\rho_{|A+\wt{\bsx}-z|}\sim 1$ for $z\in\partial\supp\rho_{A+\bsx}$, where $\wt{\bsx}$ and $\bsx$ are the operator analogues of $\wt{X}$ and $X$. The first input (i) already follows from \cite[Theorem 2.1 and Corollary 2.3]{Erdos-Kruger-Schroder2019} since they allow 
	general deformations. 
	The second input (ii), which is a purely deterministic result, needs to be checked by solving the corresponding MDE
	on the macroscopic scale.
	For example, this was explicitly done in 
	\cite[Proposition 5.15 (iv)]{Alt-Kruger2024} for the special case where $A$ is diagonal.
\end{rem}

\subsection{Organization of the proof} In Section \ref{sec:Preliminaries} we gather the gather the free probabilistic and random matrix tools which we use in the proofs of our main results. This includes the statement of the local law (see Theorem \ref{thm:local law: Local Law}). We then present the proof of Theorem \ref{theo:noput}, Corollary \ref{cor:clustrig}, and Theorem \ref{theo:llaw}. In Section~\ref{sec:proofmain} we prove Theorem \ref{thm:main correlations} using our main technical comparison result Theorem \ref{thm:comp_Gini}, which compares the joint moments of log determinants of various non-Hermitian random matrices. In Sections~\ref{sec:comp_Gini}--\ref{sec:short GFT} we provide the proof of Theorem \ref{thm:comp_Gini}. These three sections divide the proof up into the cases: Section \ref{sec:comp_Gini} comparing $A+X^{\mathrm{Gin}(\C)}$ and $X^{\mathrm{Gin}(\C)}$, Section \ref{sec:real} comparing $A+X^{\mathrm{Gin}(\R)}$ and $X^{\mathrm{Gin}(\R)}$, and finally Section~\ref{sec:short GFT} comparing $A+X$ and $A+X^{\mathrm{Gin}(\mathbb{F})}$ for both fields $\mathbb{F}=\R,\C$. Finally, in Section~\ref{sec:Proof of the local law}, we present the proof of the local law Theorem \ref{thm:local law: Local Law}. Appendix~\ref{sec:meta} compiles technical results on the macroscopic eigenvalue density.

\section{Preliminaries and the proofs of Theorems~\ref{theo:noput}--\ref{theo:llaw}}\label{sec:Preliminaries}

In this section we present the main technical ingredients.  In Section~\ref{sec:llawgen}, we list some properties of the Brown measure $\rho_{A+\bsx}$ and we state the local law close to the cusp regime of the  Hermitization of $A+X-z$, when $z$ is a regular edge point. Then, in Section~\ref{sec:proofout} we prove Theorems~\ref{theo:noput}--\ref{theo:llaw}.

\subsection{Preliminary results and the local law for $A+X$}
\label{sec:llawgen}

Before we state the local law for $A+X$ (see Theorem \ref{thm:local law: Local Law} below), we collect  some preliminary results on free probability and matrix Dyson equations.

\begin{lem}\label{lem:edge}
	The Brown measure $\rho_{A+\bsx}$ has a density that satisfies the following:
	\begin{itemize}
		\item[(i)] $\rho_{A+\bsx}(z)\leq 1/\pi$ for all $z\in\C$.
		\item[(ii)]$\supp\rho_{A+\bsx}=\ol{\{z\in\C:\brkt{\absv{A-z}^{-2}}>1\}}.$
		\item[(iii)] If $z_{0}$ satisfies Assumption \ref{assump:edge} (i)--(ii), the density $\rho_{A+\bsx}(z_{0}+w)$ follows the following asymptotics as $w\to 0$:
		\beq\label{eq:sharp}\begin{aligned}
			\rho_{A+\bsx}(z_{0}+w)=&\frac{1}{\pi}\absv{\gamma_{0}}^{2}+O(\absv{w}),\qquad  &\text{for }z_{0}+w\in \supp\rho_{A+\bsx},\\
			\rho_{A+\bsx}(z_{0}+w)\equiv& 0,\qquad &\text{for }z_{0}+w\notin\supp\rho_{A+\bsx},
		\end{aligned}\eeq
		where $\gamma_{0}$ was defined in \eqref{eq:def_gamma0}.
		The implicit constant in \eqref{eq:sharp} depends only on the $C_{0}$ from Definition \ref{defn:A} and $C_{1}$ from Assumption \ref{assump:edge} (ii).
	\end{itemize}
\end{lem}
\begin{proof}
	The first two statements of Lemma \ref{lem:edge} are due to \cite[Theorem 7.10]{Belinschi-Yin-Zhong2024}. The last one follows from Lemma~\ref{lem:sharp} and the fact that
	\begin{equation}
		{\absv{\gamma_{0}}^{2}=I_{4}^{-1}\absv{I_{3}}^{2}=\brkt{\absv{A-z_{0}}^{-4}}^{-1}\absv{\brkt{(A-z_{0})\absv{A-z_{0}}^{-4}}}^{2}.}
	\end{equation}
\end{proof}

We next connect the Brown measure $\rho_{A+\bsx}$ to the Hermitization of the operator $A+\bsx$ via the trace of the Hermitized resolvent, which in turn is characterized by a matrix Dyson equation:
\begin{lem}\label{lem:MDE}
	For each $z\in\C$ and $\wh{z}\in\C_{+}$, the matrix Dyson equation (MDE)
	\beq\label{eq:MDE0}
	-\frac{1}{M}=-\begin{pmatrix}
		0 & A-z \\ (A-z)\adj & 0
	\end{pmatrix}+\wh{z}+\brkt{M}
	\eeq
	has a unique $2N\times 2N$ solution $M\equiv M_{A-z}(\wh{z})$ with $\im M=\frac{1}{2\ii}(M-M\adj)>0$. 
	Then we have
	\beq\label{eq:free_Stiel}
	\Brkt{\begin{pmatrix}
			-\wh{z} & A+\bsx-z \\ (A+\bsx-z)\adj & -\wh{z}
		\end{pmatrix}^{-1}}_{\caM}=M_{A-z}(\wh{z})\in\C^{2N\times 2N}, \qquad \forall z\in\C,\,\,\wh{z}\in\C_{+},
	\eeq
	so that the Brown measure of the sum $A+\bsx$ can be computed as the distributional Laplacian
	\beq\begin{split}
		\rho_{A+\bsx}(z)\equiv-\frac{1}{2\pi}\Delta_{z}\int_{0}^{\infty}\left[\Brkt{\Brkt{\frac{\eta}{\absv{A+\bsx-z}^{2}+\eta^{2}}}_{\caM}}-\frac{\eta}{1+\eta^{2}}\right]\dd\eta	\\
		=-\frac{1}{2\pi}\Delta_{z}\int_{0}^{\infty}\left[\rho^{A-z}(\ii\eta)-\frac{\eta}{1+\eta^{2}}\right]\dd\eta,
	\end{split}\eeq
	where we defined the spectral density of the Hermitization of $A+\bsx-z$ as
	\begin{equation}\label{eq:rho_def}
		\rho^{A-z}(w)\deq \frac{1}{\pi}\im\brkt{M_{A-z}(w)},\qquad w\in\C_{+}.
	\end{equation}
\end{lem} 
The existence and uniqueness of solution of the MDE are consequences of \cite[Theorem 2.1]{Helton-Rashidi-Speicher2007}, 
and \eqref{eq:free_Stiel} follows from combining \cite[Chapter 9, Theorems 2 and 11]{Mingo-Speicher2017}. \normalcolor We remark that $\rho^{A-z}(\cdot+\ii 0)$ is exactly the spectral density of the Hermitization of $A+\bsx-z$. With a slight abuse of notations we write $\supp\rho^{A-z}$ for the support of this measure on $\R$.

Note, for $\wh{z}=\ii\eta$ on the imaginary axis, $\brkt{M(\ii\eta)}$ is purely (positive) imaginary and $M$ is the resolvent of the Hermitization of $(A-z)$ evaluated at $\ii\eta+\brkt{M}$. Hence for all $\eta>0$ we have
\begin{equation}
	\label{eq:M boundedness}
	\Norm{M(\ii\eta)}\leq\norm{(A-z)^{-1}}.
\end{equation}

A key tool in the proof of Theorems~\ref{thm:main correlations}, \ref{theo:noput}--\ref{theo:llaw} is an optimal local law for the Hermitization $H^z$ defined in \eqref{eq:herm}. We define the resolvent of $H^z$ by
\begin{equation}
	G^{A,z}(\wh{z})=G^z(\wh{z})=\left(H^{z}-\wh{z} \right)^{-1},\quad \wh{z}\in\C_{+}, 
\end{equation} although for our local law we will exclusively consider $\wh{z}=\ii \eta$ for $\eta\geq 0$ on the positive imaginary axis. The dependence of $G^{A,z}$ through $H^{z}$ defined in \eqref{eq:herm} will often be omitted when $A$ is to be thought of as fixed.

The local law gives an optimal comparison between $G^{z}$ and the deterministic approximation $M_{A-z}$ near an edge point $z_0$. The proof of this local law will be presented in Section~\ref{sec:Proof of the local law}.   

\begin{thm}[Local law for $H^{z_0}$ near the edge]\label{thm:local law: Local Law}
	Let $A+X$ be a deformed IID matrix as in Definition \ref{defn:A}. Let $z_{0}$ be an edge point of $A+X$ such that $\Norm{(A-z_{0})^{-1}}<C_{4}$ for some $C_{4}>0$. There exists a (small enough) constant $c>0$ such that for any deterministic vectors $\bsu,\bsy\in\C^{2N}$ and matrix $B\in\C^{2N\times 2N}$ 
	\begin{equation}\label{eq:local law goal}\begin{aligned}
			\absv{ \brkt{ B(G^{z}(\ii\eta)-M_{A-z}(\ii\eta)) } }&\prec \norm{B}\frac{1}{N\eta}, \\ \left|\bsu^*\left(G^{z}(\ii\eta)-M_{A-z}(\ii\eta)\right)\bsy \right|&\prec \Norm{\bsu }\Norm{\bsy}\left[\sqrt{\frac{\rho^{A-z}(\ii\eta)}{N\eta}}+\frac{1}{N\eta}\right],
	\end{aligned}\end{equation}
	uniformly in $A,X,z_{0}$ satisfying ~\eqref{eq:momentass}, \eqref{eq:Abound}, $\norm{(A-z_{0})^{-1}}\leq C_{4}$,  $|z-z_{0}|<c$  and $N^{-1}\leq\eta\leq N^{100}$.
\end{thm}

For $z\in\C$, we define the function $f_{A-z}:\C\setminus\mathrm{Spec}(A-z)\rightarrow (0,\infty)$ as \begin{equation}
	\label{eq:deffaz}
	f_{A-z}(w):=\Brkt{\frac{1}{\absv{A-z-w}^2 } }.
\end{equation}
\begin{lem}\label{lem:MDE_asymp}
	Let $A\in\C^{N\times N}$ and $z\in\partial\supp\rho_{A+\bsx}$ satisfy $\norm{A}\leq C_{0}$ and $\norm{(A-z)^{-1}}\leq C_{1}$ for some constants $C_{0},C_{1}>0$.
	As $\eta,\absv{w}\to 0$, we have the following asymptotics:
	\beq\label{eq:M_size}
	\im\brkt{M_{A-(z+w)}(\ii\eta)}\sim\begin{cases}
		\absv{f_{A-z}(w)-1}^{1/2}+\eta^{1/3}, & z+w\in\supp\rho_{A+\bsx},\\
		\dfrac{\eta}{\absv{f_{A-z}(w)-1}+\eta^{2/3}}, & z+w\notin\supp\rho_{A+\bsx}.
	\end{cases}
	\eeq
	Additionally, if $z\in\caD_{C_{1}}$ is a sharp edge (recall $\caD_{C}$ from \eqref{eq:regular}), then the regularized eigenvalue density of $H^{z}$, $v\equiv v(w,\eta)\deq\rho^{A-z-w}(\ii\eta)+\eta$, satisfies the following asymptotic cubic equation as $\absv{w},\eta\to0$:
	\beq
	\label{eq:cubic}
	I_{4}v^{3}-2\re[wI_{3}]v-\eta=O((\absv{w}^{1/2}+\eta^{1/3})^{5}),
	\eeq
	where $I_{3}$ and $I_{4}$ are the same quantities as in \eqref{eq:def_I3I4_0} but with $z_{0}$ replaced by $z$. The implicit constants in \eqref{eq:M_size} and \eqref{eq:cubic} depend only on $C_{0}$ and $C_{1}$.
\end{lem}
\begin{proof}
	Plugging in $\bsa=A-z$ and $z_{0}=0$ in \cite[Lemma 5.1]{Erdos-Ji2023circ}, we obtain
		\begin{equation}\label{eq:M_size1}
			v=\eta+\im\brkt{M_{A,z+w}(\ii\eta)}\sim\begin{cases}
				\absv{f_{A-z}(w)-1}^{1/2}+\eta^{1/3}, & f_{A-z}(w)>1,\\
				\dfrac{\eta}{\absv{f_{A-z}(w)-1}+\eta^{2/3}}, & f_{A-z}(w)\leq 1,
			\end{cases}
		\end{equation}
		where the implicit constants are uniform over $C_{0},C_{1}$ due to \cite[Remark 5.2]{Erdos-Ji2023circ}.
		Using Lemma~\ref{lem:supp} (ii) we can convert the cases in \eqref{eq:M_size1} to the ones in \eqref{eq:M_size}, and then subtracting $\eta$ from both sides proves \eqref{eq:M_size}.
	
	Now we turn to the proof of \eqref{eq:cubic}  With $v=\eta+\im\brkt{M_{A,z+w}(\ii\eta)}$, we may write the (trace of) the MDE \eqref{eq:MDE0} as
	\beq
	v-\eta=\Brkt{\frac{v}{\absv{A-z-w}^{2}+v^{2}}}.
	\eeq
	Then we expand the right-hand side with respect to $v$ and then $w$, so that 
	\beq\begin{split}
		v-\eta&=v\Brkt{\frac{1}{\absv{A-z-w}^{2}}}-v^{3}\Brkt{\frac{1}{\absv{A-z-w}^{4}}}+O(v^{5})	\\
		&=v\Brkt{\frac{1}{\absv{A-z}^{2}}}+2\re[wI_{3}]v-I_{4}v^{3}+O(v^{5}+v\absv{w}^{2}+v^{3}\absv{w})\\
		&=v+2\re[wI_{3}]v-I_{4}v^{3}+O((\absv{w}^{1/2}+\eta^{1/3})^{5}),
	\end{split}\eeq
	where in the last line we also used \eqref{eq:M_size}. Rearranging the terms proves the  approximate  cubic equation  \eqref{eq:cubic}.
\end{proof}

\subsection{Proof of Theorems~\ref{theo:noput}--\ref{theo:llaw}}
\label{sec:proofout}

We first present the proof of Theorem~\ref{theo:llaw}, and then we conclude this section with the proof of Theorem~\ref{theo:noput} and Corollary~\ref{cor:clustrig}. All estimates in this section are uniform over $X,A$, and $z_{0}$ satisfying \eqref{eq:momentass}, \eqref{eq:Abound}, and $\norm{(A-z_{0})^{-1}}\leq C_{4}$.

\begin{proof}[Proof of Theorem~\ref{theo:llaw}]
	
	Given the local law \eqref{eq:local law goal} for the Hermitized resolvent, the proof of \eqref{eq:llaw} is by now standard (see e.g. \cite[Section 5.2]{Alt-Erdos-Kruger2018}). Here we sketch the main steps for completeness.
	
	This proof relies on Girko's formula \eqref{eq:girkof} in the form
	\begin{equation}
		\label{eq:girnew}
		\int_{\C}F(z)\dd\rho_{A+X}(z)=-\frac{1}{4\pi N}\int_\C \Delta F(z)\int_{N^{-100}}^{N^{100}} \langle \Im G^z(\ii\eta) \rangle\,\rmd \eta\rmd^2 z +O(N^{-10}),
	\end{equation}
	where we used the bound for the smallest singular value in \cite[Theorem 3.2]{Tao-Vu2010MC} to remove the very small $\eta\le N^{-100}$ regime and simple explicit calculations to remove the regime $\eta\ge N^{100}$ (see e.g. \cite[Eq. (5.27)]{Alt-Erdos-Kruger2018}). Then using the monotonicity of $\eta\to\eta\langle\Im G^z(\ii\eta) \rangle$ to extend \eqref{eq:local law goal} down to $\eta\ge N^{-100}$, and plugging \eqref{eq:local law goal} into \eqref{eq:girnew} we immediately obtain the bound \eqref{eq:llaw}.
	
\end{proof}
Before moving on to the proof of Theorem~\ref{theo:noput}, we present a lemma that shows the value of the function $f_{A}(z)$ is away from $1$ when $z$ is away from $\supp\rho_{A+\bsx}$. We present its proof in Appendix~\ref{sec:meta} since it is purely computational.
\begin{lem}\label{lem:fc}
	If $A$ satisfies Assumption \ref{ass:newassout}, then there is a constant $c_{*}$ such that we have
	\begin{equation}
		1-f_{A}(z)\geq N^{-1/2+c_{*}\epsilon},\qquad \forall z\in\C\setminus\mathrm{Spec}_{\epsilon}(A+\bsx).
	\end{equation}
\end{lem}
\begin{proof}[Proof of Theorem~\ref{theo:noput}]
	To prove that the are no outliers in the spectrum of $A+X$ we need an improvement of the averaged local law \eqref{eq:local law goal} when $z_0$ is outside of the spectrum (see the proof of \cite[Theorem 2.1]{Alt-Erdos-Kruger2021} for a similar analysis). In particular, for any deterministic $B\in\C^{2N\times 2N}$, we claim that
	\begin{equation}
		\label{eq:imprllaw}
		\big|\langle (G^z(\ii\eta)-M_{A-z}(\ii\eta))B\rangle\big|\prec \frac{1}{N\kappa},
	\end{equation}
	with $\kappa:=\mathrm{dist}(\ii\eta,\supp(\rho^z))$; note that if $0$ is in $\supp(\rho^z)$, then $\kappa=\eta$. If $z$ is well outside of the support of 
	$\rho_{A+\bsx}$, i.e. more than the typical fluctuation scale, then we can choose even $\eta=0^+$. The proof of \eqref{eq:imprllaw} is presented at the end of Section~\ref{sec:zig}.
	
	Proceeding similarly to the proof of \cite[Theorem 2.5 (ii)]{Alt-Erdos-Kruger2018} we can show that there are no eigenvalues in the region
	\[
	\big\{z\in \C:\mathrm{dist}(z,\supp(\rho_{A+\bsx}))\ge c\big\},
	\]
	with very high probability for some small $N$--independent constant $c>0$. We now strengthen this statement by reducing $c$ to $N^\epsilon\sigma_{\mathrm{f}}$. More precisely, we show that if we pick  $z\notin \mathrm{Spec}_\epsilon(A+\bsx)$ then $z\notin \mathrm{Spec}(A+X)$. Fix an edge point $z_0$, then for $z\notin \mathrm{Spec}_\epsilon(A+\bsx)$ such that $|z-z_0|\le c$, we define (cf. \cite[Eq. (2.7)]{Erdos-Kruger-Schroder2020}) the Hermitian fluctuation scale
	\begin{equation}\label{etaf}
		\eta_{\mathrm{f}}=\eta_{\mathrm{f}}(z):=\frac{|f_{A-z_0}(z)-1|^{1/6}}{N^{2/3}},
	\end{equation}
	with $f_{A-z_0}$ being defined in \eqref{eq:deffaz}. Fix $\eta:= N^{c_*\epsilon/5}\eta_{\mathrm{f}}$, then, by \eqref{eq:M_size} it follows that
	\begin{equation}
		\label{eq:boundM}
		\rho^{A-z}(\ii\eta)
		\le \frac{N^{-c_*\epsilon/5}}{N\eta},
	\end{equation}
	where we used Lemma \ref{lem:fc} to deduce $f_{A}(z)\leq 1-N^{-1/2+c_{*}\epsilon}$.  
	
	Next, using the local law \eqref{eq:imprllaw}, together with \eqref{eq:boundM}, we have
	\begin{equation}
		\label{eq:boundG}
		\sup\Big\{ \langle \Im G^z(\ii\eta)\rangle:\,|z-z_0|\le c, \, |f_{A-z_0}(z)-1|\ge N^{-1/2+c_*\epsilon}\Big\}\prec \frac{N^{-c_*\epsilon/5}}{N\eta}.
	\end{equation}
	Here we used a standard grid argument and the Lipschitz continuity of $z\mapsto  \langle \Im G^z\rangle$ to show that \eqref{eq:imprllaw} holds simultaneously for all $|z-z_0|\le c$.  More precisely, we consider a mesh of $N^{10}$ equidistant points of the set $\{z:|z-z_0|\le c\}$. Then, for any point of this mesh \eqref{eq:imprllaw} holds with very high probability. On the other hand if $z'\in\{z:|z-z_0|\le c\}$ but it is not one of the points in the mesh, then we use that
	\[
	\big|\langle \Im G^{z'}(\ii\eta)\rangle-\langle \Im G^{z''}(\ii\eta)\rangle\big|\lesssim N^2 |z'-z''|,
	\]
	where $z''$ is the closest point to $z'$ which is also part of the mesh, and we used that $\lVert \nabla_z \Im G^z\rVert\le 1/\eta^2\le N^2$ deterministically. This shows that \eqref{eq:imprllaw} in fact holds simultaneously for all points in $\{z:|z-z_0|\le c\}$. Throughout the paper whenever we refer to a \emph{grid argument} we refer to an argument similar to the one described above, even if we do not spell out the details. 
	
	Denote the eigenvalues of $H^z$ by
	\[
	\lambda_{-N}^z\le \dots\le \lambda_{-1}^z\le 0\le \lambda_1^z\le \dots \le \lambda_N^z,
	\]
	then \eqref{eq:boundG} implies that as long as $z$ is such that $|f_{A-z_0}(z)-1|\ge N^{-1/2+c_*\epsilon}$ we have $\lambda_1^z\ne 0$  with very high probability (recall that $\lambda_{-i}^z=-\lambda_i^z$). In fact, if $\lambda_1^z=0$ then we would have $ \Im G^z(\ii\eta)\rangle\ge 1/(2N\eta)$, which contradicts \eqref{eq:boundG}. This concludes the proof showing that $z\notin \mathrm{Spec}(A+X)$, as a consequence of the fact that
	\begin{equation}
		\label{eq:hermnonherm}
		0\in\mathrm{Spec}(H^z)\Longleftrightarrow z \in \mathrm{Spec}(A+X).
	\end{equation}
	
\end{proof}

We conclude this section with the proof of cluster--rigidity.

\begin{proof}[Proof of Corollary~\ref{cor:clustrig}]
	This proof is similar to the proof of \cite[Corollary 2.9]{Alt-Erdos-Kruger-Schroder2020}. We thus only highlight the main differences and present a sketch of the proof for the convenience of the reader.
	
	As a consequence of \eqref{eq:hermnonherm}, similarly to the proof of Theorem~\ref{theo:noput} above, it is enough to show that the Hermitization $H^z$ satisfies ``band rigidity", in the sense of \cite[Corollary 2.9]{Alt-Erdos-Kruger-Schroder2020}, when its density has a gap around zero of order one size.
	
	Consider a contour $\Gamma$ which encircles $\mathcal{C}$ and it is such that $\mathrm{dist}(\Gamma,\mathrm{Spec}_\epsilon(A+\bsx))\ge \delta/10$. 
	For any 
	$z\in\Gamma$ 
	we let $H^z$ to be the Hermitization of $A+X-z$. We now consider a flow that interpolates between $H^z$ and $\mathcal{A}^z$, the Hermitization of $A-z$, which is a deterministic matrix, i.e.\ define $H_t^z$ to be the Hermitization of $A+X_t-z$, where
	\begin{equation}
		X_t:=\sqrt{1-t}X.
	\end{equation}
	Note that $H_0^z=H^z$, $H_1^z=\mathcal{A}^z$, and $H_t^z=\sqrt{1-t}W+\mathcal{A}^z$, with $W$ being the Hermitization of $X$. We now show that for any $D>0$ we have
	\begin{equation}
		\label{eq:bandrigidity}
		\P\big(0\in \mathrm{Spec}(H_t^z)\, \mathrm{for}\,\, \mathrm{some}\,\, (t,z)\in [0,1]\times\Gamma\big)\le N^{-D}.
	\end{equation}
	As a consequence of \eqref{eq:hermnonherm}, this yields the desired result \eqref{eq:boundG}.
	
	To prove \eqref{eq:bandrigidity}, we distinguish the analysis into two cases: $t$ close to one and the rest. By assumption, $0\notin \mathrm{Spec}(\mathcal{A}^z)$ for any $z\in\Gamma$, additionally $\lVert H^0\rVert \le C_*$, for some fixed $C_*>0$, we can thus choose a sufficiently small $c_*>0$ such that $0\notin \mathrm{Spec}(H_t^z)$, for any $z\in\Gamma$, for $t\in [1-c_*,1]$ with very high probability. In the regime when $t\in [0,1-c_*)$, the variance of the entries of $X_t$ is of order $1/N$, hence, modulo a simple rescaling, $X_t$ is an IID matrix in the sense of Definition~\ref{def:matrix assumption}. Hence, we can use the same argument as in \eqref{eq:boundM}--\eqref{eq:boundG} to conclude that, with very high probability, $0\notin \mathrm{Spec}(H_t^z)$ for any $z\in\Gamma$ and $t\in [0,1-c_*]$. This proves \eqref{eq:bandrigidity}.
\end{proof}

\section{Proof of Theorem~\ref{thm:main correlations}}\label{sec:proofmain}

\subsection{Proof of Theorem~\ref{thm:main correlations}}
The purpose of this section is to reduce the proof of Theorem \ref{thm:main correlations} to our main technical result, Theorem \ref{thm:comp_Gini} below, which compares joint moments of log determinants and lower tail of the smallest singular value. All estimates (including asymptotic relations and stochastic dominance) in this section are uniform over $A$, $X$, and $z_{0}$ satisfying \eqref{eq:momentass}, \eqref{eq:Abound}, and Assumption \ref{assump:edge}. We first reduce the test function $F$, i.e. we claim that Theorem \ref{thm:main correlations} follows from the following result:
\begin{thm}\label{thm:main1}
	Let $z_0\in\C$ satisfy Assumption \ref{assump:edge}, $\gamma_{0}$ be defined in \eqref{eq:def_gamma0}, and $k\in\N$. Let $F_{1},\ldots, F_{k}\in C_{c}^{\infty}(\C)$, and define their local rescalings
	\beqs
	\wt{F}_{j}\deq NF_{j}(\sqrt{N}\gamma_{0}(\cdot-z_{0})),\qquad \wt{F}_{j}^{\mathrm{Gin}}\deq NF_{j}(\sqrt{N}z_{1}^{-1}(\cdot-z_{1})),
	\eeqs
	where $z_{1}$ is chosen as in Theorem \ref{thm:main correlations}. Then, for a constant $c>0$ we have
	\beq\label{eq:main}
	\E\prod_{j=1}^{k}\left(\int_{\C}\wt{F}_{j}(w)\dd (\rho_{A+X}-\rho_{A+\bsx})(w)\right) 
	-\E\prod_{j=1}^{k}\left(\int_{\C}\wt{F}_{j}^{\mathrm{Gin}}(w)\dd (\rho_{X^{\mathrm{Gin}(\bbF)}}-\rho_{\bsx})(w)\right) =O(N^{-c}),
	\eeq
	where $c>0$ and the implicit constant depend\footnote{They depend on $F_{j}$'s only via $\norm{F_{j}}_{L^{\infty}}, \norm{\Delta F_{j}}_{L^{1}}$, and $\mathrm{diam}(\supp F_{j})$.} on $k$ and $F_{j}$'s.
\end{thm} 
\begin{proof}[Proof of Theorem \ref{thm:main correlations}]
First, we prove \eqref{eq:main_corr} when $F$ is of the product form $F(\bsw)=\prod_{j=1}^{k}F_{j}(w_{j})$. Given Theorem \ref{thm:main1}, it suffices to check that
\beq\label{eq:density bound in N}
\Absv{\int_{\C}\wt{F}_{j}(w)\dd\rho_{A+\bsx}(w)-\int_{\C}\wt{F}^{\mathrm{Gin}}_{j}(w)\dd\rho_{\bsx}(w)}=O(N^{-1/2}).
\eeq
To see \eqref{eq:density bound in N}, we write $\rho$ for the density of $\rho_{A+\bsx}$ and
\beq
\int_{\C}\wt{F}_{j}(w)\dd\rho_{A+\bsx}(w)
=\absv{\gamma_{0}}^{-2}\int_{\C}F_{j}(w)\rho(z_{0}+\gamma_{0}^{-1}N^{-1/2}w)\dd^{2}w.
\eeq
On the other hand, Lemma \ref{lem:edge} (ii) gives
\beq
\absv{\gamma_{0}}^{-2}\rho(z_{0}+\gamma_{0}^{-1}N^{-1/2}w)=\begin{cases}
	\frac{1}{\pi}+O(N^{-1/2}\absv{w}) & \text{if }z_{0}+\gamma_{0}^{-1}N^{-1/2}w\in\supp\rho, \\
	0 & \text{elsewhere}.
\end{cases}
\eeq
Note from Lemma \ref{lem:edge} (ii) that the support of $\rho_{A+\bsx}$ is a super-level set of the map $z\mapsto\brkt{\absv{A-z}^{-2}}$, and recall that $\gamma_{0}$ is exactly the gradient (written as a complex number) of this function at $z_{0}$. Thus we find that $\partial\supp\rho_{\gamma_{0}(A+\bsx-z_{0})}$ is tangent to the half-space $\{w:\re w\leq 0\}$,  so that
\beq
\int_{w\in\supp F_{j}}\absv{\lone(z_{0}+\gamma_{0}^{-1}N^{-1/2}w\in \supp\rho)-\lone(\re w\leq 0)}\dd^{2}w=O(N^{-1/2}).
\eeq
We thus conclude that
\beq
\Absv{\int_{\C} \wt{F}_{j}(w)\dd\rho_{A+\bsx}(w)-\frac{1}{\pi}\int_{\C}F_{j}(w)\lone(\re w\leq 0)\dd^{2}w}=O(N^{-1/2}).
\eeq
Repeating the same argument for $\rho_{\bsx}$ proves \eqref{eq:density bound in N}.

Now it only remains to extend the test function beyond those in the product form $F(\bsw)=\prod_{j=1}^{k}F_{j}(w_{j})$. To this end, we approximate the given smooth test function $F$ in Theorem \ref{thm:main correlations} by a linear combination of those in the product form using truncated Fourier series. We omit further details and refer the reader to \cite[Theorem 1]{Cipolloni-Erdos-Schroder2021} or \cite[Corollary 7]{Tao-Vu2015}. This concludes the proof of Theorem \ref{thm:main correlations}.
\end{proof}

We next present the key technical inputs, and use them to prove Theorem \ref{thm:main1}. The main inputs are two parallel comparison results, for the two pairs of regularized functionals\footnote{The parameter $z_{1}\in\bbS^{1}$ in \eqref{eq:def_L01} and \eqref{eq:sing_01} do not play any role when $\bbF=\C$ since a complex Ginibre matrix is rotationally invariant. Nonetheless, we keep the current form to facilitate the proof in the real case in Section \ref{sec:real}.}
\beq\begin{aligned}\label{eq:def_L01}
L_{0}(w)\deq \Tr{\log(\absv{A+X-z_{0}-\gamma_{0}^{-1}N^{-1/2}w}^{2}+\eta_{0}^{2})}-\Tr\brkt{\log(\absv{A+\bsx-z_{0}-\gamma_{0}^{-1}N^{-1/2}w}^{2}+\eta_{0}^{2})}_{\caM},	\\
L_{1}^{\mathrm{Gin}(\bbF)}(w)\deq \Tr{\log(\absv{X^{\mathrm{Gin}(\bbF)}-z_{1}-z_{1}N^{-1/2}w}^{2}+\eta_{1}^{2})}-\Tr\brkt{\log(\absv{\bsx-z_{1}-z_{1}N^{-1/2}w}^{2}+\eta_{1}^{2})}_{\caM},
\end{aligned}\eeq
and 
\beq\begin{aligned}\label{eq:sing_01}
N_{0}(w)\deq \frac{1}{c_{0}}\Brkt{\frac{\eta_{0}}{\absv{A+X-z_{0}-\gamma_{0}^{-1}N^{-1/2}w}^{2}+\eta_{0}^{2}}},\\
N_{1}^{\mathrm{Gin}(\bbF)}(w)\deq \Brkt{\frac{\eta_{1}}{\absv{X^{\mathrm{Gin}(\bbF)}-z_{1}-z_{1}N^{-1/2}w}^{2}+\eta_{1}^{2}}}.
\end{aligned}\eeq
Here, the scales $\eta_{0}$ and $\eta_{1}$ are defined by (recall the definition of $I_{4}$ from \eqref{eq:def_I3I4_0})
\beq\label{eq:def_eta}
\eta_{0}\deq c_{0}^{-1}\eta_{1},\qquad c_{0}\deq I_{4}^{-1/4},\qquad \eta_{1}\deq N^{-3/4-\delta},
\eeq
where $0<\delta<(1\wedge \frc)/100$ is a small but fixed constant, with\footnote{The $\frc$-dependent constraint for $\delta$ is not needed in the complex case, but we include it here to avoid repetition.} $\frc$ from Assumption \ref{assump:edge} (iii$\R$) and (iii$\R'$). With our choice of the scale $\eta_{1}$, we define the fundamental control parameter
\beq\label{eq:def_Psi}
\Psi\deq \frac{1}{N\eta_{1}}=N^{-1/4+\delta}.
\eeq
Our error estimates are typically of the form $O(N^{\epsilon}\Psi^{n})$ for some $n\in\N$ and arbitrarily small (but fixed) $\epsilon>0$. Now we are ready to state our comparison results:
\begin{thm}\label{thm:comp_Gini}
Let $k\in\N$ and $C>0$ be fixed. Then there exists a constant $c>0$ such that the following hold uniformly in $w\in\C$ and $\bsw=(w_{1},\ldots,w_{k})\in\C^{k}$ with $\absv{w},\norm{\bsw}\leq C$:
\begin{align}
	\E\prod_{i=1}^{k}L_{0}(w_{i})=&\E\prod_{i=1}^{k}L_{1}^{\mathrm{Gin}(\bbF)}(w_{i})+h^{(0)}(\bsw)+O(N^{-c}),\label{eq:comp_logd}\\
	\E N_{0}(w)
	=&\E N_{1}^{\mathrm{Gin}(\bbF)}+O(N^{-c}\Psi), \label{eq:comp_sing}
\end{align}
where $h^{(0)}\equiv h_{N}^{(0)}\in C^{2k}(\C^{k})$ is harmonic in at least one variable, i.e. it satisfies $(\prod_{j=1}^{k}\Delta_{w_{j}})h^{(0)}\equiv0$. 
\end{thm}
\begin{proof}[Proof of Theorem \ref{thm:main1}]
Using Girko's formula, we express the two statistics in terms of the log potential:
\beq\begin{split}
	\int_{\C}\wt{F}_{j}(w)\dd\rho_{A+X}(w)=&N\int_{\C}F_{j}(w)\dd\rho_{N^{1/2}\gamma_{0}(A+X-z_{0})}(w)	\\
	=&-\frac{N}{2\pi}\int_{\C}\Delta_{w}F_{j}(w)\brkt{\log\absv{N^{1/2}\gamma_{0}(A+X-z_{0})-w}}\dd^{2}w \\
	=&-\frac{N}{2\pi}\int_{\C}\Delta_{w}F_{j}(w)\brkt{\log \absv{A+X-z_{0}-\gamma_{0}^{-1}N^{-1/2}w}}\dd^{2}w,
\end{split}\eeq
where in the second line we used integration by parts and in the third line we subtracted $\log(N^{1/2}\absv{\gamma_{0}})$ from the integrand using harmonicity. Writing the same expression for $X^{\mathrm{Gin}(\bbF)}$ as well as for the deterministic counterparts, we decompose
\beq\begin{split}\label{eq:I1I2}
	\int_{\C}\wt{F}_{j}(w)\dd(\rho_{A+X}-\rho_{A+\bsx})(w)&=I^{(j)}_{1}+I^{(j)}_{2},	\\
	\int_{\C}\wt{F}_{j}^{\mathrm{Gin}}(w)\dd(\rho_{X^{\mathrm{Gin}(\bbF)}}-\rho_{\bsx})(w)&=I^{(j)\mathrm{Gin}(\bbF)}_{1}+I^{(j)\mathrm{Gin}(\bbF)}_{2},
\end{split}\eeq
where
\beq
\begin{split}
	I^{(j)}_{1}\deq&-\frac{1}{2\pi}\int_{\C}\Delta_{w}F_{j}(w) L_{0}(w)\dd^{2}w,\\
	I^{(j)}_{2}\deq&-\frac{1}{2\pi} \int_{\C}\Delta_{w}F_{j}(w)	\\
	&\qquad \times \int_{0}^{\eta_{0}}N\Brkt{\frac{\eta}{\absv{A+X-z_{0}-\gamma_{0}^{-1}w}^{2}+\eta^{2}}-\im M_{A-z_{0},\gamma_{0}^{-1}w}(\ii\eta)}\dd\eta\dd^{2}w,
\end{split}\eeq
and $I^{(j)\mathrm{Gin}(\bbF)}_{1}$ and $I^{(j)\mathrm{Gin}(\bbF)}_{2}$ are defined similarly with $(A,X,z_{0},\eta_0,\gamma_{0})$ replaced by $(0,X^{\mathrm{Gin}(\bbF)},z_{1},\eta_1,z_{1}^{-1})$.

The next lemma shows that the product of $I_{1}^{(j)}$'s has the leading contribution to the statistics. 
\begin{lem}\label{lem:cut_off}
	For some constant $c>0$ depending on $k$ and $F_{j}$'s, we have
	\beq\begin{split}\label{eq:cut_off}
		\E\prod_{j=1}^{k}\left(\int_{\C}\wt{F}_{j}(w)\dd (\rho_{A+X}-\rho_{A+\bsx})(w)\right)=&\E\prod_{j=1}^{k} I_{1}^{(j)}+O(N^{-c}),	\\
		\E\prod_{j=1}^{k}\left(\int_{\C}\wt{F}_{j}^{\mathrm{Gin}}(w)\dd (\rho_{X^{\mathrm{Gin}(\bbF)}}-\rho_{\bsx})(w)\right)=&\E\prod_{j=1}^{k} I_{1}^{(j)\mathrm{Gin}(\bbF)}+O(N^{-c}).
	\end{split}\eeq
\end{lem}

Given Lemma \ref{lem:cut_off} and \eqref{eq:comp_logd}, we conclude the proof of Theorem \ref{thm:main1} as 
\beq\begin{split}
	&\E\prod_{j=1}^{k}\left(\int_{\C}\wt{F}_{j}(w)\dd (\rho_{A+X}-\rho_{A+\bsx})(w)\right) 
	-\E\prod_{j=1}^{k}\left(\int_{\C}\wt{F}_{j}^{\mathrm{Gin}}(w)\dd (\rho_{X^{\mathrm{Gin}(\bbF)}}-\rho_{\bsx})(w)\right)	\\
	=&\Absv{\E\prod_{j=1}^{k} I_{1}^{(j)}-\E\prod_{j=1}^{k} I_{1}^{(j)\mathrm{Gin}(\bbF)}}+O(N^{-c})	\\
	\leq &\int_{\C^{k}}\prod_{j=1}^{k}\absv{\Delta F_{j}(w_{j})}\Absv{\E\prod_{j=1}^{k}L_{0}(w_{j})-\prod_{j=1}^{k}L_{1}^{\mathrm{Gin}(\bbF)}(w_{j})}\dd^{2k}\bsw	+O(N^{-c})\\
	\lesssim &N^{-c},
\end{split}\eeq
where in the last line we used that $F_{j}\in C_{c}^{\infty}(\C)$.
\end{proof}

\begin{proof}[Proof of Lemma \ref{lem:cut_off}]
We only prove the first estimate in \eqref{eq:cut_off}, which implies the second by taking $A=0$ and $X=X^{\mathrm{Gin}(\bbF)}$.

We first collect the inputs. As a standard corollary of local laws (see e.g. \cite[Lemma 3]{Cipolloni-Erdos-Schroder2021}\footnote{
	In \cite{Cipolloni-Erdos-Schroder2021} there is the additional assumption which we do not impose that there exist $\alpha,\beta>0$ such that the entries of $X$ have density $g$ satisfying \begin{equation}
		g\in L^{1+\alpha}(\C),\qquad \Norm{g}_{1+\alpha}\leq N^{\beta}.
	\end{equation} However, as is pointed out in \cite[Remark 1]{Cipolloni-Erdos-Schroder2021} this assumption is to simplify presentation, and could be removed using \cite[Theorem 3.2]{Tao-Vu2010MC}. Alternatively we could assume $X$ has a very small Gaussian component and use \cite[Section 6.1]{Tao-Vu2015} as was done in \cite[Section 6]{Cipolloni-Erdos-Xu2023}. 
}), we have
\beq\label{eq:cut_off_ll}
\absv{I_{1}^{(j)}}+\absv{I_{2}^{(j)}}\prec \norm{\Delta F_{j}}_{L^{1}}\lesssim 1.
\eeq
Also, we can deduce from \eqref{eq:comp_sing} and \cite[Eq.(28)]{Cipolloni-Erdos-Schroder2021} that 
\beq\begin{split}\label{eq:sing_tail_Gini}
	\P[\lambda_{1}(A+X-z_{0}-\gamma_{0}^{-1}N^{-1/2}w)\leq \eta_{0}]\leq& 2c_{0}N\eta_{0}N_{0}(w)	\\
	=& 2N\eta_{1}N_{1}^{\mathrm{Gin}(\bbF)}+O(N^{-c})	\\
	\lesssim &N^{-c},
\end{split}\eeq
where $\lambda_{1}(\cdot)$ denotes the smallest singular value. Then \eqref{eq:sing_tail_Gini} in turn implies, following the same proof as \cite[Lemma 4]{Cipolloni-Erdos-Schroder2021} and taking small enough $\delta,c>0$, that (recall $\delta$ from \eqref{eq:def_eta})
\beq\label{eq:cut_off_I2}
\E I_{2}^{(j)}=O(N^{-c}).
\eeq
Finally, we have the non-random, non-asymptotic inequalities
\beq\label{eq:abs_0}
\absv{I_{1}^{(j)}}\leq \frac{N}{\eta_{0}^{2}}\Norm{(\norm{A}^{2}+\norm{X}^{2}+N^{-1}\absv{w}^{2})\Delta F_{j}(w)}_{L^{1}(\dd^2 w)},\qquad \absv{I_{2}^{(j)}}\leq N\norm{F_{j}}_{L^{\infty}}+\absv{I_{1}^{(j)}},
\eeq
which in turn implies
\beq\label{eq:abs}
\absv{I_{1}^{(j)}}+\absv{I_{2}^{(j)}}\lesssim N^{10}(1+\norm{X}^{2}).
\eeq
The first inequality in \eqref{eq:abs_0} follows from the fact that, for any $Y\in\C^{N\times N}$, $f\in C_{c}^{2}(\C)$, and $\eta>0$,
\beq\begin{split}
	\Absv{\int_{\C}\Delta f(w)\Brkt{\log(\absv{Y}^{2}+\eta^{2})}\dd^{2}w}
	=\Absv{\int_{\C}\Delta f(w)\Brkt{\log\left(1+\frac{\absv{Y}^{2}}{\eta^{2}}\right)}\dd^{2}w}	\\
	\leq \int_{\C}\absv{\Delta f(w)}\left(1+\frac{\Brkt{\absv{Y}^{2}}}{\eta^{2}}\right)\dd^{2}w,
\end{split}\eeq
and the second inequality in \eqref{eq:abs_0} is due to \eqref{eq:I1I2}.

Computing the difference of both sides of \eqref{eq:cut_off} and using \eqref{eq:I1I2}, it suffices to prove that 
\beq\label{eq:cut_off_prf}
\E\prod_{j=1}^{k}I_{i_{j}}^{(j)}=O(N^{-c/2}),\qquad i_{j}\in\{1,2\},
\eeq
whenever $\absv{\{j:i_{j}=2\}}\geq1$. We define the event $\Xi$ on which
\beq
\max_{1\leq j\leq k}(\absv{I_{1}^{(j)}}+\absv{I_{2}^{(j)}})\leq N^{c/2k},
\eeq
Then \eqref{eq:cut_off_ll} implies $\P[\Xi^{c}]=O(N^{-D})$ for any fixed $D>0$. 

Assuming $i_{1}=2$ without loss of generality, we now estimate the expectation in \eqref{eq:cut_off_prf}. Restricting on the event $\Xi$, we have
\beq\label{eq:cut_off_prf1}
\E\lone_{\Xi}\prod_{j=1}^{k}I_{i_{j}}^{(j)}\leq N^{c(k-1)/(2k)}\E I_{1}^{(2)}\leq N^{-c/2}.
\eeq
On the complementary event, we use Cauchy-Schwarz inequality and \eqref{eq:abs} and then choose large enough $D>0$ to get 
\beq\label{eq:cut_off_prf2}
\E\lone_{\Xi^{c}}\prod_{j=1}^{k}I_{i_{j}}^{(j)}\lesssim \P[\Xi^{c}]^{1/2}N^{10}(1+\E(\norm{X}^{4})^{1/2})\lesssim N^{10-D/2}\lesssim N^{-c/2}.
\eeq
Combining \eqref{eq:cut_off_prf1} and \eqref{eq:cut_off_prf2} concludes the proof of Lemma \ref{lem:cut_off}.
\end{proof}

The remaining content of this section as well as Sections \ref{sec:comp_Gini}--\ref{sec:short GFT} are devoted to the proof of Theorem \ref{thm:comp_Gini}. The proof is first divided into two parts by considering $A+X^{\mathrm{Gin}(\bbF)}$ as an intermediate ensemble between $A+X$ and $X^{\mathrm{Gin}(\bbF)}$. To be more precise, we introduce the intermediate functions
\beq\begin{aligned}
L_{0}^{\mathrm{Gin}(\bbF)}(w)&\deq \Tr \log(\absv{A+X^{\mathrm{Gin}(\bbF)}-z_{0}-\gamma_{0}^{-1}N^{-1/2}w}^{2}+\eta_{0}^{2})	\\
&\qquad\qquad\qquad -\Tr\brkt{\log(\absv{A+\bsx-z_{0}-\gamma_{0}^{-1}N^{-1/2}w}^{2}+\eta_{0}^{2})}_{\caM},	\\
N_{0}^{\mathrm{Gin}(\bbF)}(w)&\deq \frac{1}{c_{0}}\Brkt{\frac{\eta_{0}}{\absv{A+X^{\mathrm{Gin}(\bbF)}-z_{0}-\gamma_{0}^{-1}N^{-1/2}w}^{2}+\eta_{0}^{2}}},
\end{aligned}
\eeq
obtained by replacing $X$ with $X^{\mathrm{Gin}(\bbF)}$ in $L_{0}$ and $N_{0}$. We now state the two key steps of the proof of Theorem \ref{thm:comp_Gini} as follows, whose combination immediately gives Theorem~\ref{thm:comp_Gini}: 
\begin{prop}\label{prop:comp_G}
Let $k\in\N$ and $C>0$ be fixed. Then there exists a constant $c>0$ such that the following hold uniformly over $w\in\C$ and $\bsw=(w_{1},\ldots,w_{k})\in\C^{k}$ with $\absv{w},\norm{\bsw}\leq C$:
\begin{align}
	\E\prod_{i=1}^{k}L_{0}^{\mathrm{Gin}(\bbF)}(w_{i})=&\E\prod_{i=1}^{k}L_{1}^{\mathrm{Gin}(\bbF)}(w_{i})+h^{(0)}(\bsw)+O(N^{-c}),\label{eq:comp_logd_G}\\
	\E N_{0}^{\mathrm{Gin}(\bbF)}(w)
	=&\E N_{1}^{\mathrm{Gin}(\bbF)}(w)+O(N^{-c}\Psi), \label{eq:comp_sing_G}
\end{align}
where $h^{(0)}\equiv h_{N}^{(0)}\in C^{2k}(\C^{k})$ satisfies $(\prod_{j=1}^{k}\Delta_{w_{j}})h^{(0)}\equiv0$.
\end{prop}

\begin{prop}\label{prop:comp_IID}
Let $k\in\N$ and $C>0$ be fixed. Then there exists a constant $c>0$ such that the following hold uniformly over $w\in\C$ and $\bsw=(w_{1},\ldots,w_{k})\in\C^{k}$ with $\absv{w},\norm{\bsw}\leq C$:
\begin{align}
	\E\prod_{i=1}^{k}L_{0}(w_{i})=&\E\prod_{i=1}^{k}L_{0}^{\mathrm{Gin}}(w_{i})+O(N^{-c}),\label{eq:comp_logd_IID}\\
	\E N_{0}(w)
	=&\E N_{0}^{\mathrm{Gin}(\bbF)}(w)+O(N^{-c}\Psi), \label{eq:comp_sing_IID}
\end{align}
\end{prop}
We prove Proposition \ref{prop:comp_G} for $\bbF=\C$ and $\bbF=\R$ respectively in Sections \ref{sec:comp_Gini} and \ref{sec:real} since the real symmetry class requires a quite different and more involved treatment. The proof of Proposition \ref{prop:comp_IID} is independent of the symmetry class and it is postponed to Section \ref{sec:short GFT}.

\subsection{Deterministic flow}\label{sec:flow}
In this section, we introduce further notations that are commonly used in Sections~\ref{sec:comp_Gini}--\ref{sec:real} along the proof of Proposition \ref{prop:comp_G}. We also present preliminary reductions to the proof. Since the IID matrix $X$ is always a Ginibre ensemble in this section (in fact up to Section \ref{sec:real}), we omit the superscripts $\mathrm{Gin}$ and $\mathrm{Gin}(\bbF)$.

The comparison in Proposition \ref{prop:comp_G} is performed along a piecewise $C^{1}$ flow of deterministic matrices, denoted by $\caA_{t}:[0,1]\to\C^{N\times N}$. The specific choice of the path $\caA_{t}$ differs by symmetry classes due to an additional constraint in the real case, and we make the actual choice at the beginning of Sections \ref{sec:comp_Gini} and \ref{sec:real}. Nonetheless, they all satisfy the same condition that the origin remains a sharp edge along the entire flow; we state this in the following assumption:
\begin{assump}\label{assump:path} 
	The path $[0,1]\ni t\mapsto \caA_{t}\in\C^{N\times N}$ is piecewise $C^{1}$ and it satisfies the following:
	\begin{itemize}
		\item[(i)] $\caA_{0}=A-z_{0}$ and $\caA_{1}=-z_{1}\in\bbS^{1}$.
		\item[(ii)] For all $t\in[0,1]$, we have
		\beq	\label{eq:solzt}
		\Brkt{\frac{1}{\absv{\caA_{t}}^{2}}}= 1.
		\eeq
		\item[(iii)] There exists a constant $C_{5}>1$  such that for all but finitely many $t\in[0,1]$ we have 
		\beq\label{eq:A_bdd}
		\norm{\caA_{t}}\leq C_{5}, \qquad
		\norm{\caA_{t}^{-1}}\leq C_{5},\qquad 
		C_{5}^{-1}\leq \absv{\brkt{\caA_{t}\absv{\caA_{t}}^{-4}}}\leq C_{5},	\qquad
		\Norm{\frac{\dd}{\dd t}\caA_{t}}\leq C_{5}\log N.
		\eeq
	\end{itemize}
\end{assump}
For such a choice of $\caA_{t}$, we will work with the continuous interpolation $\caA_{t}+X$ connecting $A+X-z_{0}$ to $X-z_{1}$. The scaling parameters of this interpolation are as follows; we introduce the time-dependent counterpart of the auxiliary quantities $I_{3},I_{4}$ in \eqref{eq:def_I3I4_0}:
\beq\label{eq:def_I}
I_{3}\equiv I_{3}(t)\deq \Brkt{\caA_{t}\adj\frac{1}{\absv{\caA_{t}}^{4}}},\qquad 
I_{4}\equiv I_{4}(t)\deq \Brkt{\frac{1}{\absv{\caA_{t}}^{4}}},
\eeq
and define $c_{t}$ and $\gamma_{t}$ as
\beq\label{eq:def_gamma}
c_{t}\deq I_{4}(t)^{-1/4},\qquad  \gamma_{t}\deq -I_{4}(t)^{-1/2}I_{3}(t).
\eeq
Furthermore, we point out that, uniformly in $t\in [0,1]$, $I_3(t)\sim 1$ by the third relation in \eqref{eq:A_bdd}, and that $I_4(t)\sim 1$ by
\begin{equation}
\label{eq:i4sim1}
1=\Brkt{\frac{1}{|\mathcal{A}_t|^2}}\le \Brkt{\frac{1}{|\mathcal{A}_t|^4}}\le C_5^4,
\end{equation}
where the last inequality follows by the second relation in \eqref{eq:A_bdd}.
Next, for $w\in\C$ we introduce the Hermitizations and the regularizing parameter $\eta_{t}$
\beq\begin{aligned}\label{eq:def_Herm}
	W&\deq \begin{pmatrix}
		0 & X \\ 
		X\adj & 0
	\end{pmatrix}, \\
	H_{t}^{w}&\deq\begin{pmatrix}
		0 & \caA_{t}+X-\gamma_{t}^{-1}N^{-1/2}w \\
		(\caA_{t}+X-\gamma_{t}^{-1}N^{-1/2}w)^{\adj} & 0 
	\end{pmatrix},	\\
	\bsh_{t}^{w}&\deq\begin{pmatrix}
		0 & \caA_{t}+\bsx-\gamma_{t}^{-1}N^{-1/2}w \\
		(\caA_{t}+\bsx-\gamma_{t}^{-1}N^{-1/2}w)^{\adj} & 0 
	\end{pmatrix},\\
	J_{t}^{w}&\deq \begin{pmatrix}
		0 & \gamma_{t}^{-1}w \\
		\ol{\gamma}_{t}^{-1}\ol{w} & 0
	\end{pmatrix},	\\
	\eta_{t}&\deq c_{t}^{-1}\eta_{1}.
\end{aligned}\eeq
We denote the resolvent of $H_{t}^{w}$ (restricted to the imaginary axis) by
\beq
G_{t}^{w}(\ii\eta)=(H_{t}^{w}-\ii\eta)^{-1},\qquad \eta>0.
\eeq
We also introduce the interpolation between $(L_{0}(w),N_{0}(w))$ and $(L_{1}(w),N_{1}(w))$ as 
\beq\label{eq:def_Lt}\begin{aligned}
	L_{t}(w)\deq \Tr \log(\absv{H_{t}^{w}-\ii\eta_{t}})-\Tr\brkt{\log(\absv{\bsh_{t}^{w}-\ii\eta_{t}})}_{\caM},	\\
	N_{t}(w)\deq \frac{1}{c_{t}}\Brkt{\frac{\eta_{t}}{\absv{\caA_{t}+X-\gamma_{t}^{-1}N^{-1/2}w}^{2}+\eta_{t}^{2}}},
\end{aligned}\eeq
where $\brkt{\cdot}_{\caM}$ is applied entrywise as described in Definition \ref{defn:vN}. Notice the identity
\beq\label{eq:N=G}
N_{t}(w)=\frac{1}{c_{t}}\brkt{\im G_{t}^{w}(\ii\eta_{t})},
\eeq
and that $L_{t}(w)$ can be expressed as an integral of $\brkt{G_{t}^{w}(\ii\eta)-M_{t}^{w}(\ii\eta)}$ in $\eta$, i.e.
\beq\label{eq:dlog=G}
L_{t}(w)=2N\int_{\eta_{t}}^{\infty}\im\brkt{G_{t}^{w}(\ii\eta)-M_{t}^{w}(\ii\eta)}\dd\eta.
\eeq

Note that the matrix $H_{t}^{w}$ as well as the operator $\bsh_{t}^{w}$ depend on $t$ only deterministically, i.e.
\beq\label{eq:H_t_determ}
\frac{\dd H_{t}^{w}}{\dd t}=\frac{\dd \E H_{t}^{w}}{\dd t}=\frac{\dd\brkt{\bsh_{t}^{w}}_{\caM}}{\dd t}=\frac{\dd \bsh_{t}^{w}}{\dd t}.
\eeq
The scaling parameters in \eqref{eq:def_gamma} are chosen so that the cubic equation \eqref{eq:cubic} for the rescaling of the regularized eigenvalue density of $H_{t}^{w}$,  $v_{t}^{w}\deq\im\brkt{M_{\caA_{t},N^{-1/2}\gamma_{t}^{-1}w}(\ii\eta_{t})}+\eta_{t}$, becomes time-independent, i.e. we have
\beq
\left(\frac{v_{t}^{w}}{c_{t}}\right)^{3}+2N^{-1/2}\re[w]\left(\frac{v_{t}^{w}}{c_{t}}\right)-\eta_{1}=O(\Psi^{5})
\eeq
thus $v_{t}^{w}/c_{t}$ is essentially constant in $t$. At the same time, the choice of the scaling parameters also makes the regularized, rescaled, local eigenvalue density
\beqs
w\mapsto \frac{N}{4\pi}\Delta_{w}\brkt{\log(\absv{\caA_{t}+X-\gamma_{t}^{-1}N^{-1/2}w}^{2}+\eta_{t}^{2})}
\eeqs
almost constant over $t$ for each $w$; see Lemma \ref{lem:determ}. 

Recall from Lemma \ref{lem:MDE} that the $\C^{2N\times 2N}$-valued Stieltjes transform
\beq
\Brkt{(\bsh_{t}^{w}-\wh{z})^{-1}}_{\caM}\in\C^{2N\times 2N},\qquad \wh{z}\in\C_{+}
\eeq
is the unique solution $M\equiv M_{\caA_{t},N^{-1/2}\gamma_{t}^{-1}w}(\wh{z})\in\C^{2N\times 2N}$ to the MDE
\beq\label{eq:MDE}
\frac{1}{M}=\E H_{t}^{w}-\wh{z}-\brkt{M},\qquad \im M>0.
\eeq
In accordance with \eqref{eq:def_Herm}, we write the $w$-dependent solution $M$
\beq
M_{t}^{w}(\wh{z})\deq M_{\caA_{t},N^{-1/2}\gamma_{t}^{-1}w}(\wh{z}).
\eeq
From the properties of $\caA_{t}$, we have the following time dependent version of Theorem \ref{thm:local law: Local Law}.
\begin{prop}\label{prop:time dependent local law}
	For any fixed large $C>0$, uniformly in deterministic matrices $B\in\C^{N\times N}$ and deterministic vectors $\bsx, \bsy\in\C^N$, we have
	\begin{equation}\label{eq:time dependent local law}
		\absv{ \brkt{B(G_{t}^{w}(\ii\eta)-M_{t}^{w}(\ii\eta) ) } }\prec \frac{\lVert B\rVert}{N\eta}, \qquad \left|\bsx^*\left(G_{t}^{w}(\ii\eta)-M_{t}^{w}(\ii\eta)\right)\bsy \right|\prec \Norm{\bsx }\Norm{\bsy}\sqrt{\frac{\rho^{w,t}(\ii\eta)}{N\eta}},
	\end{equation} 
	uniformly over $\caA_{t}$ satisfying Assumption \ref{assump:path}, $t\in[0,1]$, $\absv{w}\leq C$, and $\eta\in[N^{-1},N^{100}]$. Consequently, the following holds uniformly over $\absv{w}\leq C$ and $t\in[0,1]$:
	\begin{equation}\label{eq:Lt_bdd}
		\absv{L_{t}(w)}\prec 1.
	\end{equation}
\end{prop}
\begin{proof}
	For each fixed $t\in[0,1]$, as the asymptotic of Lemma \ref{lem:MDE_asymp} holds for $t$ this is a straightforward generalization of Theorem \ref{thm:local law: Local Law} by replacing $A-z_0$ by $\caA_{t}$. Taking a sufficiently large, but still polynomial sized, grid of $[0,1]$ and union bounding over the very high probability events on which \eqref{eq:time dependent local law} holds we see that \eqref{eq:time dependent local law} holds uniformly over $t\in[0,1]$. To see \eqref{eq:Lt_bdd}, we apply \eqref{eq:time dependent local law} to \eqref{eq:dlog=G}, so that 
	\beq
	L_{t}(w)=2N\int_{\eta_{t}}^{\infty}\left[\im\brkt{G_{t}^{w}(\ii\eta)-M_{t}^{w}(\ii\eta)}\right]\dd\eta
	\prec2N\int_{\eta_{t}}^{N^{100}}\frac{1}{N\eta}\dd\eta+N\int_{N^{100}}^{\infty}\frac{1}{\eta^{2}}\dd\eta
	\lesssim \log N.
	\eeq
\end{proof}

Notice from Lemma \ref{lem:MDE_asymp} that $\absv{\brkt{M_{t}^{w}(\ii\eta)}}=\im\brkt{M_{t}^{w}(\ii\eta)}=O(N^{-1/4})$ for $\absv{w}=O(1)$ and $\eta=O(N^{-3/4})$. We often use the abbreviation
\beq\label{eq:def_M_bare}
\mr{M}_{t}^{w}\deq (\E H_{t}^{w})^{-1}\equiv(\brkt{\bsh_{t}^{w}}_{\caM})^{-1}\in\C^{2N\times 2N}.
\eeq
The matrix $\mr{M}_{t}^{w}$ solves the degenerate MDE $\mr{M}^{-1}=\E H_{t}^{w}$. We may consider $\mr{M}$ as a perturbation of the solution $M$ of the true MDE, i.e. for $\absv{w}=O(1)$ and $\eta=O(N^{-3/4})$ we have
\beq\label{eq:M_bare_close}
\norm{M_{t}^{w}(\ii\eta)-\mr{M}_{t}^{w}}\lesssim \eta+\brkt{\im M_{t}^{w}(\ii\eta)}\lesssim N^{-1/4}.
\eeq
Consequently, in the regime $\eta=O(N^{-3/4})$ and $\absv{w}=O(1)$, $G_{t}^{w}(\ii\eta)-\mr{M}_{t}^{w}$ satisfies the same local law Proposition \ref{prop:time dependent local law} as $G_{t}^{w}(\ii\eta)-M_{t}^{w}(\ii\eta)$ since the errors in \eqref{eq:time dependent local law} are larger than $N^{-1/4}$.

In the next lemma, we control the time derivative of the deterministic, regularized Brown measure density.
\begin{lem}\label{lem:determ}
	The following holds uniformly over the path $\caA_{t}$ satisfying Assumption \ref{assump:path} and $w$ in a compact subset of $\C$.
	\beq\label{eq:determ}
	\frac{\dd\brkt{\brkt{\log\absv{\bsh_{t}^{w}-\ii\eta_{t}}}_{\caM}}}{\dd t}
	=-\frac{\dd\brkt{\log\absv{\mr{M}_{t}^{w}}}}{\dd t}+O(N^{-3/2}).
	\eeq
	Moreover, the function $w\mapsto \brkt{\log \absv{\mr{M}_{t}^{w}}}$, and hence $\dd\brkt{\log\absv{\mr{M}_{t}^{w}}}/\dd t$, are harmonic in $w$ over a compact subset of $\C$.
\end{lem}
\begin{proof}	
	We write the time derivative on the left-hand side of \eqref{eq:determ} as
	\beq\label{eq:determ_1}
	\begin{split}
		\frac{1}{2}\frac{\dd}{\dd t}\brkt{\brkt{\log(\absv{\bsh_{t}^{w}}^{2}+\eta_{t}^{2})}_{\caM}}
		=\Brkt{\Brkt{\frac{\dd (\bsh_{t}^{w}-\ii\eta_{t})}{\dd t}\frac{1}{\bsh_{t}^{w}-\ii\eta_{t}}}_{\caM}}
		=\Brkt{\frac{\dd(\E H_{t}^{w}-\ii\eta_{t})}{\dd t}M_{t}^{w}},
	\end{split}
	\eeq
	where we used in the first equality that $\bsh_{t}^{w}$ is supported strictly on the off-diagonal blocks, and in the second equality \eqref{eq:H_t_determ}.
	
	To see \eqref{eq:determ}, we write $\ii v\equiv \ii v_{t}^{w}\deq \ii\eta_{t}+\brkt{M_{t}^{w}}$ so that the MDE \eqref{eq:MDE} becomes 
	\beq
	\frac{1}{M_{t}^{w}}=\frac{1}{\mr{M}_{t}^{w}}-\ii v.
	\eeq
	Since $\absv{v}=O(N^{-1/4})$ by \eqref{eq:M_size}, we have the geometric expansion up to the sixth order
	\beq\label{eq:determ_expa}
	\Norm{M_{t}^{w}-\mr{M}_{t}^{w}\sum_{k=0}^{5}(\ii v\mr{M}_{t}^{w})^{k}}=O(v^{6})=O(N^{-3/2}).
	\eeq
	Substituting \eqref{eq:determ_expa} into \eqref{eq:determ_1} gives
	\beq\begin{split}\label{eq:determ_3}
		\Brkt{\frac{\dd (\E H_{t}^{w}-\ii\eta_{t})}{\dd t}M_{t}^{w}}
		=&-\frac{\dd}{\dd t}\brkt{\log\absv{\mr{M}_{t}^{w}}}+\frac{v^{2}}{2}\frac{\dd}{\dd t}\brkt{(\mr{M}_{t}^{w})^{2}}	\\
		&-\frac{v^{4}}{4}\frac{\dd}{\dd t}\brkt{(\mr{M}_{t}^{w})^{4}}+v\frac{\dd\eta_{t}}{\dd t}\Brkt{(\mr{M}_{t}^{w})^{2}}+O(N^{-3/2}\log N),
	\end{split}\eeq
	where the contribution of (most) odd powers of $v$ disappears since $\E H$ and $\mr{M}=(\E H)^{-1}$ are supported on off-diagonal blocks, and we also used $\eta_{t}=c_{t}\eta_{1}$, $\eta_{1}=O(N^{-3/4})$.
	
	Next, we expand \eqref{eq:determ_3} around $w=0$; (recall the definitions of $I_{3},I_{4},\gamma_{t}$ and $J_{t}^{w}$ from \eqref{eq:def_I}--\eqref{eq:def_Herm})
	\beq\begin{split}\label{eq:determ_4}
		\brkt{(\mr{M}_{t}^{w})^{2}}=&\brkt{(\mr{M}_{t}^{0})^{2}}
		+2\Brkt{J_{t}^{w}(\mr{M}_{t}^{0})^{3}}+O(N^{-1})	\\
		=&1+2N^{-1/2}\re\left[\frac{w}{\gamma_{t}}I_{3}\right]+O(N^{-1})=1+2N^{-1/2}I_{4}^{1/2}\re[-w]+O(N^{-1}),\\
		\brkt{(\mr{M}_{t}^{w})^{4}}=&\brkt{(\mr{M}_{t}^{0})^{4}}+O(N^{-1/2})=I_{4}+O(N^{-1/2}).
	\end{split}\eeq
	Following the same calculations for the time derivatives of the left-hand sides of \eqref{eq:determ_4}, using the extra input from the last estimate of \eqref{eq:A_bdd}, we arrive at the time derivatives of the right-hand sides of \eqref{eq:determ_4} with the same errors except with an extra factor of $\log N$. Combining this with $\eta_{t}=I_{4}^{1/4}\eta_{1}$, we obtain
	\beq\label{eq:determ_5}\begin{split}
		&\Brkt{\frac{\dd (\E H_{t}^{w}-\ii\eta_{t})}{\dd t}M_{t}^{w}}	\\
		=&-\frac{\dd}{\dd t}\brkt{\log\absv{\mr{M}_{t}^{w}}}
		+N^{-1/2}v^{2}\re[-w]\frac{\dd I_{4}^{1/2}}{\dd t}
		-\frac{v^{4}}{4}\frac{\dd I_{4}}{\dd t}+v\eta_{1}\frac{\dd I_{4}^{1/4}}{\dd t}+O(N^{-3/2}\log N)	\\
		=&-\frac{\dd}{\dd t}\brkt{\log\absv{\mr{M}_{t}^{w}}}
		-\frac{v}{4 I_{4}}\frac{\dd I_{4}}{\dd t}\left(I_{4}v^{3}+2N^{-1/2}I_{4}^{1/2}v\re[w]-\eta_{t}\right)+O(N^{-3/2}\log N).
	\end{split}\eeq
	
	To conclude \eqref{eq:determ}, notice that the second term of \eqref{eq:determ_5} is exactly the cubic equation from Lemma \ref{lem:MDE_asymp}:
	\beq
	I_{4}v^{3}+2N^{-1/2}I_{4}^{1/2}v\re[w]-\eta_{t}
	=I_{4}v^{3}-2v\re\left[{N^{-1/2}\gamma_{t}^{-1}w}I_{3}\right]-\eta_{t}=O(N^{-5/4}).
	\eeq
	Since $v=O(N^{-1/4})$ and $\dd I_{4}/\dd t=O(1)$, this proves 
	\beq
	\Brkt{\frac{\dd (\E H_{t}^{w}-\ii\eta_{t})}{\dd t}M_{t}^{w}}=-\frac{\dd}{\dd t}\brkt{\log\absv{\mr{M}_{t}^{w}}}+O(N^{-3/2}\log N),
	\eeq
	as desired.
	
	Finally, to see the harmonicity, notice from $\norm{\caA_{t}^{-1}}\leq C$ that $\min_{i}\absv{\sigma_{i}(\caA_{t})}\geq 1/C$, so that the map
	\beq
	w\mapsto \brkt{\log\absv{\mr{M}_{t}^{w}}}=\frac{1}{N}\sum_{i=1}^{N}\log\absv{\sigma_{i}(\caA_{t})-N^{-1/2}\gamma_{t}^{-1}w}
	\eeq
	is harmonic on $\{w:\absv{w}\leq C\sqrt{N}\absv{\gamma_{t}}\}$.
\end{proof}

In the last bit of this section, we further reduce the proof of Proposition \ref{prop:comp_G} into an estimate of its time derivative. To this end, we introduce the following notation\footnote{In Lemma \ref{lem:dL=dHG} we only need $\absv{I}=1$, but $\absv{I}$ can reach (at most) five later in Sections \ref{sec:comp_Gini} and \ref{sec:real}.} for simplicity: For an index set $I=\{i_{1},\ldots,i_{n}\}\subset\bbrktt{k}$, we write
\beq\begin{aligned}\label{eq:def_L}
	\bsw^{(I)}&\equiv \bsw^{(i_{1},\ldots,i_{n})}\deq(w_{j})_{j\notin I}\in\C^{\bbrktt{k}\setminus I}, \qquad & &\bsw\deq \bsw^{\emptyset},	\\
	\bsL^{(I)}_{t}(\bsw^{(I)})&\equiv\bsL_{t}^{(i_{1},\ldots,i_{n})}(\bsw^{(I)})\deq \prod_{\substack{\ell=1, \\ \ell\notin I}}^{k}L_{t}(w_{\ell}), \qquad& &\bsL_{t}(\bsw)\deq \bsL_{t}^{(\emptyset)}(\bsw).
\end{aligned}\eeq
Now we compute the time derivative of $\bsL_{t}$, with the help of Lemma \ref{lem:determ}, as follows:
\begin{lem}\label{lem:dL=dHG}
	For any fixed $\epsilon>0$, the following holds uniformly over $t\in[0,1]$ and $\absv{w},\norm{\bsw}\leq C$:
	\begin{align}
		\E\frac{\dd\bsL_{t}(\bsw)}{\dd t}=\sum_{j=1}^{k}\E\Tr\left[\frac{\dd(\E H_{t}^{w_{j}}-\ii\eta_{t})}{\dd t}G^{w_{j}}_{t}(\ii\eta_{t})\right]\bsL_{t}^{(j)}(\bsw^{(j)})
		&+h_{t}(\bsw)+O(N^{-1/2+\epsilon}),	\label{eq:dL=dHG}
	\end{align}
	where $\Delta_{w_{1}}\cdots\Delta_{w_{k}}h_{t}\equiv 0$.
\end{lem}
\begin{proof}
	From the first-order perturbation theory
	\beq
	\log\det (1+\dd Y)=\Tr\log(1+\dd Y)=\Tr \dd Y+O(\norm{\dd Y}^{2}),\qquad \dd Y\in\C^{2N\times 2N},
	\eeq
	we immediately have
	\beq\begin{split}
		\frac{\dd L_{t}(w)}{\dd t}
		=&\Tr\frac{\dd (H_{t}^{w}-\ii\eta_{t})}{\dd t} G_{t}^{w}-2N\frac{\dd\brkt{\brkt{\log\absv{\bsh_{t}^{w}-\ii\eta_{t}}}_{\caM}}}{\dd t}	\\
		=&\Tr\frac{\dd (\E H_{t}^{w_{j}}-\ii\eta_{t})}{\dd t} G_{t}^{w}-2N\frac{\dd\brkt{\log\absv{\mr{M}_{t}^{w}}}}{\dd t}+O(N^{-1/2}),
	\end{split}\eeq
	where in the second line we used \eqref{eq:H_t_determ} and Lemma \ref{lem:determ}. On the other hand, recall from \eqref{eq:Lt_bdd} that $L_{t}(w)\prec 1$.
	
	Therefore, combining with Lemma \ref{lem:determ}, we get
	\beq\label{eq:dR/dt1}
	\begin{split}
		\E \frac{\dd}{\dd t}\bsL_{t}(\bsw)=&\E\sum_{j=1}^{k}\left(\Tr \frac{\dd(\E H_{t}^{w_{j}}-\ii\eta_{t})}{\dd t}G_{t}^{w_{j}}-2N\frac{\dd\brkt{\log\absv{\mr{M}_{t}^{w_{j}}}}}{\dd t}\right) \bsL_{t}^{(j)}(\bsw^{(j)})+O(N^{-1/2+\epsilon}).
	\end{split}
	\eeq
	Denote the sum of the second terms in the summation in \eqref{eq:dR/dt1} by $h_{t}(\bsw)$. Notice that each of these summands is harmonic in $w_{j}$ by Lemma \ref{lem:determ}. This proves $\Delta_{w_{1}}\cdots\Delta_{w_{k}}h_{t}\equiv 0$.
\end{proof}

\section{Proof of Proposition \ref{prop:comp_G}, the complex case}
\label{sec:comp_Gini}
This section is devoted to the proof of Proposition \ref{prop:comp_G} when $\bbF=\bbC$. Since we only consider deformed complex Ginibre ensembles, we omit the superscripts $\mathrm{Gin}$ and $\mathrm{Gin}(\C)$ throughout the section.
We start with the construction of the path $\caA_{t}$.

\begin{prop}\label{prop:z}
	If $z_{0}$ satisfies Assumption \ref{assump:edge} (i), (ii), and (iii$\C$), then there exist a path $\caA_{t}$ satisfying Assumption \ref{assump:path}, with the constant $C_{5}$ in \eqref{eq:A_bdd} depending only on the constants from Assumption \ref{assump:edge}.
\end{prop}
\begin{proof}
	Recall the definition of the domain $\caD_{C}$ from Assumption \ref{assump:edge}. Notice from $\norm{A}\leq C$ that 
	\beq\label{eq:z_large}
	\norm{(A-z)^{-1}+z^{-1}}=O(\absv{z}^{-2}),\qquad \brkt{(A-z)\absv{A-z}^{-4}}=-z\absv{z}^{-4}+O(\absv{z}^{-4}),\qquad \text{as $\absv{z}\to\infty$}.
	\eeq
	This in particular implies that for large enough $C''>0$ the complementary ball $\{z:\absv{z}\geq C''\}$ is in $\caD_{C'}$, so that $\caD_{C'}$ has a unique unbounded component. We define the path $[0,1)\ni t\mapsto \wh{z}_{t}\in\caD_{C'}$ as follows.
	\begin{itemize}
		\item[(a)] For $0\leq t\leq 1/2$, the smooth path $\wh{z}_{t}$ connects $z_{0}$ to $z_{1/2}=C''z_{1}$
		\item[(b)] For $1/2\leq t<1$, the path $\wh{z}_{t}$ is given by $\wh{z}_{t}=C''z_{1}/(1-t)$.
	\end{itemize}
	Note that the existence of the path connecting $z_{0}$ to $C''z_1$ follows from Assumption \ref{assump:edge} (iii) and that $\{z:\absv{z}\geq C''\}\subset\caD_{C'}$. Finally, define
	\beq
	\caA_{t}\deq \brkt{\absv{A-\wh{z}_{t}}^{-2}}^{1/2}(A-\wh{z}_{t}),\quad t\in(0,1),\qquad \caA_{1}\deq -1.
	\eeq
	
	The first two conditions in Assumption \ref{assump:path} are immediate from the definition of $\caA_{t}$. For any given (small) $\delta>0$, the estimates in \eqref{eq:A_bdd}, by continuity, are trivial for $t\in[0,1-\delta)$ from the definition of $\caD_{C'}$ and $\absv{\wh{z}_{t}}=O(1)$. When $t$ is sufficiently close to $1$, the first three estimates in \eqref{eq:A_bdd} are consequences of the large-$\absv{z}$-asymptotics \eqref{eq:z_large} and
	\beq
	\norm{A-z}=\absv{z}+O(1),\qquad \brkt{\absv{A-z}^{-2}}=\absv{z}^{-2}+O(\absv{z}^{-3}).
	\eeq
	For the derivative of $\caA_{t}$ when $t$ is close to $1$, since $\wh{z}_{t}$ follows the ray $z_{1}\cdot(0,\infty)$ with $\absv{z_{1}}=1$, we get
	\beq\begin{split}
		\Norm{\frac{\dd}{\dd t}\caA_{t}}=&\Absv{\frac{\dd \wh{z}_{t}}{\dd t}}\Norm{\brkt{\absv{A-\wh{z}_{t}}^{-2}}^{-1/2}\re\Brkt{(A-\wh{z}_{t})\adj\frac{-z_{1}}{\absv{A-\wh{z}_{t}}^{-4}}}(A-\wh{z}_{t})+\brkt{\absv{A-\wh{z}_{t}}^{-2}}^{1/2}z_{1}}	\\
		=&\Absv{\frac{\dd \wh{z}_{t}}{\dd t}}O(\absv{\wh{z}_{t}}^{-2})(1+\norm{A})=O(1).
	\end{split}\eeq
	This concludes the proof of Proposition \ref{prop:z}.
\end{proof}

With \eqref{eq:N=G}, Lemma \ref{lem:dL=dHG}, and Proposition \ref{prop:z}, we immediately find that the next result implies Theorem \ref{thm:comp_Gini} when $\bbF=\bbC$. 
\begin{prop}\label{prop:dt}
	Let $\bbF=\bbC$, $k\in\N$, $C>0$, and $\epsilon>0$ be fixed. Then the following holds uniformly over $\caA_{t}$ satisfying Assumption \ref{assump:path}, $t\in[0,1]$, and $\absv{w},\norm{\bsw}\leq C$:
	\begin{align}
		\sum_{j=1}^{k}\E\Brkt{\frac{\dd(\E H_{t}^{w_{j}}-\ii\eta_{t})}{\dd t}G^{w_{j}}_{t}(\ii\eta_{t})}\bsL_{t}^{(j)}(\bsw^{(j)})&=h_{t}(\bsw)+O(N^{\epsilon}\Psi^{5}),	\label{eq:comp_logd_d}\\
		\frac{\dd}{\dd t}\E \frac{\brkt{G_{t}^{w}}}{c_{t}}&=O(N^{1+\epsilon}\Psi^{6}), \label{eq:comp_sing_d}
	\end{align}
	where $\Delta_{w_{1}}\cdots\Delta_{w_{k}}h_{t}\equiv 0$. 
\end{prop}
The rest of this section is devoted to the proof of Proposition \ref{prop:dt}, and all estimates depend only on its parameters $k,C,\epsilon$, and $C_{5}$ from Assumption \ref{assump:path}.
\subsection{General framework}\label{sec:comp_prep}
Recall that in \eqref{eq:determ_expa} we Taylor-expanded the \emph{matrix} $M_{t}^{w}$ in the trace
\beq
\Brkt{\frac{\dd \E H_{t}^{w}}{\dd t}M_{t}^{w}}
\eeq 
around $\mr{M}_{t}^{w}$ with respect to the \emph{scalar} $v^{2}$. The main goal of this section is to emulate the same expansion for \eqref{eq:comp_logd_d}.  While all our estimates shall be uniform over $t\in[0,1]$, we will consider $t$ to be fixed throughout this section. Henceforth we abbreviate $\eta\equiv\eta_{t}$, $G^{w_{j}}\equiv G_{t}^{w_{j}}$, $M^{w_{j}}\equiv M_{t}^{w_{j}}$, etc. Also, we omit the spectral parameters of $G^{w}$ and $M^{w}$ when they are $\ii\eta$.

Recall the definition of $W$ from \eqref{eq:def_Herm}. In what follows, we often use Stein's lemma. For a function $f:\C^{2N\times 2N}\to \C^{2N\times 2N}$, Stein's lemma states
\beq\label{eq:Stein}
\E Wf(W)=\E\E_{\wt{W}}\wt{W}(\partial_{\wt{W}}f(W)),
\eeq
where $\wt{W}$ is an i.i.d.\ copy of $W$ (hence has Gaussian entries) and $\partial_{\wt{W}}$ denotes the directional partial derivative along $\wt{W}$. As a simple consequence of \eqref{eq:Stein} together with the identity
\beq\label{eq:sc00}
G-M=M(M^{-1}-G^{-1})G=M(-W-\brkt{M})G,
\eeq
which follows from the defining equation of $M$ \eqref{eq:MDE0}, we have
\beq\label{eq:sc0}
\E[G-M]=\E[M\caS[G-M]G],
\eeq
where $\caS:\C^{2N\times 2N}\to\C^{2N\times 2N}$ is defined by
\beq
\caS[B]\deq \E_{\wt{W}}\wt{W}B\wt{W} =\begin{pmatrix} \brkt{2E_{2}B} & 0 \\ 0 & \brkt{2E_{1}B}\end{pmatrix}=\sum_{E=E_{1},E_{2}}E\brkt{2\check{E}B}.
\eeq
Here, the block diagonal matrices $E_{1},E_{2}$ are defined in \eqref{eq:def_E1} and the check denotes the involution $\check{E}_{1}\deq E_{2}$, $\check{E}_{2}=E_{1}$. Also note the identity
\beq\label{eq:quad}
\E_{\wt{W}}\wt{W}\Tr[\wt{W}B]=\frac{1}{N}\sum_{E=E_{1},E_{2}}EB\check{E}.
\eeq

For the matrix $\mr{M}=(\E H)^{-1}$, we have the following variant of \eqref{eq:sc00};
\begin{align}
	G-\mr{M}&=\mr{M}(\E H-G^{-1})G=\mr{M}(\ii\eta-W)G.	\label{eq:sc1}
\end{align}
Note also that
\beq\label{eq:expa_J}
\Norm{\mr{M}^{w}-\left(\mr{M}^{0}+\mr{M}^{0}J^{w}\mr{M}^{0}\right)}=O(\norm{J^{w}}^{2})=O(N^{-1}),
\eeq
where $J^{w}$ was defined in \eqref{eq:def_Herm}. We point out that to obtain \eqref{eq:expa_J} we used a second order resolvent expansion and the fact that $\norm{J^{w}}^{2}=O(N^{-1})$, since each of its element are smaller than $N^{-1/2}$ (recall that $\eta\ll N^{-3/4}$ and $|w|=O(N^{-1/2})$.

Our `Taylor expansion' expresses \eqref{eq:comp_logd_d} as a deterministic linear combination with each term of the form
\beq
\E\brkt{GE\cdots GE}\cdots\brkt{GE\cdots GE}\bsL^{(I)},
\eeq
for $G=G^{w_{i}}$ with varying $i$'s, $E\in\{E_{1},E_{2}\}$, and $I\subset\bbrktt{k}$. We refer to such an expression as `decoupled', in the sense that it involves only the plain trace $\brkt{\cdot}$ as opposed to a generic weighted trace $\brkt{B\,\,\cdot }$. To rigorously formulate our method, we introduce a few notations. In what follows, bold-faced symbols $\bsw$, and $\bsF$ always denote finite tuples (with appropriate dimensions) of elements of
\beq\label{eq:symbol}
\{w_{1},\ldots,w_{k}\} \qquad\text{and}\quad \{E_{1},E_{2}\},
\eeq
respectively. Here $w_{i}$'s are free (in the literal sense) complex variables (e.g. in Lemma \ref{lem:decoup}), unless specified to be in some subset of $\C$ (e.g. in Lemma \ref{lem:BP_rough}).
\begin{defn}
	For $n\in\N$, $\bsw=(\wt{w}_{j})_{j\in\bbrktt{n}}$, and $\bsF=(F_{j})_{j\in\bbrktt{n}}$, define
	\beq\label{eq:P_def}
	P_{n}(\bsw;\bsF)\deq \frac{1}{N^{n-1}}G^{\wt{w}_{1}}F_{1}\cdots G^{\wt{w}_{n}}F_{n}-\lone(n=1)\mr{M}^{\wt{w}_{1}}F_{1}.
	\eeq
	Similarly, for a multi-index $\bsm\equiv (m_{j})$ in $\N$, $\bsw\equiv (\bsw_{j})$, and $\bsF\equiv (\bsF_{j})$, define
	\beq\label{eq:Q_def}
	\bsQ_{\bsm}(\bsw;\bsF)\deq\prod_{j}\brkt{P_{m_{j}}(\bsw_{j};\bsF_{j})},\qquad \absv{\bsm}\deq \sum_{j}m_{j},
	\eeq
	with the convention $\bsQ_{\bsm}\equiv 1$ for $\absv{\bsm}=0$. For $\ell\in\N\cup\{0\}$ define 
	\beq
	\bsL_{\ell}(\bsw)\deq \prod_{j=1}^{\ell}L(\wt w_{j}).
	\eeq
	We often omit the subscript $\ell$ to write $\wt{\bsL}\equiv \wt{\bsL}_{\ell}$. Finally, $\bsQ_{\bsm}^{(I)}$ and $\wt{\bsL}^{(I)}$ are defined in the same fashion as \eqref{eq:def_L}. 
\end{defn}
\begin{defn}\label{def:caX_caY}
	For $n\in\N$, we define $\caX_{n}$ by the set of expressions 
	\beq\begin{split}\label{eq:def_caX}
		\{\bsQ_{\bsm}(\bsw_{1};\bsF_{1})\wt{\bsL}_{\ell}(\bsw_{2}):\,\,\, \absv{\bsm}=n, \ell\geq0,\bsw_{1},\bsw_{2},\bsF_{1}\},
	\end{split}\eeq
	where the components of $\bsw_{1},\bsw_{2}$, and $\bsF_{1}$ may vary among the symbols in \eqref{eq:symbol}.
	For a given matrix $B\in\C^{N\times N}$ and $n\in\N$, we define $\caY_{n}(B)$ to be the set of expressions
	\beq\begin{split}\label{eq:def_caY}
		\{\brkt{BP_{n_{1}}(\bsw_{1};\bsF_{1})}\bsQ_{\bsm}(\bsw_{2};\bsF_{2})\wt{\bsL}_{\ell}(\bsw_{3}):\,\,\, n_{1}+\absv{\bsm}=n,\ell\geq 0,\bsw_{1},\bsw_{2},\bsw_{3},\bsF_{1},\bsF_{2}\},
	\end{split}\eeq
	where the components of $\bsw_{1},\bsw_{2},\bsw_{3},\bsF_{1}$, and $\bsF_{2}$ may vary among the symbols in \eqref{eq:symbol}.
	Notice that $\caY_{n}(I)=\caX_{n}$. With a slight abuse of notations, we write $\frX_{n}$ (resp. $\frY_{n}(B)$) for generic finite sums (with the number of summand independent of $N$) of expressions in $\caX_{n}$ (resp. $\caY_{n}(B)$), whose exact form may vary by lines.
\end{defn}
The notations $\frX_{n}$ and $\frY_{n}(B)$ in Definition \ref{def:caX_caY} are introduced purely for the sake of brevity. For all instances of $\frX_{n}$ and $\frY_{n}(B)$, we can track down their exact forms as sums over elements of $\caX_{n}$ or $\caY_{n}(B)$. However we do not write them down simply because the explicit formulas play no role in the proof while they are excessively lengthy.

We may now present the two inputs that form the basis of our expansion: Their proofs are postponed to the end of this section.
\begin{lem}\label{lem:BP_rough}
	For any fixed $B\in\C^{N\times N}$ and $n\in\N$, we have
	\beq\label{eq:BP_lead}
	\absv{\wt{\bsL}_{n}(\bsw)} \prec1, \qquad \absv{\brkt{BP_{n}(\bsw;\bsF)}}\prec \norm{B}\Psi^{n}
	\eeq
	uniformly over $(w_{1},\ldots,w_{k})$ in a compact set.
\end{lem}
\begin{lem}[Decoupling]\label{lem:decoup}
	Let $B_{d}$ and $B_{o}\in\C^{N\times N}$ be supported respectively on diagonal and off-diagonal blocks. Then
	\begin{align}\begin{split}\label{eq:decoup_diag}
			\E\brkt{B_{d}P_{n}}\bsQ_{\bsm}\wt{\bsL}
			=&\,\,\brkt{2B_{d}\mr{M}^{\wt{w}_{1}}\mr{M}^{\wt{w}_{n}}F_{n}}\E\brkt{P_{n}}\bsQ_{\bsm}\wt{\bsL}	\\
			&+\left(\lone(n=1)\ii\eta_{t}\brkt{B_{d}(\mr{M}^{\wt{w}_{1}})^{2}F_{1}}+\lone(n=2)\frac{1}{N}\brkt{B_{d}\mr{M}^{\wt{w}_{1}}F_{1}\mr{M}^{\wt{w}_{2}}F_{2}}\right)\E\bsQ_{\bsm}\wt{\bsL}	\\
			&+\frac{1}{2}\E\frY_{\absv{\bsm}+n+1}(B_{d}\mr{M}^{\wt{w}_{1}})\\
			&+\ii\eta\E\frY_{\absv{\bsm}+n}(B_{d}\mr{M}^{\wt{w}_{1}})
			+\frac{1}{N}\E\frY_{\absv{\bsm}+n-1}(B_{d}\mr{M}^{\wt{w}_{1}}),
		\end{split}	\\
		\begin{split}\label{eq:decoup_off}
			\E\brkt{B_{o}P_{n}}\bsQ_{\bsm}\wt{\bsL}
			=&\,\,\frac{1}{2}\E\frY_{\absv{\bsm}+n+1}(B_{o}\mr{M}^{\wt{w}_{1}})	\\
			&+\ii\eta\E\frY_{\absv{\bsm}+n}(B_{o}\mr{M}^{\wt{w}_{1}})
			+\frac{1}{N}\E\frY_{\absv{\bsm}+n-1}(B_{o}\mr{M}^{\wt{w}_{1}}),
	\end{split}\end{align}
	where we abbreviated $P_{n}\equiv P_{n}(\wt{w}_{1},\ldots,\wt{w}_{n};\bsF)$. The exact forms of the expressions $\frY$ on right-hand sides do not depend on $B_{d}$($B_{o}$, resp.) except that they are supported on diagonal (off-diagonal, resp.) blocks.
	\end{lem}

While we can formally recover \eqref{eq:decoup_off} from \eqref{eq:decoup_diag} upon plugging in $B_{d}=B_{o}$, we separated the two identities since $\frY$ on the right-hand side depends on the support of $B_{d}$ or $B_{o}$.

As an immediate consequence of Lemma \ref{lem:BP_rough}, we have 
\beq\label{eq:XY_rough}
\absv{\E\frX_{n}}\lesssim N^{\epsilon}\Psi^{n},\qquad
\absv{\E\frY_{n}(B_{d})}\lesssim N^{\epsilon}\norm{B_{d}}\Psi^{n},\qquad \text{and }\qquad \absv{\E\frY_{n}(B_{o})}\lesssim N^{\epsilon}\norm{B_{o}}\Psi^{n+1}. 
\eeq
In this regard, Lemma \ref{lem:decoup} decouples the expression on the left-hand side to the leading order, but does not provide a useful estimate in terms of a fully deterministic quantity. This in fact is the general theme of our proof; say, we never try to estimate $\E\bsQ_{\bsm}$ but rather claim that the contribution of each $\E\bsQ_{\bsm}$ cancels out algebraically.

Notice that $B_{d}\mr{M}^{w_{1}}$ and $B_{o}\mr{M}^{w_{1}}$ are supported respectively on off-diagonal and diagonal blocks. Thus, if necessary, one may iteratively apply Lemma \ref{lem:decoup} to gain an extra precision of $\Psi$ at each step. Indeed, this will be our main task in the next section; we start from the case when $\absv{\bsm}=0$ and $n=1$, more precisely
\beq
\E\Brkt{\frac{\dd \E H_{t}^{w}}{\dd t}G_{t}^{w}}\bsL^{(j)},
\eeq
and apply (a more explicit version of) Lemma \ref{lem:decoup} four times so that the resulting decoupled expression has the precision $\Psi^{5}$. Indeed, the result of the last three consecutive applications of Lemma~\ref{lem:decoup} can be combined into the following generic lemma. We prove the lemma at the end of this section.
\begin{lem}\label{lem:decoup_3_C}
	Let $n\in\N$ and $\bsm$ be a multi-index with $n+\absv{\bsm}=2$, and $B_{d}\in\C^{2N\times 2N}$ be supported on diagonal blocks. Then we have, uniformly over $(w_{1},\ldots,w_{k})$ in a compact set,
	\beq\begin{split}\label{eq:decoup3_C}
		\E\brkt{B_{d}P_{n}}\bsQ_{\bsm}\wt{\bsL}
		=&\brkt{2B_{d}\mr{M}^{\wt{w}_{1}}\mr{M}^{\wt{w}_{n}}F_{n}}\E\brkt{P_{n}}\bsQ_{\bsm}\wt{\bsL}	\\
		&+\left(\lone(n=1)\ii\eta_{t}\brkt{B_{d}(\mr{M}^{\wt{w}_{1}})^{2}F_{1}}+\lone(n=2)\frac{1}{N}\brkt{B_{d}\mr{M}^{\wt{w}_{1}}F_{1}\mr{M}^{\wt{w}_{2}}F_{2}}\right)\E\bsQ_{\bsm}\wt{\bsL}	\\
		&+\frac{1}{2}\brkt{B_{d}(\mr{M}^{0})^{4}F_{n}}\E\frX_{4}
		+O(N^{\epsilon}\norm{B_{d}}\Psi^{5}),
	\end{split}\eeq
	where $P_{n}=P_{n}(\wt{w}_{1},\ldots,\wt{w}_{n};\bsF)$ and the expression $\frX_{4}$ does not depend on $B_{d}$.
\end{lem}
Taking $B_{d}=I$, we can solve \eqref{eq:decoup3_C} for the unspecified expression $\E\frX_{4}$. Consequently, we have the following immediate corollary of Lemma \ref{lem:decoup_3_C}. We omit the proof and only point out that by \eqref{eq:i4sim1} we have
\beqs
\brkt{2(\mr{M}^{0})^{4}E_{i}}=\brkt{{\absv{\caA_{t}}^{-4}}}=I_{4}\sim 1.
\eeqs
\begin{lem}\label{lem:decoup_3_C_c}
	Under the same notations as in Lemma \ref{lem:decoup_3_C} , the following holds uniformly over $(w_{1},\ldots,w_{k})$ in a compact set:
	\beq\begin{split}\label{eq:decoup3_C_c}
		\E\brkt{B_{d}P_{n}}\bsQ_{\bsm}\wt{\bsL}
		=&\left(\brkt{2B_{d}\mr{M}^{\wt{w}_{1}}\mr{M}^{\wt{w}_{n}}F_{n}}-(\brkt{2\mr{M}^{\wt{w}_{1}}\mr{M}^{\wt{w}_{n}}F_{n}}-1)\frac{\brkt{2B_{d}(\mr{M}^{0})^{4}F_{n}}}{I_{4}}\right)\E\brkt{P_{n}}\bsQ_{\bsm}\wt{\bsL}	\\
		&+\lone(n=1)\left(\brkt{B_{d}(\mr{M}^{\wt{w}_{1}})^{2}F_{1}}-\brkt{(\mr{M}^{\wt{w}_{1}})^{2}F_{1}}\frac{\brkt{2B_{d}(\mr{M}^{0})^{4}F_{n}}}{I_{4}}\right)\ii\eta_{t}\E\bsQ_{\bsm}\wt{\bsL}\\
		&+\lone(n=2)\left(\brkt{B_{d}\mr{M}^{\wt{w}_{1}}F_{1}\mr{M}^{\wt{w}_{2}}F_{2}}-\brkt{\mr{M}^{\wt{w}_{1}}F_{1}\mr{M}^{\wt{w}_{2}}F_{2}}\frac{\brkt{2B_{d}(\mr{M}^{0})^{4}F_{n}}}{I_{4}}\right)\frac{1}{N}\E\bsQ_{\bsm}\wt{\bsL}	\\
		&+O(N^{\epsilon}\norm{B_{d}}\Psi^{5}).
	\end{split}\eeq
\end{lem}

\subsection{Proof of Proposition \ref{prop:dt}}\label{sec:comp_comput}
We first prove \eqref{eq:comp_logd_d} and then \eqref{eq:comp_sing_d}.

\begin{proof}[Proof of Proposition \ref{prop:dt}, \eqref{eq:comp_logd_d}]
	By \eqref{eq:sc1} and Stein's lemma, we obtain (recall $\mr{M}=(\E H)^{-1}$)
	\beq\begin{split}\label{eq:KG_1}
		&\E\Brkt{\frac{\dd \E H_{t}^{w_{j}}}{\dd t}G^{w_{j}}E_{1}}\bsL^{(j)}	\\
		=&\Brkt{\frac{\dd \E H_{t}^{w_{j}}}{\dd t}\mr{M}^{w_{j}}E_{1}}\E\bsL^{(j)}	
		+\ii\eta_{t}\E\Brkt{\frac{\dd \E H_{t}^{w_{j}}}{\dd t}\mr{M}^{w_{j}}G^{w_{j}}E_{1}}\bsL^{(j)}	\\
		&-\E\Brkt{\frac{\dd \E H_{t}^{w_{j}}}{\dd t}\mr{M}^{w_{j}}WG^{w_{j}}E_{1}}\bsL^{(j)}	\\
		=&\Brkt{B_{d}^{w_{j}}E_{1}}\E\bsL^{(j)}
		+\ii\eta_{t}\E\Brkt{B_{d}^{w_{j}}G^{w_{j}}E_{1}}\bsL^{(j)}	\\
		&+\E\brkt{G^{w_{j}}}\Brkt{B_{d}^{w_{j}}G^{w_{j}}E_{1}}\bsL^{(j)}	
		-\frac{1}{N}\sum_{\ell}^{(j)}\E\Brkt{B_{d}^{w_{j}}G^{w_{\ell}}E_{2}G^{w_{j}}E_{1}}\bsL^{(j,\ell)},
	\end{split}\eeq
	where we abbreviated
	\beqs
	B_{d}^{w_{j}}\deq \frac{\dd \E H_{t}^{w_{j}}}{\dd t}\mr{M}^{w_{j}}
	\eeqs
	and used for the last term that $B_{d}^{w_{j}}$ and $E_{i}$ commute. Then we directly apply Lemma \ref{lem:decoup} to the second term of \eqref{eq:KG_1} with the choices $B_{d}=B_{d}^{w_{j}}$, $P_{n}=P_{1}(w_{j};E_{1})$, and $\wt{\bsL}=\bsL^{(j)}$, so that
	\beq\begin{split}\label{eq:KG}
		&\E\Brkt{2\frac{\dd \E H_{t}^{w_{j}}}{\dd t}G^{w_{j}}E_{1}}\bsL^{(j)}	\\
		=&\Brkt{2B_{d}^{w_{j}}E_{1}}\E\bsL^{(j)}
		+\ii\eta\Brkt{2B_{d}^{w_{j}}(\mr{M}^{w_{j}})^{2}E_{1}}\E\brkt{G^{w_{j}}}\bsL^{(j)}	\\
		&+\E\brkt{G^{w_{j}}}\Brkt{2B_{d}^{w_{j}}G^{w_{j}}E_{1}}\bsL^{(j)}
		-\frac{1}{N}\sum_{\ell}^{(j)}\E\Brkt{2B_{d}^{w_{j}}G^{w_{\ell}}E_{2}G^{w_{j}}E_{1}}\bsL^{(j,\ell)}
		+O(N^{\epsilon}\Psi^{5}),
	\end{split}\eeq
	where we used $\eta=O(\Psi^{3})$ and $N=O(\Psi^{4})$. 
	
	Likewise, applying Lemma \ref{lem:decoup_3_C_c} to the remaining two undecoupled terms in \eqref{eq:KG}, we find that
	\beq\begin{split}\label{eq:cancel_0}
		&\E\Brkt{2\frac{\dd (\E H_{t}^{w_{j}}-\ii\eta_{t})}{\dd t}G^{w_{j}}E_{1}}\bsL^{(j)}	\\
		=&\Brkt{2B_{d}^{w_{j}}E_{1}}\E\bsL^{(j)} \\
		&+\left(\brkt{2B_{d}^{w_{j}}(\mr{M}^{w_{j}})^{2}E_{1}}-\frac{\brkt{2B_{d}^{w_{j}}(\mr{M}^{0})^{4}E_{1}}}{I_{4}}\left(\brkt{(\mr{M}^{w_{j}})^{2}}-1\right)\right)\E\brkt{G^{w_{j}}}^{2}\bsL^{(j)}	\\
		&+\left(\brkt{4B_{d}^{w_{j}}(\mr{M}^{0})^{2}E_{1}}-\frac{\brkt{2B_{d}^{w_{j}}(\mr{M}^{0})^{4}E_{1}}}{I_{4}}+\frac{\dd c_{t}}{\dd t}\frac{1}{c_{t}}\right)\ii\eta\E\brkt{G^{w_{j}}}\bsL^{(j)}	\\
		&+\sum_{\ell}^{(j)}\left(\brkt{2B_{d}^{w_{j}}\mr{M}^{w_{\ell}}\mr{M}^{w_{j}}E_{1}}-\frac{\brkt{2B_{d}^{w_{j}}(\mr{M}^{0})^{4}E_{1}}}{I_{4}}(\brkt{2\mr{M}^{w_{\ell}}\mr{M}^{w_{j}}E_{1}}-1)\right)\E\frac{\brkt{G^{w_{\ell}}E_{2}G^{w_{j}}E_{1}}}{N}\bsL^{(j,\ell)}	\\
		&+\sum_{\ell}^{(j)}\frac{1}{2N}\left(\brkt{2B_{d}^{w_{j}}(\mr{M}^{w_{j}})^{2}E_{1}}-\frac{\brkt{2B_{d}^{w_{j}}(\mr{M}^{0})^{4}E_{1}}}{I_{4}}\right)\E\bsL^{(j,\ell)}
		+O(N^{\epsilon}\Psi^{5}),
	\end{split}\eeq
	where we also used the definition $\eta_{t}=c_{t}^{-1}\eta_{1}$.
	
	In what follows, we prove that the real part of \eqref{eq:cancel_0} is harmonic up to an $O(\Psi^{6})$ error after summing over $j=1,\ldots,k$. We plug in the definition of $B_{d}^{w}$ to compute the deterministic prefactors. For the first term, we simply have (recall $\mr{M}=(\E H)^{-1}$)
	\beq
	\re\brkt{2B_{d}^{w_{j}}E_{1}}=\Brkt{\frac{\dd \E H_{t}^{w_{j}}}{\dd t}\mr{M}_{t}^{w_{j}}}=-\frac{\dd}{\dd t}\brkt{\log\absv{\mr{M}_{t}^{w_{j}}}},
	\eeq
	which is harmonic in $w_{j}$ by Lemma \ref{lem:determ}. Secondly, we have from \eqref{eq:determ_4} that
	\beq\begin{split}\label{eq:can1}
		\re\brkt{2B_{d}^{w}(\mr{M}^{w})^{2}E_{1}}
		=&\re\Brkt{2\frac{\dd \E H_{t}^{w}}{\dd t}(\mr{M}_{t}^{w})^{3}E_{1}}=\Brkt{\frac{\dd \E H_{t}^{w}}{\dd t}(\mr{M}_{t}^{w})^{3}}=-\frac{1}{2}\frac{\dd}{\dd t}\Brkt{(\mr{M}_{t}^{w})^{2}}	\\
		=&-\frac{1}{2}\frac{\dd I_{4}}{\dd t}I_{4}^{-1/2}\re[-N^{-1/2}w]+O(N^{-1}\log N).
	\end{split}\eeq
	Likewise we have
	\beq\label{eq:can2}
	\re\brkt{2B_{d}^{w_{j}}(\mr{M}^{0})^{4}E_{1}}=\Brkt{\frac{\dd \E H_{t}^{w}}{\dd t}\mr{M}^{w}_{t}(\mr{M}^{0})^{4}}=-\frac{1}{4}\frac{\dd}{\dd t}I_{4}+O(N^{-1/2}\log N).
	\eeq
	Plugging \eqref{eq:determ_4}, \eqref{eq:can1}, and \eqref{eq:can2}, we have
	\beq\label{eq:can3}\begin{split}
		\re\left[\brkt{2B_{d}^{w_{j}}(\mr{M}^{w_{j}})^{2}E_{1}}-\frac{\brkt{2B_{d}^{w_{j}}(\mr{M}^{0})^{4}E_{1}}}{I_{4}}\left(\brkt{(\mr{M}^{w_{j}})^{2}}-1\right)\right]=O(N^{-1}\log N)=O(\Psi^{4}),
	\end{split}\eeq
	so that the second term of \eqref{eq:cancel_0} is $O(\Psi^{6})$.
	For the third term of \eqref{eq:cancel_0}, we combine $c_{t}=I_{4}^{-1/4}$, \eqref{eq:can1} up to the leading order, and \eqref{eq:can2} to obtain
	\beq\begin{split}
		&\brkt{4B_{d}^{w_{j}}(\mr{M}^{0})^{2}E_{1}}-\frac{\brkt{2B_{d}^{w_{j}}(\mr{M}^{0})^{4}E_{1}}}{I_{4}}+\frac{\dd c_{t}}{\dd t}\frac{1}{c_{t}}	\\
		=&O(N^{-1/2})+\frac{1}{4 I_{4}}\frac{\dd I_{4}}{\dd t}-\frac{1}{4 I_{4}}\frac{\dd I_{4}}{\dd t}=O(N^{-1/2})=O(\Psi^{2}).
	\end{split}\eeq
	Thus we conclude that the third term of \eqref{eq:cancel_0} is $O(\Psi^{6})$. The last term of \eqref{eq:cancel_0} is constant hence harmonic with respect to the variable $w_{j}$.
	
	Now it only remains to deal with the fourth term of \eqref{eq:cancel_0}. Using the facts that 
	\beq\begin{split}
		\Norm{B_{d}^{w_{j}}-B_{d}^{0}}=O(N^{-1/2})&=O(\Psi^{2}),	\\
		\brkt{2\mr{M}^{w_{j}}\mr{M}^{w_{\ell}}E_{1}}-1=O(N^{-1/2})&=O(\Psi^{2}),	\\
		\brkt{2B_{d}^{w_{j}}\mr{M}^{w_{\ell}}\mr{M}^{w_{j}}E_{1}}
		=\Brkt{2\mr{M}^{w_{j}}\frac{\dd \E H_{t}^{w_{j}}}{\dd t}\mr{M}^{w_{j}}\mr{M}^{w_{\ell}}E_{2}}&=-\Brkt{2\frac{\dd \mr{M}^{w_{j}}}{\dd t}\mr{M}^{w_{\ell}}E_{2}},
	\end{split}\eeq
	each summand in the fifth term of \eqref{eq:cancel_0} is equal to, up to an additive error of $O(\Psi^{6})$,
	\beq\label{eq:cancel_2pt}\begin{split}
		\re\left(-\Brkt{2\frac{\dd \mr{M}^{w_{j}}}{\dd t}\mr{M}^{w_{\ell}}E_{2}}-\frac{\brkt{2B_{d}^{0}(\mr{M}^{0})^{4}E_{1}}}{I_{4}}(\brkt{2\mr{M}^{w_{j}}\mr{M}^{w_{\ell}}E_{2}}-1)\right)\E\frac{\brkt{G^{w_{\ell}}E_{2}G^{w_{j}}E_{1}}}{N}\bsL^{(j,\ell)}.
	\end{split}\eeq

We take the sum of \eqref{eq:cancel_2pt} over $j\neq \ell$ in the following order; for each fixed pair of indices $j>\ell$, take the sum of \eqref{eq:cancel_2pt} with the same quantity with $j$ and $\ell$ interchanged, and then we take the sum over all pairs $j>\ell$. Consequently, the sum becomes
	\beqs\begin{split}\label{eq:cancel_2pt_0}
		&2\re\sum_{j>\ell}\left(-\Brkt{\frac{\dd (\mr{M}^{w_{j}}\mr{M}^{w_{\ell}})}{\dd t}E_{2}}-\frac{\brkt{B_{d}^{0}(\mr{M}^{0})^{4}}}{I_{4}}\left(\brkt{2\mr{M}^{w_{j}}\mr{M}^{w_{\ell}}E_{2}}-1\right)\right)\frac{\brkt{G^{w_{\ell}}E_{2}G^{w_{j}}E_{1}}}{N}\bsL^{(j,\ell)},
	\end{split}\eeqs
	where we used that
	\beq\begin{aligned}
		\ol{\Brkt{\frac{\dd \mr{M}^{w_{\ell}}}{\dd t}\mr{M}^{w_{j}}E_{2}}}&=\Brkt{\mr{M}^{w_{j}}\frac{\dd\mr{M}^{w_{\ell}}}{\dd t}E_{2}},	\quad&
		\ol{\brkt{\mr{M}^{w_{\ell}}\mr{M}^{w_{j}}E_{2}}}&=\brkt{\mr{M}^{w_{j}}\mr{M}^{w_{\ell}}E_{2}},		\\
		\ol{\brkt{B_{d}^{0}(\mr{M}^{0})^{4}E_{1}}}&=\brkt{B_{d}^{0}(\mr{M}^{0})^{4}E_{2}},	\quad &
		\ol{\brkt{G^{w_{j}}E_{2}G^{w_{\ell}}E_{1}}}&=\brkt{G^{w_{\ell}}E_{2}G^{w_{j}}E_{1}},	
	\end{aligned}\eeq
	due to $\mr{M}^{w}$ being Hermitian and $EG\check{E} =EG\adj \check{E}$ for $E\in\{E_{1},E_{2}\}$.
	
	Now we notice that the prefactor can be written as
	\beq\begin{split}\label{eq:cancel_2pt_1}
		\left(-\Brkt{\frac{\dd (\mr{M}^{w_{j}}\mr{M}^{w_{\ell}})}{\dd t}E_{2}}-\frac{\brkt{B_{d}^{0}(\mr{M}^{0})^{4}}}{I_{4}}\left(\brkt{2\mr{M}^{w_{j}}\mr{M}^{w_{\ell}}E_{2}}-1\right)\right)	\\
		=-\frac{1}{2}I_{4}^{1/2}\frac{\dd}{\dd t}\left[I_{4}^{-1/2}(\brkt{2\mr{M}^{w_{j}}\mr{M}^{w_{\ell}}E_{2}}-1)\right].
	\end{split}\eeq
	On the other hand, we have
	\beq\begin{split}\label{eq:cancel_2pt_2}
		\brkt{2\mr{M}^{w_{j}}\mr{M}^{w_{\ell}}E_{1}}&=\brkt{2(\mr{M}^{0})^{2}}+\brkt{J^{w_{j}}(\mr{M}^{0})^{3}E_{2}}+\brkt{J^{w_{\ell}}(\mr{M}^{0})^{3}E_{1}}+O(N^{-1})	\\
		&=1+N^{-1/2}\left(\ol{\gamma}^{-1}\ol{w}_{j}\ol{I}_{3}+\gamma^{-1}w_{\ell}I_{3}\right)+O(N^{-1})		\\
		&=1-I_{4}^{1/2}N^{-1/2}(w_{j}+w_{\ell})+O(N^{-1}).
	\end{split}\eeq
	As before, the error remains $O(N^{-1}\log N)$ after taking the time derivative, hence we conclude
	\beq\label{eq:cancel_2pt_c}
	\frac{\dd}{\dd t}\left[I_{4}^{-1/2}(\brkt{2\mr{M}^{w_{j}}\mr{M}^{w_{\ell}}E_{2}}-1)\right]=O(N^{-1}\log N).
	\eeq
	Plugging this back into \eqref{eq:cancel_2pt_1} and then \eqref{eq:cancel_2pt_0}, we conclude that the contribution of the fifth term of \eqref{eq:cancel_0} is $O(N^{-1}\Psi^{2})=O(N^{-1-c})$.
\end{proof}

\begin{proof}[Proof of Proposition \ref{prop:dt}, \eqref{eq:comp_sing_d}]
	We first make the following two observations about \eqref{eq:cancel_0} in the previous proof for $j=1$. Firstly, upon close inspection, the $O(\Psi^{5})$ error in \eqref{eq:cancel_0} can be tracked down algebraically. Indeed, we may write the error as $\E\frY$, where $\frY$ is a finite sum of terms of the following form:
	\beq\label{eq:sing_err_1}
	C_{\bsw}N^{-k_{1}/2}\cdot \eta^{k_{2}}N^{-k_{3}}\frY_{n}(B).
	\eeq
	Here $k_{1}$ is either $0$ or $1$, $k_{2}$ and $k_{3}$ are nonnegative integers, $\norm{B}=O(1)$ and $C_{\bsw}=O(1)$ are deterministic and independent of $\eta$, and
	\beq
	2k_{1}+3k_{2}+4k_{3}+n\geq 5.
	\eeq
	Secondly, the proof of \eqref{eq:cancel_0} is insensitive to the choice $\eta_{1}= N^{-3/4-\delta}$, but only requires that $\eta_{1}\sim N^{-3/4-\delta}$.
	
	Differentiating \eqref{eq:comp_sing_d} in $\ii\eta$ (more precisely $\partial_{\ii\eta}=\ii^{-1}\partial_{\eta}$), we find that
	\begin{equation}
		\label{eq:newder}
		\frac{\dd}{\dd t}\frac{\E \langle G_t\rangle}{c_t}=\partial_{\ii\eta_{t}}\E \left\langle \frac{\dd (\E H_t^w-\ii\eta_t)}{\dd t} G_t \right\rangle.
	\end{equation}
	Our goal is to show that the right-hand side of \eqref{eq:newder} is $O(N^{1+\epsilon}\Psi^{6})$. 
	
	Since \eqref{eq:cancel_0} for $j=1$ remains valid over varying $\eta_{1}$, we may apply $\partial_{\ii\eta_{t}}$ so that
	\begin{equation}\label{eq:can_sing}
		\begin{split}
			\partial_{\eta_{t}}\E\Brkt{\frac{\dd (\E H_t^w-{\ii\eta_{t}})}{\dd t}G_t^w}	&=\left(\brkt{B_{d}^w(\mr{M}_{t}^w)^{2}}-\frac{\brkt{B_{d}^w(\mr{M}_{t}^{0})^{4}}}{I_{4}}\left(\brkt{(\mr{M}_{t}^w)^{2}}-1\right)\right)\partial_\eta\E\brkt{G_t^w}^2	\\
			&\quad+\left(\brkt{4B_{d}^w(\mr{M}_{t}^{0})^{2}}-\frac{\brkt{B_{d}^w(\mr{M}_{t}^{0})^{4}}}{I_{4}}+{\frac{\dd c_{t}}{\dd t}\frac{1}{c_{t}}}\right)\ii\partial_{\eta_{t}}\big[\eta\E\brkt{G_t^w}\big]+\partial_{\eta_{t}}\E\frY,
		\end{split}
	\end{equation}
	where we used that $B_d^w$ does not depend on $\eta$, and $\frY$ is a finite sum of terms of the form given in \eqref{eq:sing_err_1}. Then, notice that the derivative $\partial_{\eta_{t}}\frY$ is a finite sum of
	\beq
	C_{\bsw}N^{-k_{1}/2}\eta^{k_{2}-1}N^{-k_{3}}\frY_{n}(B),\qquad
	C_{\bsw}N^{-k_{1}/2}\eta^{k_{2}}N^{-k_{3}+1}\frY_{n+1}(B).
	\eeq
	This in turn implies, by Lemma \ref{lem:BP_rough}, that
	\beq
	\partial_{\eta_{t}}\E\frY =O(N^{1+\epsilon}\Psi^{6})=O(N^{-c}\Psi).
	\eeq
	Finally, the result follows from substituting \eqref{eq:can2} and \eqref{eq:can3} into \eqref{eq:can_sing}.
\end{proof}
\subsection{Proof of Lemmas \ref{lem:BP_rough}--\ref{lem:decoup_3_C}}
Throughout this section, we use the following shorthand notation; for a given $\bsw\in\C^{n}$ and $\bsF\in\{E_{1},E_{2}\}^{n}$, we write
\beq\label{eq:tP}
\wt{P}_{n}\equiv \wt{P}_{n}(\bsw,\bsF)\deq N^{-n+1}G^{w_{1}}F_{1}\cdots G^{w_{n}}F_{n}\equiv P_{n}+\lone(n=1)\mr{M}^{w_{1}}F_{1},	
\eeq
and for a subset $I\subset\bbrktt{N}$ write 
\beq\label{eq:P_I}\begin{split}
	P_{I}\deq P_{\absv{I}}(\bsw_{I};\bsF_{i}), \quad \wt{P}_{I}&\deq \wt{P}_{\absv{I}}(\bsw_{I};\bsF_{I}),\qquad \bsw_{I}=(w_{i})_{i\in I},\quad \bsF_{I}=(F_{i})_{i\in I}.
\end{split}\eeq
Since $\brkt{\mr{M}^{w}F}=0$, we have $\brkt{P_{n}}=\brkt{\wt{P}_{n}}$ and similarly $\caS[P_{n}]=\caS[\wt{P}_{n}]$.
\begin{proof}[Proof of Lemma \ref{lem:BP_rough}]
	For $n=1$, \eqref{eq:BP_lead} is an immediate consequence of the local law. To prove \eqref{eq:BP_lead} for $n\geq 2$, we write (with the convention $\wt{P}_{0}=N$)
	\beq\begin{split}
		\absv{\brkt{BP_{n}(\bsw;\bsF)}}=\frac{1}{N^{2}}\absv{\brkt{G^{w_{n}} F_{n}B\wt{P}_{\bbrktt{n-2}}G^{w_{n-1}}F_{n-1}}}	\\
		\leq \frac{1}{N}\norm{B\wt{P}_{\bbrktt{n-2}}}\frac{\brkt{G^{w_{n}}}^{2}+\brkt{G^{w_{n-1}}}^{2}}{N},
	\end{split}\eeq
	where we used the general inequality
	\beq
	\absv{\brkt{B_{1}B_{2}B_{3}}}\leq \norm{B_{2}}(\brkt{\absv{B_{1}}^{2}}+\brkt{\absv{B_{2}}^{2}}),\qquad B_{1},B_{2},B_{3}\in\C^{2N\times 2N}.
	\eeq
	Then by $\norm{\wt{P}_{n}}\leq (N\eta)^{-n+1}\eta^{-1}$, we get
	\beq
	\frac{1}{N}\norm{\brkt{B\wt{P}_{\bbrktt{n-2}}}}\leq \frac{\norm{B}}{(N\eta)^{n-2}},
	\eeq
	so that by Ward identity and local laws Proposition \ref{prop:time dependent local law} it follows that
	\beq
	\absv{\brkt{BP_{n}(\bsw;\bsF)}}\leq \norm{B}\frac{\im\brkt{G^{w_{n}}}+\im\brkt{G^{w_{n-1}}}}{(N\eta)^{n-1}}\prec \norm{B}\Psi^{n}.
	\eeq
\end{proof}
\begin{proof}[Proof of Lemma \ref{lem:decoup}]
	We only prove \eqref{eq:decoup_diag}, and \eqref{eq:decoup_off} can be proved analogously by changing the parities of few $F$'s. Expanding the first factor $G^{\wt{w}_{1}}$ in $\brkt{B_{d}P_{n}}$, we get
	\beq\begin{split}\label{eq:BP_1}
		\E\brkt{B_{d}P_{n}}\bsQ_{\bsm}\wt{\bsL}
		=&\E\left(-\brkt{B_{d}\mr{M}^{w_{1}}W\wt{P}_{n}}+\lone(n\geq 2)\frac{\brkt{B_{d}\mr{M}^{w_{1}}F_{1}\wt{P}_{\bbrktt{2,n}}}}{N}+\ii\eta_{t}\brkt{B_{d}\mr{M}^{w_{1}}\wt{P}_{n}}\right)\bsQ_{\bsm}\wt{\bsL}	\\
		=&-\E\brkt{B_{d}\mr{M}^{w_{1}}W\wt{P}_{n}}\bsQ_{\bsm}\wt{\bsL}\\
		&+\left(\lone(n=1)\ii\eta_{t}\brkt{B_{d}(\mr{M}^{\wt{w}_{1}})^{2}F_{1}}+\lone(n=2)\frac{1}{N}\brkt{B_{d}\mr{M}^{\wt{w}_{1}}F_{1}\mr{M}^{\wt{w}_{2}}F_{2}}\right)\E\bsQ_{\bsm}\wt{\bsL}	\\
		&+\ii\eta\E\frY_{\absv{\bsm}+n}(B_{d}\mr{M}^{\wt{w}_{1}})
		+\frac{1}{N}\E\frY_{\absv{\bsm}+n-1}(B_{d}\mr{M}^{\wt{w}_{1}}).
	\end{split}\eeq
	By Stein's lemma and Leibniz rule, we may write the first term of \eqref{eq:BP_1} as
	\beq\label{eq:BP_2}
	-\E\bsQ_{\bsm}\brkt{B_{d}\mr{M}^{\wt{w}_{1}}W\wt{P}_{n}}\wt{\bsL}=-\E\E_{\wt{W}}\brkt{B_{d}\mr{M}^{\wt{w}_{1}}\wt{W}\partial_{\wt{W}}[\wt{P}_{n}]}\bsQ_{\bsm}\wt{\bsL}-\E\E_{\wt{W}}\brkt{B_{d}\mr{M}^{\wt{w}_{1}}\wt{W}\wt{P}_{n}}\partial_{\wt{W}}[\bsQ_{\bsm}\wt{\bsL}].
	\eeq
	For the first term of \eqref{eq:BP_2}, we write
	\beq\label{eq:BP_2_0}
	-\partial_{\wt{W}}\wt{P}_{n}=\frac{1}{N}\sum_{i=1}^{n}\wt{P}_{\bbrktt{i-1}}G^{\wt{w}_{i}}\wt{W}\wt{P}_{\bbrktt{i,n}},
	\eeq
	so that
	\beq
	-\E\E_{\wt{W}}\brkt{B_{d}\mr{M}^{\wt{w}_{1}}\wt{W}\partial_{\wt{W}}[\wt{P}_{n}]}\bsQ_{\bsm}\wt{\bsL}
	=\frac{1}{N}\E\sum_{i=1}^{n}\Brkt{B_{d}\mr{M}^{\wt{w}_{1}}\caS[\wt{P}_{\bbrktt{i-1}}G^{\wt{w}_{i}}]\wt{P}_{\bbrktt{i,n}}}\bsQ_{\bsm}\wt{\bsL}.
	\eeq
	For each $i\in\bbrktt{1,n}$, we have
	\beq\label{eq:BP_2_1}\begin{split}
		&\frac{1}{N}\brkt{B_{d}\mr{M}^{\wt{w}_{1}}\caS[\wt{P}_{\bbrktt{i-1}}G^{\wt{w}_{i}}]\wt{P}_{\bbrktt{i,n}}}\bsQ_{\bsm}\wt{\bsL}
		=\brkt{2B_{d}\mr{M}^{\wt{w}_{1}}\wt{P}_{\bbrktt{i,n}}}\brkt{P_{i}(\wt{w}_{\bbrktt{i}};\bsF_{\bbrktt{i-1}},F_{n})}\bsQ_{\bsm}\wt{\bsL}	\\
		=&\lone(i=n)\brkt{2B_{d}\mr{M}^{\wt{w}_{1}}\mr{M}^{\wt{w}_{n}}}\brkt{P_{n}}\bsQ_{\bsm}\wt{\bsL}+\frac{1}{2}\frY_{\absv{\bsm}+n+1}(B_{d}\mr{M}^{\wt{w}_{1}}).
	\end{split}\eeq
	Summing \eqref{eq:BP_2_1} over $i\in\bbrktt{1,n}$ and then substituting into \eqref{eq:BP_2}, it only remains to prove
	\beq
	-\E_{\wt{W}}\brkt{B_{d}\mr{M}^{\wt{w}_{1}}\wt{W}\wt{P}_{n}\partial_{\wt{W}}[\bsQ_{\bsm}\wt{\bsL}]}=\frY_{\bsm+n+1}(B_{d}\mr{M}^{\wt{w}_{1}}).
	\eeq
	
	We first consider the derivative of $\bsQ_{\bsm}$. For each factor $\brkt{P_{m_{j}}}\equiv \brkt{P_{m_{j}}(\wt{\bsw}',\bsF')}$ in $\bsQ_{\bsm}$, we use \eqref{eq:quad} and \eqref{eq:BP_2_0} to write
	\beq\begin{split}\label{eq:BP_2_2}
		&-\E_{\wt{W}}\brkt{B_{d}\mr{M}^{\wt{w}_{1}}\wt{W}\wt{P}_{n}}\partial_{\wt{W}}[\brkt{P_{m_{j}}(\wt{\bsw}';\bsF')}]\bsQ^{(j)}\wt{\bsL}	\\
		=&\sum_{i}\E_{\wt{W}}\brkt{B_{d}\mr{M}^{\wt{w}_{1}}\wt{W}\wt{P}_{n}}\brkt{\wt{W}\wt{P}_{m_{j}}(\wt{\bsw}_{i}';\bsF'_{i})G^{\wt{w}_{i}'}}\bsQ^{(j)}\wt{\bsL}	\\
		=&\frac{1}{2}\sum_{i}\brkt{B_{d}\mr{M}^{\wt{w}_{1}}P_{m_{j}+n+1}(\wt{\bsw}_{i}',\wt{w}_{i}',\wt{\bsw};\bsF_{i}',\check{F}_{n},\bsF)}\bsQ^{(j)}\wt{\bsL}	\\
		=&\frac{1}{2}\frY_{\bsm+n+1}(B_{d}\mr{M}^{\wt{w}_{1}}),
	\end{split}\eeq
	for some cyclic rearrangements $\wt{\bsw}'_{i}$ and $\bsF'_{i}$ of $\wt{\bsw}'$ and $\bsF'$ depending on $i\in\bbrktt{m_{j}}$. Summing \eqref{eq:BP_2_2} over $j$ proves
	\beq
	-\E_{\wt{W}}\brkt{B_{d}\mr{M}^{\wt{w}_{1}}\wt{W}\wt{P}_{n}}\partial_{\wt{W}}[\bsQ_{\bsm}]\wt{\bsL}=\frac{1}{2}\frY_{\bsm+n+1}(B_{d}\mr{M}^{\wt{w}_{1}}).
	\eeq
	Similarly, the derivative of $\wt{\bsL}$ can be handled using
	\beq
	\partial_{\wt{W}}L(w)=\Tr G^{w}\wt{W}.
	\eeq
	This completes the proof of Lemma \ref{lem:decoup}.
\end{proof}

\begin{proof}[Proof of Lemma \ref{lem:decoup_3_C}]
	First of all, we replace $B_{d}$ by $F_{n}B_{d}$ using
	\beq
	\brkt{B_{d}P_{n}}=\brkt{B_{d}P_{n}F_{n}}=\brkt{F_{n}B_{d}P_{n}}.
	\eeq
	Then we apply \eqref{eq:decoup_diag} and estimate the last line therein using Lemma \ref{lem:BP_rough}, $\eta\leq\Psi^{3}$, and $N^{-1}\leq \Psi^{4}$;
	\beq
	\eta\absv{\frY_{2}(F_{n}B_{d}\mr{M}^{\wt{w}_{1}})}+N^{-1}\absv{\frY_{1}(F_{n}B_{d}\mr{M}^{\wt{w}_{1}})}\prec\norm{B_{d}}\Psi^{5}.
	\eeq 
	Thus it suffices to express the the third line of \eqref{eq:decoup_diag} as $\frX_{4}$, i.e. to prove that
	\beq
	\E\frY_{3}(F_{n}B_{d}\mr{M}^{\wt{w}_{1}})=\brkt{F_{n}B_{d}(\mr{M}^{0})^{4}}\frX_{4}.
	\eeq
	
	Since $F_{n}B_{d}\mr{M}^{\wt{w}_{1}}$ is supported on off-diagonal blocks, we may apply \eqref{eq:decoup_off} to each summand in $\E\frY_{3}(F_{n}B_{d}\mr{M}^{\wt{w}_{1}})$ to write it as a finite sum of
	\beq
	\E\frY_{4}(F_{n}B_{d}\mr{M}^{\wt{w}_{1}}\mr{M}^{w'_{1}})+O(N^{\epsilon}\norm{B_{d}}\Psi^{6}),
	\eeq
	where $w'_{1}$ may vary in $\{w_{1},\ldots,w_{k}\}$. Likewise, applying \eqref{eq:decoup_diag} to each summand in $\E\frY_{4}(F_{n}B_{d}\mr{M}^{\wt{w}_{1}}\mr{M}^{w'_{1}})$, we conclude that 
	$\E \frY_{3}(F_{n}B_{d}\mr{M}^{\wt{w}_{1}})$ is a finite sum of the expressions 
	\beq\label{eq:decoup_3_1}
	\brkt{2F_{n}B_{d}\mr{M}^{\wt{w}_{1}}\mr{M}^{w_{1}'}\mr{M}^{w'_{2}}\mr{M}^{w'_{3}}E}\E\frX_{4}+O(N^{\epsilon}\norm{B_{d}}\Psi^{6}),
	\eeq
	where $E$ can be either $E_{1}$ or $E_{2}$ and $w'_{1},w'_{2},w'_{3}$ may vary in $\{w_{1},\ldots,w_{k}\}$. Furthermore, Lemma \ref{lem:decoup} implies that the exact form of the sum is determined by the expressions $P_{n}$, $\bsQ_{\bsm}$, and $\wt{\bsL}$ on the left-hand side of \eqref{eq:decoup3_C} that we started from.
	
	Now it only remains to notice that \eqref{eq:decoup_3_1} can be further expressed as
	\beq
	\brkt{2F_{n}B_{d}\mr{M}^{\wt{w}_{1}}\mr{M}^{w_{1}'}\mr{M}^{w'_{2}}\mr{M}^{w'_{3}}E}\E\frX_{4}=\lone(E=F_{n})\brkt{2B_{d}(\mr{M}^{0})^{4}F_{n}}+O(N^{\epsilon}\Psi^{6}),
	\eeq
	where we used $\norm{\mr{M}^{w}-\mr{M}^{0}}=O(N^{-1/2})=O(\Psi^{2})$ and $\absv{\frX_{4}}\prec\Psi^{4}$.
	
\end{proof}

\section{Proof of Proposition \ref{prop:comp_G}, the real case}\label{sec:real}
The main goal of this section is to prove Proposition \ref{prop:comp_G} when $\bbF=\bbR$.
Since we only consider deformed real Ginibre ensembles, we drop the superscripts $\mathrm{Gin}$ and $\mathrm{Gin}(\bbF)$.

The fundamental difference in the real case is that the Hermitization
\beq
W\deq \begin{pmatrix} 0 & X \\ X\tp & 0 \end{pmatrix}
\eeq
is now a real symmetric random matrix. Thus, each application of Stein's lemma results in an additional term as follows (recall the definitions of $E_{1},E_{2}$ from \eqref{eq:def_E1}):
\beq\begin{aligned}\label{eq:S_R}
	\caS_{\R}[B]\deq\E_{\wt{W}}\wt{W}B\wt{W}=&\caS[B]+\frac{1}{N}\sum_{E=E_{1},E_{2}}EB\tp\check{E},	\\
	\E_{\wt{W}}\wt{W}\Tr[\wt{W}B]=&\frac{1}{N}\sum_{E=E_{1},E_{2}}E(B+B\tp)\check{E}.
\end{aligned}\eeq
Notice that the additional terms involve the transpose, so that we often encounter $G\tp$ and its deterministic counterparts $M\tp$ and $\mr{M}\tp$. In particular, we will shortly see that the size of the quantity
\beq
\theta_{0}\deq 1-\Brkt{\frac{1}{(A-z_{0})(A-z_{0})\tp}}\in\C
\eeq
determines the universality class of our statistics. Notice from Assumption \ref{assump:edge} (i) that $\absv{1-\theta_{0}}\leq 1$, and that if $(A-z_{0})\tp$ is equal to $(A-z_{0})\adj$, i.e. if $A-z_{0}$ is real, then $\theta_{0}$ is identically zero.
Furthermore, writing $\caA_{0}\deq (A-z_{0})$, (recall from \eqref{eq:def_REIM} that $\Im[\caA_{0}]$ stands for entrywise imaginary part of $\caA_{0}$)
\beq\begin{split}\label{eq:theta_asymp}
	0\leq \re \theta_{0} =\frac{1}{2}\Brkt{\frac{1}{\absv{\caA_{0}}^{2}}-\frac{1}{\caA_{0}\caA_{0}\tp}-\frac{1}{\caA_{0}\adj\ol{\caA}_{0}}+\frac{1}{\absv{\ol{\caA}_{0}}^{2}}}
	=\frac{1}{2}\Brkt{\Absv{\frac{1}{\caA_{0}}-\frac{1}{\ol{\caA_{0}}}}^{2}}
	\sim\brkt{\absv{\Im[\caA_{0}]}^{2}},	\\
	\im\theta_{0}= \Brkt{\frac{1}{\caA_{0}\caA_{0}\tp}(\Im[\caA_{0}]\caA_{0}\tp+\ol{\caA_{0}}\Im[\caA_{0}]\tp)\frac{1}{\ol{\caA}_{0}\caA_{0}\adj}}
	=2\Brkt{\Im[\caA_{0}]\caA_{0}\adj\frac{1}{\absv{\caA_{0}}^{4}}}
	+O(\brkt{\absv{\Im[\caA_{0}]}^{2}}),
\end{split}\eeq
where we used $\norm{\caA_{0}^{-1}},\norm{\caA_{0}}=O(1)$ and Cauchy-Schwarz to get asymptotic relations. In particular, since $A$ is real, the asymptotics \eqref{eq:theta_asymp} reads 
\beq\label{eq:theta_asymp_R}
\re \theta_{0}\sim\absv{\im z_{0}}^{2},\qquad \im \theta_{0}=-2I_{3}\im z_{0}+O(\absv{\im z_{0}}^{2})
\eeq
as $\absv{\im z_{0}}\to 0$.
With this in mind, our proof of Proposition \ref{prop:comp_G} in the real case can be divided into two regimes depending on the size of $\absv{\theta_{0}}\sim\absv{\im z_{0}}$ in which we construct a path $\caA_{t}$ differently. When $\absv{\im z_{0}}\geq N^{-1/2+\frc}$, we construct a path $\caA_{t}$ that makes the corresponding $\theta$ always much larger than $N^{-1/2}$. When $\absv{\im z_{0}}\leq N^{-1/2+\frc}$, our choice of $\caA_{t}$ keeps the value of $I_{4}^{-1/2}\theta$ almost constant. We formalize this intuition in the next assumption.
\begin{assump}\label{assu:real_flow}
	Assume that $\caA_{t}$ is a path satisfying Assumption \ref{assump:path}.
	Defining 
	\beq\label{eq:def_theta}
	\theta_{t}\deq 1-\Brkt{\frac{1}{\caA_{t}\caA_{t}\tp}},
	\eeq
	we further assume that either of the following holds.
	\begin{itemize}
		\item[(A1)] for some small constants $\frc_{1},\frc_{2}>0$, $\absv{\theta_{t}}\geq N^{-1/2+\frc_{1}}$ and $\absv{\theta_{t}^{-2}{\dd\theta_{t}}/{\dd t}}\leq N^{1/2-\frc_{2}}$ for all $t\in[0,1]$.
		\item[(A2)] for some small constants $\frc_{1},\frc_{2}>0$, $\absv{\theta_{t}}\leq N^{-1/2+\frc_{1}}$, $\absv{\dd[I_{4}(t)^{-1/2}\theta_{t}]/{\dd t}}\leq N^{-1/2-\frc_{2}}$, and $\norm{\dd\Im\caA_{t}/\dd t}\leq N^{-\frc_{2}}$ for all $t\in[0,1]$.
	\end{itemize}
\end{assump}
We next construct such a path $\caA_{t}$ in each of the two cases.
\begin{prop}\label{prop:real_path}
	Suppose that $A$ is real and $z_{0}\in\C$ satisfies Assumption \ref{assump:edge}. Recall the small constant $\frc>0$ from Assumption \ref{assump:edge}.
	\begin{itemize}
		\item[(i)] If $\absv{\im z_{0}}\geq N^{-1/2+\frc}$ i.e. if $z_{0}$ satisfies Assumption \ref{assump:edge} (iii$\R$), then there exists a path $\caA_{t}$ satisfying Assumption \ref{assu:real_flow} (A1).
		\item[(ii)] If $\absv{\im z_{0}}\leq N^{-1/2+\frc}$ i.e. if $z_{0}$ satisfies Assumption \ref{assump:edge} (iii$\R'$), then there exists a path $\caA_{t}$ satisfying Assumption \ref{assu:real_flow} (A2). 
	\end{itemize}
	Moreover, the constants $C_{5}$ in \eqref{eq:A_bdd} and $\frc_{1},\frc_{2}$ in (A1) -- (A2) depend only on those in Assumption \ref{assump:edge}.
\end{prop} 

\begin{figure}
	\centering
\begin{subfigure}{0.24\textwidth}
	\includegraphics[width=\linewidth]{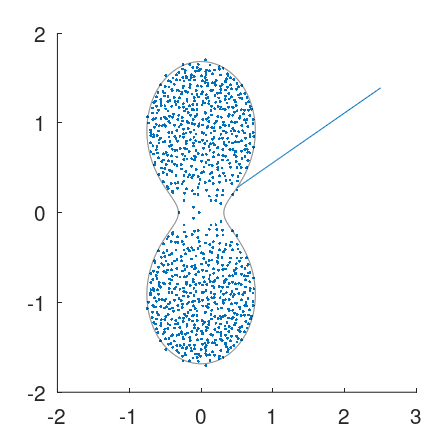}
\end{subfigure}
\begin{subfigure}{0.24\textwidth}
	\includegraphics[width=\linewidth]{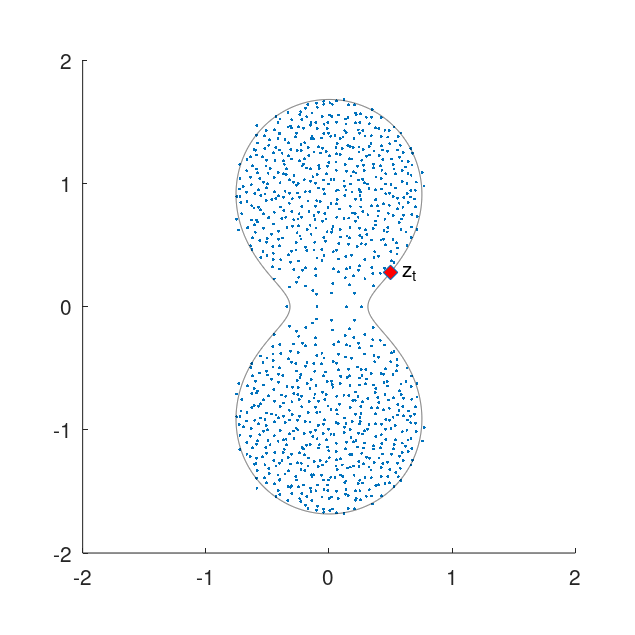}
\end{subfigure}
\begin{subfigure}{0.24\textwidth}
	\includegraphics[width=\linewidth]{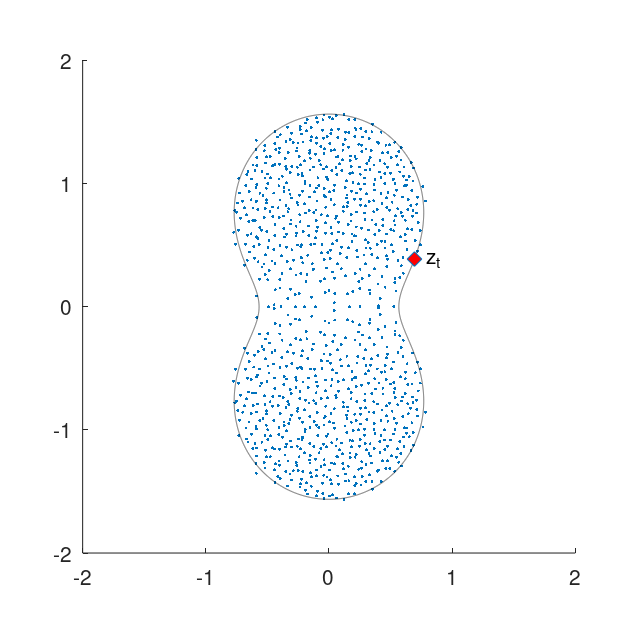}
\end{subfigure}	
\begin{subfigure}{0.24\textwidth}
	\includegraphics[width=\linewidth]{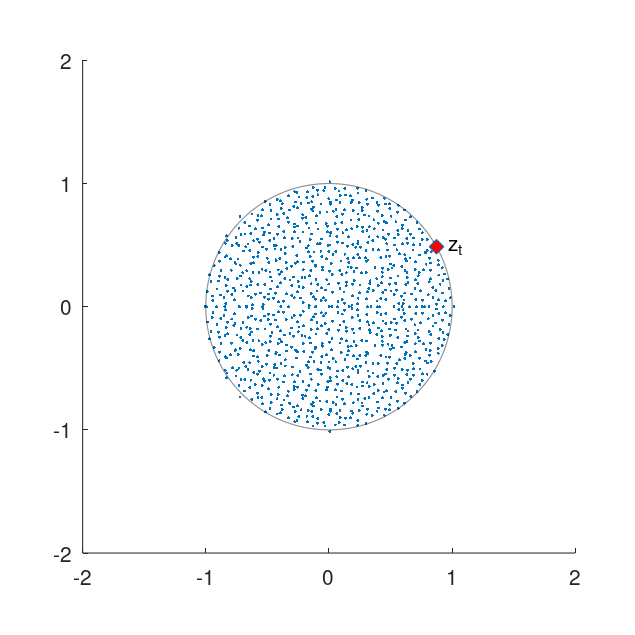}
\end{subfigure}
\caption{Case (A2): Leftmost picture shows the path $\wh{z}_{t}$, the other pictures show the eigenvalues of $\caA_{t}+X$
and the evolution of the base point $z_t = \brkt{\absv{A-\wh{z}_{t}}^{-2}}^{1/2} \wh{z}_{t}$.}\label{fig:real_path}
\end{figure}
\begin{proof}
	In what follows, we construct two different paths $\wh{z}_{t}\in\C$ depending on $\im z_{0}$, i.e. the cases (i) and (ii), and show that the matrix-valued path 
	\beq\label{eq:A_t}
	\caA_{t}\deq \brkt{\absv{A-\wh{z}_{t}}^{-2}}^{1/2}(A-\wh{z}_{t})
	\eeq
	satisfies Assumption \ref{assu:real_flow}, see Figure~\ref{fig:real_path}. 
	
	We start with the proof of (i), and without loss of generality assume $\im z_{0}>0$. Using Assumption \ref{assump:edge} (iii$\R$), we take a smooth path $[0,1]\ni t\to \wt{z}_{t}\in\caD_{C'}\cap\C_{+}$ such that $\wt{z}_{0}=z_{0}$ and $\wt{z}_{1}=Cz_{1}$ where $C>0$ is a large enough constant. Then we notice that, for a small positive constant $\tau$, the path $\wt{z}_{t}+\ii\tau$ is in $\caD_{2C'}$. We then define the path $\wh{z}_{t}$ as follows:
	\beq
	\wh{z}_{t}\deq \begin{cases} \re z_{0}+\ii (\tau+\im z_{0})\left(\frac{\im z_{0}}{\tau+\im z_{0}}\right)^{1-3t}, & \qquad t\in[0,1/3];	\\
		\wt{z}_{3t-2}+\ii\tau	& \qquad t\in[2/3,2/3];	\\
		\frac{2}{3(1-t)}\wt{z}_{1}+\ii\tau	& \qquad t\in[2/3,1].
	\end{cases}
	\eeq
	The path $\wh{z}_{t}$ here is constructed in a similar way to Proposition \ref{prop:z}, except that in the first part we lift (translate) the path by $\ii\tau$. Since $\tau$ is small but fixed, this guarantees that $|\Im[A-\wh{z}_{t}]|$ is bounded from below afterwards. 
	We omit verifying Assumption \ref{assump:path} since it is analogous to the complex case, but only remark that the last inequality in \eqref{eq:A_bdd}, for $t\in[0,1/3]$, is due to
	\beq\label{eq:dt_log}
	\frac{\dd}{\dd t}a^{t}=a^{t}\log a,\qquad a>0.
	\eeq
	Note that by \eqref{eq:theta_asymp_R} we have $\absv{\theta_{0}}\sim \absv{\im z_{0}}\geq N^{-1/2+\frc}$, henceforth we verify that $\caA_{t}$ satisfies (A1).
	
	By the definition of $\caA_{t}$, we have
	\beq
	\theta_{t}=1-\brkt{\absv{A-\wh{z}_{t}}^{-2}}^{-1}\Brkt{\frac{1}{(A-\wh{z}_{t})(A-\wh{z}_{t})\tp}}.
	\eeq
	Thus it immediately follows $\theta_{t}=2+O(1-t)$ when $t$ is close enough to $1$. Thus we restrict to $t\in[0,1-\delta)$, in which case $\absv{\wh{z}_{t}}$ remains bounded so that
	\beq\label{eq:theta_C}
	\re \theta_{t}\sim \absv{\im \wh{z}_{t}}^{2}, \qquad 
	\im \theta_{t}= 2\Brkt{(A-\wh{z}_{t})\adj\frac{1}{\absv{A-\wh{z}_{t}}^{4}}}\im \wh{z}_{t}+O(\absv{\im \wh{z}_{t}}^{2})
	\eeq
	by the same argument as in \eqref{eq:theta_asymp_R}.
	
	Explicitly computing $\dd\theta_{t}/\dd t$, we find that
	\beq\begin{split}
		&\brkt{\absv{A-\wh{z}_{t}}^{-2}}\frac{\dd \theta_{t}}{\dd t}	\\
		=&(1-\theta_{t})\frac{\dd}{\dd t}\Brkt{\frac{1}{\absv{A-\wh{z}_{t}}^{2}}}
		-\frac{\dd}{\dd t}\Brkt{\frac{1}{(A-\wh{z}_{t})(A-\wh{z}_{t})\tp}}\\
		=&2(1-\theta_{t})\re\left[\frac{\dd\wh{z}_{t}}{\dd t}\Brkt{\frac{1}{(A-\wh{z}_{t})^{2}(A-\wh{z}_{t})\adj}}\right] -2\frac{\dd\wh{z}_{t}}{\dd t}\Brkt{\frac{1}{(A-\wh{z}_{t})^{2}(A-\wh{z}_{t})\tp}},
	\end{split}\eeq
	so that writing down the real and imaginary parts we get
	\beq
	\Absv{\frac{\dd\theta_{t}}{\dd t}}\lesssim\Absv{\frac{\dd \im \wh{z}_{t}}{\dd t}}+(\absv{\theta}+\im \wh{z}_{t})\Absv{\frac{\dd\re\wh{z}_{t}}{\dd t}}.
	\eeq
	 From $\absv{\theta_{t}}\sim \absv{\im \wh{z}_{t}}$ due to \eqref{eq:theta_C} together with $\absv{\dd\re \wh{z}_{t}/\dd t}=O(1)$, we thus conclude
	\beq\label{eq:dtheta_C}
	\Absv{\frac{1}{\theta_{t}^{2}}\frac{\dd \theta_{t}}{\dd t}}\lesssim \Absv{\frac{1}{(\im \wh{z}_{t})^{2}}\frac{\dd\im \wh{z}_{t}}{\dd t}}+\frac{1}{\absv{\im \wh{z}_{t}}}.
	\eeq
	Hence (A1) immediately follows from \eqref{eq:theta_C} and \eqref{eq:dtheta_C} since $\im \wh{z}_{t}\geq \tau\vee \im z_{0}$ and
	\beq
	\Absv{\frac{1}{\im \wh{z}_{t}}\frac{\dd\im\wh{z}_{t}}{\dd t}}\sim \begin{cases}1+\absv{\log \im z_{0}}, & t\in[0,1/3], \\
		1, & t\in[1/3,1],
	\end{cases}
	\eeq
	where we used \eqref{eq:dt_log} for $t\in[0,1/3]$.
	
	We move on to the proof of (ii), and assume without loss of generality $[\re z_{0},\infty)$ is in $\caD_{C'}$. Here we will build a different path $\wh{z}_{t}$ as follows. Define
	\beq
	\re \wh{z}_{t}=\re z_{0}+\frac{t}{1-t},
	\eeq
	and we take $\im \wh{z}_{t}=\frac{f(\re \wh{z}_{t})}{f(\re z_{0})}\im z_{0}$ where the function $f:[\re z_{0},\infty)\to \R\setminus\{0\}$ is defined by
	\beq
	f(x)\deq\brkt{\absv{A-x}^{-4}}^{-1/2}\brkt{\absv{A-x}^{-2}}^{-1}\Brkt{\frac{1}{(A-x)^{2}(A-x)\tp}}.
	\eeq
	We then define our matrix-valued path $\caA_{t}$ again by the same formula \eqref{eq:A_t}. We can easily verify Assumption \ref{assump:path} for this path, modulo that the end point is not exactly $-z_{1}$ but is another point $\caA_{1}=-z_{1}'\in\bbS^{1}$ with
	\beq\label{eq:z_1}
	\im z_{1}'=\absv{\gamma}\im z_{0}+o(N^{-1/2})=\im z_{1}+o(N^{-1/2}).
	\eeq
	We then concatenate to $\caA_{t}$ the (shorter) arc along the unit circle connecting $-z_{1}'$ to $-z_{1}$, which is of length $o(N^{-1/2})$. Henceforth we focus on checking (A2) and \eqref{eq:z_1}.
	
	Notice that $\caA_{t}$ can be written as
	\beq\label{eq:At_R}
	\caA_{t}=\Re \caA_{t}+\ii k_{t}\im z_{0},
	\eeq
	where $k_{t}>0$ is defined by
	\beq
	k_{t}\deq\brkt{\absv{A-\wh{z}_{t}}^{-2}}^{1/2}\frac{f(\re \wh{z}_{t})}{f(\re z_{0})} =\frac{1}{f(\re z_{0})}\brkt{\absv{\Re\caA_{t}}^{-4}}^{-1/2}\brkt{\absv{\Re \caA_{t}}^{-2}}^{-1}\Brkt{\frac{1}{(\Re \caA_{t})^{2}\Re\caA_{t}\tp}}.
	\eeq
	By direct computation, we easily find that
	\beq\label{eq:dk/dt}
	\begin{aligned}
		&\norm{\Re\caA_{t}}=O(1),&\qquad 
		&\Norm{\frac{1}{\Re\caA_{t}}}=O(1), &\qquad 
		&\Norm{\frac{\dd \Re\caA_{t}}{\dd t}}=O(1), \\
		&\Absv{\Brkt{\frac{1}{(\Re\caA_{t})^{2}\Re\caA_{t}\tp}}}\sim 1,  &
		&\absv{k_{t}}\sim 1,&
		&\Absv{\frac{\dd k_{t}}{\dd t}}=O(1).
	\end{aligned}\eeq
	
	Since $\caA_{t}$ is always of the form \eqref{eq:At_R}, the asymptotics \eqref{eq:theta_asymp_R} remains true, so that $\im \theta_{t}\sim \im z_{0}$ and  $\re \theta_{t}=O(\absv{\im z_{0}}^{2}).$
	This immediately proves $\absv{\theta_{t}}\leq N^{-1/2+\frc_{1}}$. For the second assertion in (A2), we recall from the definition of $\theta_{t}$ that
	\beq
	I_{4}^{-1/2}\theta_{t}
	=I_{4}^{-1/2}\Brkt{\frac{1}{\absv{\caA_{t}}^{2}}-\frac{1}{\caA_{t}\caA_{t}\tp}}
	=\ii I_{4}^{-1/2}k_{t}\im z_{0}\Brkt{\frac{1}{\absv{\caA_{t}}^{2}\caA_{t}\tp}}.
	\eeq
	By the definition of $k_{t}$, we have 
	\beq
	I_{4}^{-1/2}\Brkt{\frac{k_{t}}{\absv{\caA_{t}}^{2}\caA_{t}\tp}}=\frac{\brkt{\absv{\caA_{t}}^{-2}}^{2}}{\brkt{\absv{\Re\caA_{t}}^{-2}}^{2}}\frac{\Brkt{(\Re\caA_{t}\tp)^{-1}\absv{\Re\caA_{t}}^{-2}}}{\brkt{(\caA_{t}\tp)^{-1}\absv{\caA_{t}}^{-2}}}\frac{\Brkt{\absv{\Re \caA_{t}}^{-4}}^{1/2}}{\Brkt{\absv{\caA_{t}}^{-4}}^{1/2}}.
	\eeq
	Thus, after several applications of the resolvent identity 
	\beq
	\frac{1}{\Re\caA_{t}}-\frac{1}{\caA_{t}}=\ii\im z_{0}\frac{1}{\Re[\caA_{t}]\caA_{t}}
	\eeq
	and \eqref{eq:dk/dt}, we arrive at
	\beq
	\frac{\dd}{\dd t}\Brkt{\frac{k_{t}}{\absv{\caA_{t}}^{2}\caA_{t}\tp}}=O(\absv{\im z_{0}}).
	\eeq
	Therefore we conclude 
	\beq
	\Absv{\frac{\dd\theta_{t}}{\dd t}}=O(\absv{\im z_{0}}^{2})=O(N^{-1+2\frc})=O(N^{-1/2-\frc_{2}}),
	\eeq 
	as desired. The last assertion of (A2) immediately follows from \eqref{eq:dk/dt} as
	\beq
	\Norm{\frac{\dd\Im \caA_{t}}{\dd t}}=\absv{\im z_{0}}\Absv{\frac{\dd k_{t}}{\dd t}}\lesssim \im z_{0}.
	\eeq
	To show \eqref{eq:z_1}, we only need to notice that
	\beq\begin{split}
		&I_{4}(0)^{-1/2}\theta_{0}=2\ii \im z_{0}I_{4}^{-1/2}\Brkt{\frac{1}{(A-z_{0})(A-z_{0})\tp(A-z_{0})\adj}}=2\ii \im z_{0}(\absv{\gamma_{0}}+O(\absv{\im z_{0}}))	\\
		=&\,I_{4}(1)^{-1/2}\theta_{1}+O(N^{-1/2-\frc_{2}})=1-\frac{1}{z_{1}^{2}}+O(N^{-1/2-\frc_{2}})=2\ii\im z_{1}+o(N^{-1/2})
		.
	\end{split}\eeq
	This completes the proof of Proposition \ref{prop:real_path}. 
\end{proof}

Now, as in the complex case, Proposition \ref{prop:comp_G} in the real case is an immediate consequence of the following counterpart of Proposition \ref{prop:dt}. Recall the definition of $c_{t}$ from \eqref{eq:def_gamma}.
\begin{prop}\label{prop:dt_R}
	Let $\bbF=\R$, fix $k\in\N$ and $C>0$, and suppose that $\caA_{t}$ satisfies Assumption \ref{assu:real_flow}. Let $0<\delta<(\frc_{1}\wedge\frc_{2})/10$, where we recall $\eta_{0}= N^{-3/4+\delta}$. Then there exists a constant $c>0$ for which the following holds uniformly over $\caA_{t}$, $t\in[0,1]$, and $\absv{w},\norm{\bsw}\leq C$.
	\begin{align}\label{eq:comp_logd_d_R}
		\sum_{j=1}^{k}\E\Brkt{\frac{\dd(\E H_{t}^{w_{j}}-\ii\eta_{t})}{\dd t}G_{t}^{w_{j}}}\bsL^{(j)}_{t}(\bsw^{(j)})&=h_{t}(\bsw)+O(N^{-1-c}),	\\
		\frac{\dd}{\dd t}\E\frac{\brkt{G_{t}^{w}}}{c_{t}}&=O(N^{-c}\Psi),\label{eq:comp_sing_d_R}
	\end{align}
	where $h_{t}(\bsw)$ satisfies $(\prod_{j=1}^{k}\Delta_{w_{j}})h_{t}(\bsw)\equiv0$. 
\end{prop}

The rest of this section is organized as follows: In Section \ref{sec:real_intro}, we present a general framework of the proof in the real case, pointing out the main difference from the complex case. In particular, we present the decoupling lemma, Lemma \ref{lem:decoup}, adapted to the real case. Then in Sections \ref{sec:im L} and \ref{sec:im S} we prove Proposition \ref{prop:dt_R} under the assumptions (A1) and (A2), respectively. We emphasize that the rest of Section \ref{sec:real} applies even when the path $\caA_{t}$ is complex, as long as it satisfies Assumption \ref{assu:real_flow}. In particular, all estimates in the rest of this section depend only on $k,C,C_{5},\frc_{1},\frc_{2}$.
\subsection{General framework}\label{sec:real_intro}

In what follows, we introduce notations parallel to those in Section \ref{sec:comp_prep} in the complex case. Virtually the only difference is that, on top of symbols $\bsw$ and $\bsF$ in \eqref{eq:symbol}, we also have another tuple $\bsiota$ whose components vary in $\{1,\intercal\}$ (recall \eqref{eq:S_R}). As we will shortly see, the parallel notations allow us to verbatim reuse algebraic expressions appeared in the complex case. To avoid repetition, when we encounter an expression that already appeared in Sections \ref{sec:comp_prep}--\ref{sec:comp_comput}, we simply refer to an equation therein and do not rewrite the same expression again; see e.g. \eqref{eq:cancel_R_C} and \eqref{eq:cancel_R_R}.

The following definitions integrate the newly appeared transposes into the same framework as in Section \ref{sec:comp_Gini}. 
\begin{defn}
	For a given $n\in\N$ and $n$-tuples $\bsw,\bsF,\bsiota$, define
	\beq
	P_{n}(\bsw;\bsF;\bsiota)\deq (G^{w_{1}})^{\iota_{1}}F_{1}\cdots (G^{w_{n}})^{\iota_{n}}F_{n}.
	\eeq
	For a multi-index $\bsm=(m_{j})_{j}$ in $\N$ and $\absv{\bsm}$-tuples $\bsw=(\bsw_{j})_{j}$, $\bsF=(\bsF_{j})_{j}$, $\bsiota=(\bsiota_{j})_{j}$, define
	\beq
	\bsQ_{\bsm}(\bsw;\bsF;\bsiota)=\prod_{j}\brkt{P_{m_{j}}(\bsw_{j};\bsF_{j};\bsiota_{j})}.
	\eeq
	We say an $n$-tuple $\bsiota$ is \emph{mixed} if $\bsiota\neq (1,\ldots,1),(\intercal,\ldots,\intercal)$, in which case we say $P_{n}(\bsw;\bsF;\bsiota)$ is \emph{mixed} and write $P_{n}^{\mathrm{mix}}\equiv P_{n}$. The $n$-tuple $\bsiota=(1,\ldots,1)$ and $P_{n}(\bsw;\bsF;1,\ldots,1)$ are called \emph{pure}, and we write $P_{n}^{\mathrm{pure}}\equiv P_{n}$. Similarly, $\bsQ_{\bsm}$ is called \emph{mixed} if at least one of the involved factors of $P_{m_{j}}$ are \emph{mixed}, and we write $\bsQ_{\bsm}^{\mathrm{mix}}\equiv\bsQ_{\bsm}$ to emphasize this fact. The expression $\bsQ_{\bsm}$ is called \emph{pure} if all involved factors of $P_{m_{j}}$ are pure, in which case we write $\bsQ_{\bsm}^{\mathrm{pure}}\equiv\bsQ_{\bsm}$. Notice that $P_{n}(\bsw;\bsF;\intercal,\ldots,\intercal)$ is not called pure nor mixed.
\end{defn}
\begin{defn}  We naturally extend the definitions of $\caX_{n}$ and $\caY_{n}(B)$ from Definition \ref{def:caX_caY} by introducing the $\bsiota$-variables.  We define $\caX_{n}^{\mathrm{mix}}$ (resp. $\caX_{n}^{\mathrm{pure}}$) to be the subset of $\caX_{n}$ consisting of the same expressions as in \eqref{eq:def_caX} such that $\bsQ_{\bsm}$ is \emph{mixed} (resp. \emph{pure}). 
	
	Also, we write $\caY_{n}^{\mathrm{mix}}(B)$ (resp. $\caY_{n}^{\mathrm{pure}}(B)$) for the subset of $\caY_{n}(B)$ such that at least one (resp. both) of the factors $P_{n_{1}}$ and $\bsQ_{\bsm}$ in \eqref{eq:def_caY} is \emph{mixed} (resp. are \emph{pure}).
	As in Definition \ref{def:caX_caY}, we have $\caY_{n}^{\mathrm{mix}}(I)=\caX_{n}^{\mathrm{mix}}$ and $\caY_{n}^{\mathrm{pure}}(I)=\caX_{n}^{\mathrm{pure}}$.
\end{defn}
\begin{defn}
	With a light abuse of notations, we denote generic finite (independent of $N$) sums of elements of $\caX_{n},\caX_{n}^{\mathrm{mix}}$, and $\caX_{n}^{\mathrm{pure}}$ as $\frX_{n},\frX_{n}^{\mathrm{mix}}$, and $\frX_{n}^{\mathrm{pure}}$, respectively. Similarly we define $\frY_{n},\frY_{n}^{\mathrm{mix}},$ and $\frY_{n}^{\mathrm{pure}}$ to be generic finite sums over $\caY_{n}$, $\caY_{n}^{\mathrm{mix}}$, and $\caY_{n}^{\mathrm{pure}}$.
\end{defn}

We now state the decoupling lemma incorporating $G\tp$. Its proof is completely analogous to that of Lemma \ref{lem:decoup} hence omitted.
\begin{lem}\label{lem:decoup_R}
	Let $B_{d}$ and $B_{o}\in\C^{N\times N}$ be supported respectively on diagonal and off-diagonal blocks. Then we have
	\beq\begin{split}\label{eq:decoup_diag_R}
		&\E\brkt{B_{d}P_{n}^{\mathrm{mix}}}\bsQ_{\bsm}\wt{\bsL}	\\
		=&\,\,\brkt{2B_{d}(\mr{M}^{\wt{w}_{1}})^{\iota_{1}}(\mr{M}^{\wt{w}_{n}})^{\iota_{n}}F_{n}}\E\brkt{P_{n}^{\mathrm{mix}}}\bsQ_{\bsm}\wt{\bsL}	\\
		&+\left(\lone(n=1)\ii\eta_{t}\brkt{B_{d}((\mr{M}^{\wt{w}_{1}})^{\iota_{1}})^{2}F_{1}}+\lone(n=2)\frac{1}{N}\brkt{B_{d}(\mr{M}^{\wt{w}_{1}})^{\iota_{1}}F_{1}(\mr{M}^{\wt{w}_{2}})^{\iota_{2}}F_{2}}\right)\E\bsQ_{\bsm}\wt{\bsL}	\\
		&+\frac{1}{2}\E\frY^{\mathrm{mix}}_{\absv{\bsm}+n+1}(B_{d}(\mr{M}^{\wt{w}_{1}})^{\iota_{1}})\\
		&+\ii\eta\E\frY_{\absv{\bsm}+n}(B_{d}(\mr{M}^{\wt{w}_{1}})^{\iota_{1}})
		+\frac{1}{N}\E\frY_{\absv{\bsm}+n-1}(B_{d}(\mr{M}^{\wt{w}_{1}})^{\iota_{1}}),
	\end{split}
	\eeq	
	and
	\beq
	\begin{split}\label{eq:decoup_off_R}
		\E\brkt{B_{o}P_{n}^{\mathrm{mix}}}\bsQ_{\bsm}\wt{\bsL}	
		=&\,\,\frac{1}{2}\E\frY_{\absv{\bsm}+n+1}^{\mathrm{mix}}(B_{o}(\mr{M}^{\wt{w}_{1}})^{\iota_{1}})	\\
		&+\ii\eta\E\frY_{\absv{\bsm}+n}(B_{o}(\mr{M}^{\wt{w}_{1}})^{\iota_{1}})
		+\frac{1}{N}\E\frY_{\absv{\bsm}+n-1}(B_{o}(\mr{M}^{\wt{w}_{1}})^{\iota_{1}}),
	\end{split}\eeq
	where we abbreviated $P_{n}^{\mathrm{mix}}\equiv P_{n}^{\mathrm{mix}}(\wt{\bsw};\bsF;\bsiota)$. The same holds true if we remove all superscripts $\mathrm{mix}$, or if we replace $P_{n}^{\mathrm{mix}}$ on the left-hand sides by $P_{n}^{\mathrm{pure}}$ and $\frY_{\absv{\bsm+n+1}}^{\mathrm{mix}}$ on the right-hand side by $\frY_{\absv{\bsm}+n+1}^{\mathrm{pure}}+\frY_{\absv{\bsm}+n+1}^{\mathrm{mix}}$.
\end{lem}

The counterpart of \eqref{eq:KG_1} is as follows:
\beq\begin{split}\label{eq:KG_1_R}
	&\E\Brkt{\frac{\dd \E H_{t}^{w_{j}}}{\dd t}G^{w_{j}}E_{1}}\bsL^{(j)}	\\
	=&\Brkt{B_{d}^{w_{j}}E_{1}}\E\bsL^{(j)}
	+\ii\eta_{t}\E\Brkt{B_{d}^{w_{j}}G^{w_{j}}E_{1}}\bsL^{(j)}	\\
	&+\E\brkt{G^{w_{j}}}\Brkt{B_{d}^{w_{j}}G^{w_{j}}E_{1}}\bsL^{(j)}	
	+\frac{1}{N}\E\brkt{B_{d}^{w_{j}}E_{1}(G^{w_{j}})\tp E_{2}G^{w_{j}}E_{1}}\bsL^{(j)}	\\
	&-\frac{1}{N}\sum_{\ell}^{(j)}\E\Brkt{B_{d}^{w_{j}}(G^{w_{\ell}}+(G^{w_{\ell}})\tp)E_{2}G^{w_{j}}E_{1}}\bsL^{(j,\ell)}.
\end{split}\eeq
In what follows, we take care of the new terms in \eqref{eq:KG_1_R} involving $G\tp$ as well as similar new terms that appears when we expand
\beq
\E\brkt{B_{d}^{w_{j}}G^{w_{j}}E_{1}}\bsL^{(j)},\qquad \E\brkt{B_{d}^{w_{j}}G^{w_{\ell}}E_{2}G^{w_{j}}E_{1}}\bsL^{(j,\ell)}.
\eeq
In particular, we will shortly see that their contributions are negligible when $\absv{\theta_{t}}\gtrsim N^{-1/2+\frc_{1}}$, and cancel with each other when $\absv{\theta_{t}}\lesssim N^{-1/2+\frc_{1}}$. In the latter regime, we use an additional cancellation originating from Assumption \ref{assu:real_flow} (A2) that $\absv{\dd\theta_{t}/\dd t}$ is much smaller than $N^{-1/2}$.

\subsection{Large $\theta$ regime}\label{sec:im L}
Throughout this section, we assume that the path $\caA_{t}$ satisfies Assumption~\ref{assu:real_flow}~(A1).
\begin{lem}\label{lem:R_C}
	Fix a multi-index $\bsm$ in $\N$. Then we have
	\beq\label{eq:R_C}
	\absv{\E\bsQ_{\bsm}^{\mathrm{mix}}\wt{\bsL}}\lesssim N^{\epsilon}\absv{\theta}^{-1}\Psi^{\absv{\bsm}+2}.
	\eeq
	uniformly over $(w_{1},\ldots,w_{k})$ in a compact set.
\end{lem}
\begin{proof}
	Without loss of generality, suppose that the first factor $P_{m_{1}}$ in $\bsQ_{\bsm}$ is mixed with $\iota_{11}=1$ and $\iota_{1m_{1}}=\intercal$. Taking $B_{d}=I$ and applying \eqref{eq:decoup_diag_R}, we have
	\beq\label{eq:R_C_0}
	\E\bsQ_{\bsm}\wt{\bsL}=\brkt{2\mr{M}^{w_{11}}(\mr{M}^{w_{1n}})\tp F_{n}}\E\bsQ_{\bsm}\wt{\bsL}
	+\E\frY_{\absv{\bsm}+n+1}(\mr{M}^{w_{11}})+O(N^{\epsilon}\Psi^{\absv{\bsm}+2}).
	\eeq
	
	Notice that
	\beq\label{eq:MMT=theta}
	1-\brkt{2\mr{M}^{0}(\mr{M}^{0})\tp F_{n}}
	=\begin{cases}
		\ol{\theta} & \text{if }F_{n}=E_{1} \\
		\theta & \text{if }F_{n}=E_{2},
	\end{cases}
	\eeq
	so that by (A1)
	\beq\label{eq:R_C_1}
	\Absv{1-\brkt{2\mr{M}^{w_{1}}(\mr{M}^{w_{1n}})\tp F_{n}}}=\absv{\theta}+O(N^{-1/2})=\absv{\theta}(1+O(N^{-\frc_{1}})).
	\eeq
	On the other hand, since $\mr{M}^{w_{11}}$ is supported on off-diagonal blocks, \eqref{eq:decoup_off_R} implies
	\beq\label{eq:R_C_2}
	\E\frY_{\absv{\bsm}+n+1}(\mr{M}^{w_{11}})=O(N^{\epsilon}\Psi^{\absv{\bsm}+2}).
	\eeq
	Plugging \eqref{eq:R_C_1} and \eqref{eq:R_C_2} into \eqref{eq:R_C_0} proves the result.
\end{proof}
As an immediate consequence of Lemmas \ref{lem:decoup_R} and \ref{lem:R_C}, we have that
\beq\label{eq:R_C_c}
\absv{\E\frY_{n}^{\mathrm{mix}}(B_{d})}\lesssim N^{\epsilon}\norm{B_{d}}\absv{\theta}^{-1}\Psi^{n+2},\qquad 
\absv{\E\frY_{n}^{\mathrm{mix}}(B_{o})}\lesssim N^{\epsilon}\norm{B_{o}}\absv{\theta}^{-1}\Psi^{n+3},
\eeq
for any fixed $B_{d}$ and $B_{o}\in\C^{2N\times 2N}$ supported respectively on diagonal and off-diagonal blocks.

Armed with \eqref{eq:R_C} and \eqref{eq:R_C_c}, we may now prove the counterpart of Lemma \ref{lem:decoup_3_C_c}:
\begin{lem}\label{lem:R_C_1}
	Let $n,\bsm,$ and $B_{d}$ be as in Lemma \ref{lem:decoup_3_C}, and assume that $P_{n}\equiv P_{n}^{\mathrm{pure}}$ and $\bsQ_{\bsm}\equiv\bsQ_{\bsm}^{\mathrm{pure}}$. Then \eqref{eq:decoup3_C_c} remains true, modulo an additional error of $O(N^{\epsilon}\norm{B_{d}}\absv{\theta}^{-1}\Psi^{6})$.
\end{lem}
\begin{proof}
	Similarly to Lemma \ref{lem:decoup_3_C_c}, we only need to prove the following counterpart of \eqref{eq:decoup3_C}: The only difference is that the error has an additional factor of $\absv{\theta}^{-1}$.
	\beq\begin{split}\label{eq:decoup3_R}
		\E\brkt{B_{d}P_{n}}\bsQ_{\bsm}\wt{\bsL}
		=&\brkt{2B_{d}\mr{M}^{\wt{w}_{1}}\mr{M}^{\wt{w}_{n}}F_{n}}\E\brkt{P_{n}}\bsQ_{\bsm}\wt{\bsL}	\\
		&+\left(\lone(n=1)\ii\eta_{t}\brkt{B_{d}(\mr{M}^{\wt{w}_{1}})^{2}F_{1}}+\lone(n=2)\frac{1}{N}\brkt{B_{d}\mr{M}^{\wt{w}_{1}}F_{1}\mr{M}^{\wt{w}_{2}}F_{2}}\right)\E\bsQ_{\bsm}\wt{\bsL}	\\
		&+\frac{1}{2}\brkt{B_{d}(\mr{M}^{0})^{4}F_{n}}\E\frX_{4}^{\mathrm{pure}}
		+O(N^{\epsilon}\absv{\theta}^{-1}\norm{B_{d}}\Psi^{6}).
	\end{split}\eeq
	The proof of \eqref{eq:decoup3_R} follows that of \eqref{eq:decoup3_C} almost verbatim; we apply Lemma \ref{lem:decoup_R} three times instead of Lemma \ref{lem:decoup}. The only difference is that, in the first application of \eqref{eq:decoup_diag_R}, we decompose the second order term as
	\beq
	\E\frY_{\absv{\bsm}+n+1}(B_{d}\mr{M}^{\wt{w}_{1}}) =\E\frY_{\absv{\bsm}+n+1}^{\mathrm{pure}}(B_{d}\mr{M}^{\wt{w}_{1}})+\E\frY_{\absv{\bsm}+n+1}^{\mathrm{mix}}(B_{d}\mr{M}^{\wt{w}_{1}}),
	\eeq
	and estimate the mixed part using \eqref{eq:R_C_c}.
\end{proof}
Lemma \ref{lem:R_C_1} shows that the expansions of the `pure' terms in \eqref{eq:KG_1_R} is identical to those of \eqref{eq:KG_1}. Thus, it only remains to show that the `mixed' terms therein have no contribution. Here we record the application of Lemma \ref{lem:decoup_R} and \eqref{eq:R_C_c} to them:
\begin{lem}\label{lem:GGT}
	Let $B_{d}$ be supported on diagonal blocks. Then we have
	\beq\label{eq:GGT}\begin{split}
		\E\brkt{B_{d}(G^{\wt{w}_{1}})\tp E_{2}G^{\wt{w}_{2}}E_{1}}\wt{\bsL}
		=\brkt{2B_{d}(\mr{M}^{\wt{w}_{1}})\tp \mr{M}^{\wt{w}_{2}}E_{1}}\E\left(\brkt{(G^{\wt{w}_{1}})\tp E_{2}G^{\wt{w}_{2}}E_{1}}+\frac{1}{2N}\right)\wt{\bsL}	\\
		+O(N^{\epsilon}\norm{B_{d}}\absv{\theta}^{-1}\Psi^{6}).
	\end{split}\eeq
\end{lem}
\begin{proof}[Proof of Proposition \ref{prop:dt_R}, (A1)]
	We prove \eqref{eq:comp_logd_d_R} assuming that $\caA_{t}$ satisfies Assumption \ref{assu:real_flow} (A1).
	The proof of \eqref{eq:comp_sing_d_R} is omitted since it follows from taking the $\eta$-derivative of \eqref{eq:comp_logd_d_R}.
	
	Plugging the result of Lemmas \ref{lem:R_C_1} and \ref{lem:GGT} into \eqref{eq:KG_1_R}, we find that
	\beq\begin{split}\label{eq:cancel_R_C}
		\E\Brkt{\frac{\dd (\E H_{t}^{w_{j}}-\ii\eta_{t})}{\dd t}G^{w_{j}}E_{1}}&\bsL^{(j)}	
		=\text{(right-hand side of \eqref{eq:cancel_0}}	\\
		&+\brkt{2B_{d}^{w_{j}}(\mr{M}^{w_{j}})\tp \mr{M}^{w_{j}}E_{1}}\E\left(\brkt{(G^{w_{j}})\tp E_{2}G^{w_{j}}E_{1}}+\frac{1}{2N}\right)\bsL^{(j)}\\
		&+\sum_{\ell}^{(j)}\brkt{2B_{d}^{w_{j}}(\mr{M}^{w_{\ell}})\tp \mr{M}^{w_{j}}E_{1}}\E\left(\brkt{(G^{w_{\ell}})\tp E_{2}G^{w_{j}}E_{1}}+\frac{1}{2N}\right)\bsL^{(j,\ell)}\\
		&+O(N^{\epsilon}\absv{\theta}^{-1}\Psi^{6}).
	\end{split}\eeq
	To estimate the newly appeared terms, note that applying Lemma \ref{lem:GGT} to $B_{d}=I$ and then using \eqref{eq:MMT=theta} gives
	\beq
	\E\brkt{(G^{w_{\ell}})\tp E_{2}G^{w_{j}}E_{1}}\bsL^{(j,\ell)}
	=\frac{1+\theta}{\theta}\frac{1}{2N}\E\bsL^{(j,\ell)}+O(N^{\epsilon}\absv{\theta}^{-2}\Psi^{6}).
	\eeq
	On the other hand, the deterministic prefactor can be written as
	\beq
	\Brkt{2B_{d}^{w_{j}}(\mr{M}^{w_{\ell}})\tp \mr{M}^{w_{j}}E_{1}}
	=\Brkt{2\frac{\dd \mr{M}^{w_{j}}}{\dd t}(\mr{M}^{w_{\ell}})\tp E_{2}}
	=\frac{\dd\theta}{\dd t}+O(N^{-1/2}).
	\eeq
	Thus the second and third lines of \eqref{eq:cancel_R_C} can be written as
	\beq\label{eq:cancel_R_C1}
	\frac{\dd \theta}{\dd t}\frac{1+2\theta}{\theta}\frac{1}{2N}\left(\bsL^{(j)}+\sum_{\ell}^{(j)}\bsL^{(j,\ell)}\right)
	+O\left(N^{\epsilon}\absv{\theta}^{-1}\left(1+\Absv{\frac{1}{\theta}\frac{\dd\theta}{\dd t}}\right)\Psi^{6}\right).
	\eeq
	Notice that the first term of \eqref{eq:cancel_R_C1} is independent of $w_{j}$ hence harmonic. Finally, by Assumption \ref{assu:real_flow} (A1) we conclude that
	\beq
	\absv{\theta}^{-1}\left(1+\Absv{\frac{1}{\theta}\frac{\dd\theta}{\dd t}}\right)\Psi^{6}\lesssim (N^{1/2-\frc_{1}}+N^{1/2-\frc_{2}})\Psi^{6}\lesssim N^{-1-c}
	\eeq
	by taking $c<\min(\frc_{1},\frc_{2})-6\delta$. This concludes the proof of Proposition \ref{prop:dt_R} when $\caA_{t}$ satisfies (A1).
\end{proof}

\subsection{Small $\theta$ regime}\label{sec:im S}
In this section we assume that $\caA_{t}$ satisfies Assumption \ref{assu:real_flow} (A2).
\begin{lem}\label{lem:R_R}
	Let $n,\bsm,$ and $B_{d}$ be as in Lemma \ref{lem:decoup_3_C}. Then \eqref{eq:decoup3_C_c} remains true if we replace $\mr{M}^{\wt{w}_{i}}$ with $(\mr{M}^{\wt{w}_{i}})^{\iota_{i}}$ and add an additional error of $O(N^{\epsilon}\norm{B_{d}}\absv{\theta}^{1/2}\Psi^{4})$.
\end{lem}
\begin{proof}
	The proof is again three applications of Lemma \ref{lem:decoup_R}. The main difference from those of Lemmas~\ref{lem:decoup_3_C_c} and \ref{lem:R_C} is that, in the analogue of \eqref{eq:decoup3_R}, we keep all $\E\frY_{\absv{\bsm}+n+1}^{\mathrm{mix}}(B_{d}\mr{M}^{\wt{w}_{1}})$ terms by expressing it as a finite sum of
	\beq
	\brkt{B_{d}(\mr{M}^{\wt{w}_{1}})^{\iota_{1}}(\mr{M}^{w'_{1}})^{\iota'_{1}}(\mr{M}^{w'_{2}})^{\iota'_{2}}(\mr{M}^{w'_{3}})^{\iota'_{4}}F_{n}}\E\frX_{4},
	\eeq
	where $w'$ and $\iota'$ vary respectively in $\{w_{1},\ldots,w_{k}\}$ and $\{1,\intercal\}$.
	Then the deterministic prefactor is
	\beq
	\brkt{B_{d}(\mr{M}^{\wt{w}_{1}})^{\iota_{1}}(\mr{M}^{w'_{1}})^{\iota'_{1}}(\mr{M}^{w'_{2}})^{\iota'_{2}}(\mr{M}^{w'_{3}})^{\iota'_{4}}F_{n}}=\brkt{B_{d}(\mr{M}^{0})^{4}F_{n}}+O(\norm{B_{d}}(\sqrt{\re\theta}+N^{-1/2})),
	\eeq
	where we used Cauchy-Schwarz and 
	\beq
	\brkt{\absv{\mr{M}^{0}-(\mr{M}^{0})\tp}^{2}}=2\re\theta
	\eeq
	due to \eqref{eq:theta_asymp}. This proves the analogue of \eqref{eq:decoup3_C} with an additional error of $\absv{\theta}^{1/2}\Psi^{4}$, from which the result follows by solving for $\E\frX_{4}$.
\end{proof}

\begin{proof}[Proof of Proposition \ref{prop:dt_R}, (A2)]
	As before, we focus on the proof of \eqref{eq:comp_logd_d_R}. Applying Lemma \ref{lem:R_R} to each term of \eqref{eq:KG_1_R}, we have
	\beq\begin{split}\label{eq:cancel_R_R}
		&\E\Brkt{\frac{\dd(\E H_{t}^{w_{j}}-\ii\eta_{t})}{\dd t}G^{w_{j}}E_{1}}\bsL^{(j)}	\\
		=&\text{(right-hand side of \eqref{eq:cancel_0})}	\\
		&+\sum_{\ell}\left(\brkt{2B_{d}^{w_{j}}(\mr{M}^{w_{\ell}})\tp\mr{M}^{w_{j}}E_{1}}-\frac{\brkt{2B_{d}^{w_{j}}(\mr{M}^{0})^{4}E_{1}}}{I_{4}}(\brkt{2(\mr{M}^{w_{\ell}})\tp\mr{M}^{w_{j}}E_{1}}-1)\right)	\\
		&\qquad\qquad\times\E\frac{\brkt{(G^{w_{\ell}})\tp E_{2}G^{w_{j}}E_{1}}}{N}\bsL^{(j,\ell)}	\\
		&+\sum_{\ell}^{(j)}\frac{1}{2N}\left(\brkt{2B_{d}^{w_{j}}(\mr{M}^{w_{j}})\tp \mr{M}^{w_{j}}E_{1}}-\frac{\brkt{2B_{d}^{w_{j}}(\mr{M}^{0})^{4}E_{1}}}{I_{4}}\right)\E\bsL^{(j,\ell)}	\\
		&+O(N^{\epsilon}\Psi^{5})+O(N^{-1/4-\frc_{1}/2+\epsilon}\Psi^{4}),
	\end{split}\eeq
	where we used the first assertion of (A2) to estimate $\sqrt{\re \theta}\leq N^{-1/4+\frc_{1}/2}$. By choosing small enough $\delta,c>0$, we find that the error in \eqref{eq:cancel_R_R} is $O(N^{-1-c})$. Also, the third term of \eqref{eq:cancel_R_R} is independent of $w_{\ell}$, hence its real part is harmonic. 
	
	For the second term of \eqref{eq:cancel_R_R}, we first focus on the case when $\ell=j$. Note that
	\beq\begin{split}
		\brkt{2B_{d}^{w}(\mr{M}^{w})\tp\mr{M}^{w}E_{1}} 
		&=-\frac{1}{2}\frac{\dd}{\dd t}\Brkt{\frac{1}{(\caA-N^{-1/2}\gamma^{-1}w)(\caA-N^{-1/2}\gamma^{-1}w)\tp}}	\\
		&=\frac{1}{2}\frac{\dd}{\dd t}\theta -N^{-1/2}w\frac{\dd}{\dd t}\left[\gamma^{-1}\Brkt{\frac{1}{\caA^{2}\caA\tp}}\right]+O(N^{-1}\norm{\dd\caA_{t}/\dd t})	\\
		&=\frac{1}{2}\frac{\dd}{\dd t}\theta +N^{-1/2}w\frac{\dd [I_{4}^{1/2}]}{\dd t}+O(N^{-1/2-\frc_{2}}),
	\end{split}\eeq
	where in the third equality we used the last assertion of (A2) and that
	\beq
	\Brkt{\frac{1}{\caA^{2}\caA\tp}}-I_{3}=\Brkt{\frac{1}{\caA\adj}\Im[\caA]\frac{1}{\caA^{2}\caA\tp}}.
	\eeq
	Similarly, we have the two estimates
	\beq\begin{split}
		\brkt{2(\mr{M}^{w})\tp\mr{M}^{w}E_{1}}-1
		&=-\theta-2N^{-1/2}I_{4}^{1/2}w+O(N^{-1}),	\\
		\brkt{2B_{d}^{w}(\mr{M}^{0})^{4}E_{1}}
		&=-\frac{1}{4}\frac{\dd I_{4}}{\dd t}+O(N^{-1/2}).
	\end{split}\eeq
	Therefore we have
	\beq\begin{split}
		&\brkt{2B_{d}^{w_{j}}(\mr{M}^{w_{\ell}})\tp\mr{M}^{w_{j}}E_{1}}-\frac{\brkt{2B_{d}^{w_{j}}(\mr{M}^{0})^{4}E_{1}}}{I_{4}}(\brkt{2(\mr{M}^{w_{\ell}})\tp\mr{M}^{w_{j}}E_{1}}-1)	\\
		=&\frac{1}{2}\frac{\dd\theta}{\dd t}-\frac{1}{4 I_{4}}\frac{\dd I_{4}}{\dd t}\theta
		+N^{-1/2}\left(\frac{\dd [I_{4}^{1/2}]}{\dd t}-\frac{1}{2}I_{4}^{-1/2}\frac{\dd I_{4}}{\dd t}\right)w
		+O(N^{-1/2-\frc_{2}})	\\
		=&\frac{1}{2I_{4}^{1/2}}\frac{\dd (I_{4}^{-1/2}\theta)}{\dd t}+O(N^{-1/2-\frc_{2}})=O(N^{-1/2-\frc_{2}}).
	\end{split}\eeq
	This shows that the second term of \eqref{eq:cancel_R_R} for $\ell=j$ is bounded by
	\beq
	O(N^{-1/2-\frc_{2}})\Absv{\E\frac{\brkt{(G^{w_{j}})\tp E_{2}G^{w_{j}}E_{1}}}{N}\bsL^{(j)}}\lesssim N^{-1/2-\frc_{2}+\epsilon}\Psi^{2}=N^{-1-\frc_{2}+2\delta+\epsilon}\leq N^{-1-c}.
	\eeq
	The sum over $\ell\neq j$ can be dealt with a similar symmetrization as in \eqref{eq:cancel_2pt} -- \eqref{eq:cancel_2pt_c}, leading to $O(N^{-1-c})$. This concludes the proof of Proposition \ref{prop:dt_R} under Assumption \ref{assu:real_flow} (A2).
\end{proof}

\section{Proof of Proposition \ref{prop:comp_IID}}\label{sec:short GFT}
In this section we prove Proposition \ref{prop:comp_IID}, which compares the two deformed non-Hermitian ensembles $A+X$ and $A+X^{\mathrm{Gin}(\bbF)}$. For the sake of brevity we explain only how the proof follows exactly along the same lines as \cite{Cipolloni-Erdos-Schroder2021}. All estimates in this section are uniform over $z_{0}\in\caD_{C}$ defined in \eqref{eq:regular} and $X,A$ satisfying \eqref{eq:momentass}, \eqref{eq:Abound}.
\begin{proof}[Proof of Proposition \ref{prop:comp_IID}]
	Since the same argument applies to both $\bbF=\bbR$ and $\bbF=\bbC$, we omit the superscript $(\bbF)$ throughout the proof.
	
	We consider yet another interpolation $A+X_{t}$ between $A+X$ and $A+X^{\mathrm{Gin}}$, where $(X_{t})_{t>0}$ is the matrix-valued Ornstein--Uhlenbeck process (running from $t=0$ to $\infty$, unlike $(\caA_{t})_{t\in[0,1]}$ in Section \ref{sec:flow})
	\beq
	\dd X_{t}=-\frac{1}{2}X_{t}+\frac{1}{\sqrt{N}}\dd B_{t},\qquad X_{0}=X,
	\eeq
	and the entries of $B_{t}\in\bbF^{N\times N}$ are i.i.d. standard Brownian motions. As in Section \ref{sec:flow}, we write
	\beq\begin{aligned}\label{eq:L_OU}
		G_{t}^{\mathrm{OU},w}(\ii\eta)&\deq 	\begin{pmatrix}-\ii\eta & (A+X_{t}-z_{0}-N^{-1/2}w) \\
			(A+X_{t}-z_{0}-N^{-1/2}w)\adj & -\ii\eta\end{pmatrix}^{-1},\\
		L_{t}^{\mathrm{OU}}(w)&\deq \Tr\log (\absv{A+X_{t}-z_{0}-N^{-1/2}w}^{2}+\eta_{0}^{2})
		-\Tr\brkt{\log(\absv{A+\bsx-z_{0}-N^{-1/2}w}^{2}+\eta_{0}^{2})}_{\caM},	\\
		N_{t}^{\mathrm{OU}}(w)&\deq \frac{1}{c_{0}}\Brkt{\frac{\eta_{0}}{\absv{A+X_{t}-z_{0}-\gamma_{0}^{-1}N^{-1/2}w}^{2}+\eta_{0}^{2}}} =\frac{1}{c_{0}}\brkt{\im G_{t}^{\mathrm{OU},w}(\ii\eta_{0})},
	\end{aligned}\eeq
	where the last equality in the last line is the same identity as in \eqref{eq:N=G}.
	Since the first moment $A-z_{0}-N^{-1/2}w$ of our ensemble never moves, there is no need to define any dynamics in $z,A,\gamma$ or $c$ unlike in Section \ref{sec:flow}.
	
	Since the local law Theorem \ref{thm:local law: Local Law} immediately applies to $G_{t}^{\mathrm{OU},w}$ and the first two joint cumulants of $(A+X_{t})_{ij}$ remains constant over $t\geq 0$, the proof of Proposition \ref{prop:comp_IID} follows nearly exactly \cite{Cipolloni-Erdos-Schroder2021}. To be precise, following the same\footnote{The quantity $L_{t}^{\mathrm{OU}}(w)$ is a slight variant of $I_{3}(t)$ defined in \cite[Eq. (36)]{Cipolloni-Erdos-Schroder2021}; in $I_{3}(t)$ the variable $w$ is integrated out and that the logarithm is cut-off at $N^{100}$. These differences do not affect the proof.} cumulant expansions as in the proof of \cite[Lemma 2]{Cipolloni-Erdos-Schroder2021}, we find that for some constant $c>0$ and any $\xi>0$
	\beq\label{eq:short GFT log bound}
	\Absv{\frac{\dd}{\dd t}\E\prod_{j=1}^{k}L_{t}^{\mathrm{OU}}(w_{j})}
	\lesssim \e{-3t/2}\frac{N^{\xi}}{N^{5/2}\eta_{0}^{3}}
	\eeq
	holds uniformly over $t\geq 0$. Similarly, the same proof as \cite[Proposition 3]{Cipolloni-Erdos-Schroder2021} gives for any $\xi>0$ that
	\begin{equation}\label{eq:short GFT trace G bound}
		\Absv{\frac{\dd}{\dd t}\E N_{t}^{\mathrm{OU}}(w)}=\frac{1}{c_{0}}\Absv{\frac{\dd}{\dd t}\E\brkt{G^{\mathrm{OU},w}_{t}(\ii\eta_{0})}}\lesssim \e{-3t/2}\frac{N^{\xi}}{N^{7/2}\eta_{0}^4},
	\end{equation}
	uniformly over $t\geq 0$. 
	Integrating \eqref{eq:short GFT log bound} and \eqref{eq:short GFT trace G bound} over $t$ completes the proof.
\end{proof}

	\begin{rem}[Universality of the least singular value]
		The same argument as in the proof of Theorem \ref{thm:main1} can be used to prove the universality for the smallest singular value $\lambda_{1}$ of $A+X-z_{0}$. Again, we break the comparison up into:\begin{enumerate}
			\item\label{item:sing value step one}  Comparing the smallest singular values $A+X^{\mathrm{Gin}}$ and $X^{\mathrm{Gin}}$.
			
			\item\label{item:sing value step two} Comparing the smallest singular values of $A+X$ and $A+X^{\mathrm{Gin}}$.
		\end{enumerate} For \eqref{item:sing value step one} our proofs of \eqref{eq:comp_sing_d} and \eqref{eq:comp_sing_d_R} can be adapted, beyond the lower-tail estimate as in \eqref{eq:comp_sing}. More precisely, we can prove that for each fixed $w=O(N^{-1/2})$ and continuous and compactly supported test function $F$ on $\R$
		\beq
		\E F(\lambda_{1}(c_{0}(A+X-z_{0}-\gamma_{0}^{-1}w)))=\E F(\lambda_{1}(X-1-w))+o(1).
		\eeq
		A rough sketch of the proof is as follows:  Instead of $\E c_{0}^{-1}\brkt{G_{0}(\ii\eta_{0})}$, we take an arbitrary moment (or a more general function) of the regularized counting function
		\beq
		\int_{0}^{E_{1}}c_{0}^{-1}\im\Tr[G_{0}^{w}(c_{0}^{-1}(E+\ii\eta_{1})]\dd E,
		\eeq
		where $E_{1}=O(N^{-3/4+\epsilon})$ lives on the same scale as the typical eigen-spacing around $\lambda_{1}$. Then it suffices to prove that the time-derivative of this statistic is $o(1)$. To this end, we need to generalize the local law and the content of Section \ref{sec:comp_Gini} for spectral parameter $\wh{z}_{t}\deq c_{t}^{-1}(E+\ii\eta_{1})$ beyond $\eta_{t}\deq c_{t}^{-1}\eta_{1}$. Such a generalization does not require a significant changes in Section \ref{sec:comp_Gini}, since $\wh{z}_{t}$ is close enough to the imaginary axis that the cubic equation \eqref{eq:cubic} associated to the MDE remains almost intact. For \eqref{item:sing value step two} similarly to the proof of Proposition \ref{prop:comp_IID} we can compare along an Ornstein--Uhlenbeck process and use the slightly extended local law described above to control the time derivative of moments of the regularized counting function. This same approach can also be applied to the $p$-th singular value for fixed $p\in\N$.

		To adapt the Zig-Zag approach to Theorem \ref{thm:local law: Local Law} to spectral parameters $\wh{z}=E+\ii\eta$ off the imaginary axis the key is to correctly adjust the Zig step for the spectral parameter. The dynamics would \begin{equation*}
			\frac{\dd \wh{z}_s}{\dd s}=-\frac{\wh{z}_s}{2}-\brkt{M(\wh{z}_s)}
		\end{equation*} Again, because \eqref{eq:cubic} remains essentially intact for small $E$ the flow in the Zig step will actually only slightly change.
	\end{rem}

	\section{Zig--zag strategy: Proof of Theorem \ref{thm:local law: Local Law}}
	\label{sec:Proof of the local law}
	In this section we prove Theorem \ref{thm:local law: Local Law}. All estimates in this section are uniform over $X,A,z_{0}$ satisfying \eqref{eq:momentass}, \eqref{eq:Abound}, and $\norm{(A-z_{0})^{-1}}\leq C_{4}$. It is worth noting, as we only consider the local law for the argument of $G$ and $M$ on the imaginary axis, that from \eqref{eq:M boundedness} we have $\Norm{M}\le C$, for some $N$--independent $C>0$. For any $z$ such that $|z-z_0|<c$ and for every $T>0$ we define the following dynamics for $t\in[0,T]$.

	\begin{defn}[Matrix Ornstein--Uhlenbeck flow and the $\eta$-flow]\label{def:zig zag dynamics}
		Define the matrix-valued Ornstein--Uhlenbeck (OU) process $X_{t}$ as
		\beq
		\dd X_{t}=-\frac{1}{2}X_{t}\dd t+\frac{1}{\sqrt{N}}\dd B_{t},\qquad X_{0}=X,
		\eeq
		where the entries of $B_{t}\in\C^{N\times N}$ are i.i.d.\ standard Brownian motions. We define $\caA_{t}$ as \footnote{The evolutions $\caA_{t}$ and $\eta_{t}$ in this section are not related to those in Section \ref{sec:flow}.}
		\beq
		\caA_{t}\deq \e{-\frac{t}{2}}\caA_{0},\qquad \caA_{0}=A-z,
		\eeq and define the resolvent \begin{equation}\label{eq:G_t}
			G_t(\wh{z})\deq\left(\caH_{\caA_t+X_t}-\wh{z} \right)^{-1}, \qquad\quad \widehat{z}\in \C\setminus\R,
		\end{equation}
		where for a matrix $B$\begin{equation*}
			\caH_{B}\deq\begin{pmatrix}
				0&B\\
				B\adj&0
			\end{pmatrix}.
		\end{equation*} Note then that $\caA_{t}+X_{t}$ satisfies the SDE
		\beq
		\label{eq:OUflow}
		\dd (\caA_{t}+X_{t})=-\frac{1}{2}(\caA_{t}+X_{t})\dd t+\frac{1}{\sqrt{N}}\dd B_{t}.
		\eeq
		For $\eta>0$, we define $\eta_{t}\in\R$ to be the solution of the ODE
		\beq\begin{aligned}\label{eq:local law:eta dynamics}
			\frac{\dd \eta_{t}}{\dd t}=-\frac{\eta_{t}}{2}-\brkt{\im M_{\caA_{t}}(\ii\eta_{t})},\qquad \eta_{0}=\eta.
		\end{aligned}\eeq
	\end{defn}

	\begin{lem}
		$M_t\equiv M_{\caA_{t}}(\ii\eta_{t})$ satisfies
		\beq\label{eq:dM/dt}
		\frac{\dd \brkt{M_t}}{\dd t}=\frac{\brkt{M_t}}{2}.
		\eeq
	\end{lem} 

	\begin{proof}
		To see \eqref{eq:dM/dt} note that $M_t$ and the Hermitization $ \caH_{\caA_{t}}$ of $\caA_{t}$ satisfy
		\beq
		1=(\caH_{\caA_t}-\ii\eta_t)M_t-\brkt{M_t}M_t,
		\eeq
		and then we take the time derivative to find that
		\beq
			\left(\frac{\brkt{M_t}}{2}-\Brkt{\frac{\dd M_t}{\dd t}}\right)M_t-\frac{1}{M_t}\left(\frac{M_t}{2}-\frac{\dd M_t}{\dd t}\right)=0.
\eeq
		Multiplying both sides by $M_t$, taking the trace, and using $|\brkt{M_t^{2}}|\leq\brkt{|M_t|^{2}}<1$, we obtain \eqref{eq:dM/dt}. 
	\end{proof}
	
We define the time at which the characteristic curve $\eta_{t}$ touches the real axis
\begin{equation}
\label{eq:fintime}
T^*(\mathcal{A}_0,\eta_0):=\sup\{t: \mathrm{sgn}(\eta_t)=\mathrm{sgn}(\eta_0)\}.
\end{equation}
We will also use the short--hand notation $\rho_t:=\rho^{\mathcal{A}_t}(\ii\eta_t)$.
	
\begin{prop}[Zig step]
\label{pro:zignew}
Fix small constants $c, \epsilon>0$ and an arbitrary small $\xi>0$. Fix $\mathcal{A}_0\in \C^{N\times N}$ with $\lVert \mathcal{A}_0^{-1}\rVert\le C$ for some $C>0$, 
and pick $\eta_0$ such that $0<|\eta_0|/\rho_0\le c$. For $t\le T^*(\mathcal{A}_0,\eta_0)$, let $\mathcal{A}_t:=e^{-t/2}\mathcal{A}_0$, and let $\eta_t$ be the solution of \eqref{eq:local law:eta dynamics} with initial condition $\eta_0$. Define $G_t(\rmi\eta_t):=(\mathcal{H}_{\mathcal{A}_t+X_t}-\rmi\eta_t)^{-1}$. 
For some small $0<\gamma\le \epsilon/10$, assume that, with very high probability, we have
\begin{equation}
\label{eq:incond0}
\absv{ \brkt{B(G_0(\ii\eta_0)-M_{\mathcal{A}_0}(\ii\eta_0)) } }\le \frac{N^\gamma }{N\eta_0}, \qquad \left|\bsu^* \left(G_{0}(\ii\eta_0)-M_{\mathcal{A}_0}(\ii\eta_0)\right)\bsy \right|\le N^\gamma\sqrt{\frac{\rho_0}{N\eta_0}}.
\end{equation}
uniformly in unit vectors $\bsu,\bsy\in\C^{2N}$ and matrices $B\in\C^{2N\times 2N}$ with $ \lVert B\rVert\le 1$. Then\footnote{We point out that as a result of the Zig step we would get \eqref{eq:finalt} with $N^\xi$ replaced with the smaller $(\log N)^2$. We presented \eqref{eq:finalt} directly with the worse $N^\xi$ since this is the error that we will get by the application of the Zag step in Proposition~\ref{prop:new zag} below.},
\begin{equation}
\label{eq:finalt}
\absv{ \brkt{B(G_t(\ii\eta_t)-M_{\mathcal{A}_t}(\ii\eta_t)) } }\le  \frac{N^{\gamma+\xi}}{N\eta_t}, \qquad \left|\bsu^*\left(G_{0}(\ii\eta_t)-M_{\mathcal{A}_t}(\ii\eta_t)\right)\bsy \right|\le N^{\gamma+\xi}\sqrt{\frac{\rho_t}{N\eta_t}},
\end{equation}
with very high probability, uniformly in $t\le T^*(\mathcal{A}_0,\eta_0)$, $N\eta_t\rho_t\ge N^\epsilon$, unit vectors $\bsu,\bsy\in\C^{2N}$ and matrices $B\in\C^{2N\times 2N}$ with $ \lVert B\rVert\le 1$.
\end{prop}

We remark that, even if stated together, the average local law bounding $\absv{ \brkt{B(G_t(\ii\eta_t)-M_{\mathcal{A}_t}(\ii\eta_t)) } }$ in Proposition \ref{pro:zignew} can be proven independently of the isotropic local law bounding $\left| \bsu^*\left(G_{0}(\ii\eta_t)-M_{\mathcal{A}_t}(\ii\eta_t)\right)\bsy \right|$.

Fix $\mathcal{A}=A-z\in\C^{N\times N}$, and recall that $\lVert \mathcal{A}^{-1}\rVert\le C$. In the following for given $\eta>0$ and $r\ge 1$ we define the \emph{scale}
\begin{equation}
\mathfrak{s}(\eta,r):=\frac{\eta}{\rho^{r\mathcal{A}}(\ii\eta)}=\Brkt{\Brkt{\frac{1}{\absv{r\caA+\bsx}^{2}+\eta^{2}}}_{\caM}}^{-1}.
\end{equation}
Note that the map $\eta\mapsto \frs(\eta,r)$ is monotone increasing. Thus, given $r\ge 1$, there is a one to one correspondence between $\eta$ and $\mathfrak{s}(\eta,r)$. In particular, when we say that a local law holds on some given scale
$\mathfrak{s}^*$, we mean that it holds for spectral parameter $\ii\eta$ determined by  $\mathfrak{s}(\eta,r)=\mathfrak{s}^*$.
 Fix any $\epsilon>0$, and define the scale $\mathfrak{s}_\mathrm{fin}=\mathfrak{s}(\eta_\mathrm{fin},1)$ so that
\begin{equation}
\label{eq:finscale}
N\eta_\mathrm{fin}\rho^{\mathcal{A}}(\ii\eta_\mathrm{fin}):=N^\epsilon.
\end{equation} 
Since $\rho^{r\mathcal{A}}(\ii\eta)$ is uniformly $O(1)$ over $r\in[0,1]$ and $\eta\in[0,1]$, we easily find that $\eta_{\mathrm{fin}}\gtrsim N^{-1+\epsilon}$ and that $\frs_{\mathrm{fin}}\gtrsim \eta_{\mathrm{fin}}\gtrsim N^{-1+\epsilon}$.

\begin{prop}[Zag step]\label{prop:new zag}
Fix an arbitrary small $\xi>0$, let $\epsilon>0$ from \eqref{eq:finscale}, fix any small $0<\gamma\le\epsilon/10$,
and pick $\epsilon_0>0$. Fix a scale $\mathfrak{s}_\mathrm{fin}<\widetilde{\mathfrak{s}}_0=\mathfrak{s}(\widetilde{\eta_0},r) \le N^{-\epsilon_0}$,  and let $t\le \widetilde{\mathfrak{s}}_0$. Consider an $\widetilde{\eta_1}$ such that 
$\eta_{\mathrm{fin}}\vee \widetilde{\eta_0} N^{-\epsilon_0/100}\le \widetilde{\eta_1}\le \widetilde{\eta_0}$ and 
its  corresponding scale $\widetilde{\mathfrak{s}}_1=\mathfrak{s}(\widetilde{\eta_1},r)$ also satisfies
$\mathfrak{s}_\mathrm{fin}\vee N^{-\epsilon_0/200}\widetilde{\mathfrak{s}}_0\le \widetilde{\mathfrak{s}}_1\le \widetilde{\mathfrak{s}}_0$. 
  Define $G^t(\ii\eta):=(\mathcal{H}_{r\mathcal{A}+X_t}-\ii\eta)^{-1}$, and assume that, uniformly in $r\in [1,2]$, unit vectors $\bsu,\bsy\in\C^{2N}$, and matrices $B\in\C^{2N\times 2N}$ with $\lVert B\rVert\le 1$, the local law
\begin{equation}
\label{eq:llawzag} 
\absv{ \brkt{B(G^s(\ii \eta)-M_{r\mathcal{A}}(\ii \eta)) }}\le \frac{N^{\gamma} }{N\eta}, \qquad \left|\bsu^*\big(G^{s}(\ii\eta)-M_{r\mathcal{A}}(\ii\eta)\big)\bsy \right|\le N^{\gamma} \sqrt{\frac{\rho^{r\mathcal{A}}(\ii\eta)}{N\eta}}
\end{equation}
holds with very high probability for any time $s\in [0,t]$ on the scale $\widetilde{\mathfrak{s}}_0$, and at time $s=t$ on the scale $\widetilde{\mathfrak{s}}_1$. Then, \eqref{eq:llawzag} also holds on the scale $\widetilde{\mathfrak{s}}_1$ for any $s\in [0,t]$ with $\gamma$ replaced by $\gamma+\xi$. 
\end{prop}
Note that the deformation $\mathcal{A}$ and the spectral parameter $\ii\eta$  
are time-dependent in Proposition \ref{pro:zignew}, 
while they are time-independent in Proposition \ref{prop:new zag} but $\mathcal{A}$ carries an extra rescaling parameter $r$ 
and the estimate is uniform in $r$. This important subtlety is expressed by the lower and upper time indices $G_t$
and $G^t$ in the two propositions above. 
 In the application of Proposition \ref{prop:new zag} for proving Theorem \ref{thm:local law: Local Law}, $r$ will be a frozen time parameter from the previous
 step.

Our goal is now to prove \eqref{eq:llawzag} for $r=1$ and any $\eta$ such that $\mathfrak{s}(\eta,1)\ge \mathfrak{s}_\mathrm{fin}$; this would give Theorem~\ref{thm:local law: Local Law}. 
Fix $\epsilon_0>0$, and define $T:=N^{-\epsilon_0}$. Then, by standard ODE theory,
 there exists $\eta_\mathrm{in}=\eta_\mathrm{in}(\mathcal{A},T)$ such that $\mathfrak{s}_\mathrm{in}=\mathfrak{s}(\eta_\mathrm{in},e^{T/2})\sim T$ and $\mathfrak{s}(\eta_T,1)=\mathfrak{s}_\mathrm{fin}$, with $\eta_t$ being the solution of \eqref{eq:local law:eta dynamics} with initial condition $\eta_\mathrm{in}$. In particular, $T\le T^*(e^{T/2} \mathcal{A}, \eta_\mathrm{in})$; here we also used that the map $t\mapsto\eta_t$ is decreasing. Note that the $\eta_t$--dynamics depends on $\mathcal{A},T$ through $M_{e^{(T-t)/2}\mathcal{A}}$. Next, we define a decreasing sequence of  scales
\begin{equation}
\mathfrak{s}_k:=\mathfrak{s}_\mathrm{in} N^{-k\epsilon_0/100} \vee \mathfrak{s}_\mathrm{fin}.
\end{equation}
Note that $\mathfrak{s}_0=\mathfrak{s}_\mathrm{in}$, and let $K$ be the smallest index such that $\mathfrak{s}_K=\mathfrak{s}_\mathrm{fin}$, i.e. the last step may be a bit smaller than all the previous ones. We emphasize that $K\leq 100/\epsilon_{0}$ hence $K=O(1)$, which follows from $\frs_{\mathrm{fin}}\gtrsim N^{-1}$. We  define $0\le t_k\le T$ implicitly via
\begin{equation}\label{eq:def_tk}
\mathfrak{s}(\eta_{t_k},e^{(T-t_k)/2})= \mathfrak{s}_k;
\end{equation}
in particular, $t_K=T$ as a consequence of $\mathfrak{s}_K=\mathfrak{s}_\mathrm{in}=\mathfrak{s}(\eta_T,1)$. Since $\eta_{t}$ is decreasing and $\rho^{\e{-t/2}\caA}(\ii\eta_{t})$ is increasing in $t$ by \eqref{eq:local law:eta dynamics} and \eqref{eq:dM/dt}, we find that
\beq
	\mathfrak{s}(\eta_{t},\e{(T-t)/2})=\frac{\eta_{t}}{\rho^{\e{(T-t)/2}\caA}(\ii\eta_{t})}
\eeq
is decreasing in $t$. In particular, the implicit equation \eqref{eq:def_tk} has a unique solution and $t_k$ is an increasing sequence.
Further, we define
 domains
\begin{equation}
\mathcal{D}_k:=\big\{(\eta,r) : \eta\ge\eta_{t_k}, \,\, r\in [1, e^{(T-t_k)/2}] \,\big\}, \qquad\quad k=0,\dots, K.
\end{equation}
Note that
\begin{equation}\label{eq:eta_rat}
	\frac{\eta_{t_{k+1}}}{\eta_{t_{k}}}=\frac{\frs_{k+1}}{\frs_{k}}\frac{\rho^{\e{(T-t_{k})/2}\caA}(\ii\eta_{t_{k}})}{\rho^{\e{(T-t_{k+1})/2}\caA}(\ii\eta_{t_{k+1}})}=N^{-\epsilon_{0}/100}\e{-(t_{k+1}-t_{k})/2}
	\geq N^{-\epsilon_{0}/100}\e{-T}.
\end{equation}
We say that the local law holds at a tuple $(\eta,r)$ with error $\gamma$ if \eqref{eq:llawzag} holds. We are now ready to prove Theorem~\ref{thm:local law: Local Law}. 

\begin{proof}[Proof of Theorem~\ref{thm:local law: Local Law}]
		
We fix some $z$ such that $|z-z_{0}|<c$ and $\epsilon>0$. As above, we define the scale $\mathfrak{s}_\mathrm{fin}=\mathfrak{s}(\eta_\mathrm{fin},1)$ with $\caA=A-z$.  We additionally fix $0<\epsilon_0<\epsilon/100$, and define $T:=N^{-\epsilon_0}$, $\eta_{\mathrm{in}}$,  $\mathfrak{s}_{\mathrm{in}}$, $\mathfrak{s}_{k}$, $K$, and $\caD_{k}$ as above. We first establish that \eqref{eq:llawzag} holds on $\caD_{K}$. Extending \eqref{eq:llawzag} to hold uniformly over $|z-z_{0}|<c$ then follows from a simple union bound, and extending to $\eta\geq N^{-1}$ can be proved completely analogously to \cite[Appendix A]{Cipolloni-Erdos-Schroder2021}. 
	
From \cite[Theorem 2.1]{Erdos-Kruger-Schroder2019} we know that \eqref{eq:llawzag} holds with very high probability uniformly over $\caD_{0}$ with error $(\epsilon_0,0)$. We now want to conclude a local law on $\caD_{k}$ by induction. Fix the tolerance parameter $\xi$ so that $\epsilon_{0}+2K\xi<\epsilon/10$. The main inductive step is formalized in the following proposition, whose proof is presented after we finish the proof of Theorem~\ref{thm:local law: Local Law}.
\begin{prop}\label{prop:zigzag induction}
 Fix $k\in \{0,\dots, K-1\}$, and assume that \eqref{eq:llawzag} holds on $\mathcal{D}_k$ with error $\epsilon_0+2k\xi$. Then  \eqref{eq:llawzag} holds on $\mathcal{D}_{k+1}$ with error $\epsilon_0+2(k+1)\xi$.
\end{prop}

It then follows from $K$ applications of Proposition \ref{prop:zigzag induction} that \eqref{eq:llawzag} holds with very high probability uniformly over $\caD_{K}$ with error $\epsilon_0+2K\xi$. This completes the proof.
\end{proof}

\begin{proof}[Proof of Proposition~\ref{prop:zigzag induction}]
Fix some $(\eta,r)\in\mathcal{D}_{k+1}\setminus \mathcal{D}_{k}$ and set $\wh\eta=\eta_{t_{k}}$. 
We note, by the definition of our domains, for the solution $\wh{\eta}_{t}$ to \eqref{eq:local law:eta dynamics} with initial condition $\wh\eta_{0}=\wh\eta$ satisfies $(\eta,r)=(\wh\eta_{\wh t}, r)$ and $(\wh\eta,\e{\wh t/2} r)\in\caD_{k}$ for some $0<\wh t\leq t_{k+1}-t_{k}$. By assumption, \eqref{eq:llawzag} holds with very high probability on the scale $\mathfrak{s}(\wh\eta,\e{\wh t}r)$ with error $\epsilon_{0}+2k\xi$. Then, by Proposition~\ref{pro:zignew}, \eqref{eq:llawzag} holds on the scale $\mathfrak{s}(\eta,r)$ for $s=\wh t$ with error $\epsilon_0+(2k+1)\xi$. Additionally, \eqref{eq:llawzag} holds on the scale $\mathfrak{s}(\wh\eta,r)$ uniformly over $s\in[0,\wh t]$.
We then aim at applying Proposition \ref{prop:new zag} with the choices $(\wt{\eta}_{0},\wt{\eta}_{1})=(\wh{\eta},\eta)$ and $(\wt{\frs}_{0},\wt{\frs}_{1})=(\frs(\wh{\eta},r),\frs(\eta,r))$, for which we immediately have 
\beq
1\geq \frac{\wt{\eta}_{1}}{\wt{\eta}_{0}}=\frac{\eta}{\wh{\eta}}\geq \frac{\eta_{t_{k+1}}}{\eta_{t_{k}}}\gtrsim N^{-\epsilon_{0}/100},\eeq
where the second and third inequalities are due to $\eta\geq \eta_{t_{k+1}}$ and \eqref{eq:eta_rat}. Thus we only need to check $\wt{\frs}_{1}\geq N^{-\epsilon_{0}/200}\wt{\frs}_{0}$. Indeed, we have $1\geq(\wt{\frs}_{1}/\wt{\frs}_{0})\gtrsim (\eta/\wh{\eta})^{2/3}\geq (\eta_{t_{k+1}}/\eta_{t_{k}})^{2/3}\gtrsim N^{-\epsilon_{0}/150}$, where the only non-trivial inequality $(\wt{\frs}_{1}/\wt{\frs}_{0})\gtrsim (\eta/\wh{\eta})^{2/3}$ follows from Lemma \ref{lem:MDE_asymp} as 
\beq
\frac{\wt{\frs}_{1}}{\wt{\frs}_{0}}=\frac{\frs(\eta,r)}{\frs(\wh{\eta},r)}=\frac{\eta}{\wh{\eta}}\cdot\frac{\rho^{r\caA}(\ii\wh{\eta})}{\rho^{r\caA}(\ii\eta)}\sim\begin{cases}
	\dfrac{\eta}{\wh{\eta}}\cdot\dfrac{\brkt{\absv{r\caA}^{-2}-1}^{1/2}+\wh{\eta}^{1/3}}{\brkt{\absv{r\caA}^{-2}-1}^{1/2}+\eta^{1/3}} & \text{if }\brkt{\absv{r\caA}^{-2}}\geq 1, \\
	\dfrac{\brkt{1-\absv{r\caA}^{-2}}+\eta^{2/3}}{\brkt{1-\absv{r\caA}^{-2}}+\wh{\eta}^{2/3}}& \text{if }\brkt{\absv{r\caA}^{-2}}<1.
\end{cases}
\eeq
Consequently, \eqref{eq:llawzag} (with $N^{\gamma}$ replaced by $N^{\epsilon_{0}+2(k+1)\xi}$) holds with very high probability on the scale $\mathfrak{s}(\eta,r)$ uniformly over $s\in[0,\wh t]$. By Proposition \ref{prop:new zag} now the local law holds at $(\eta,r)=(\wt{\eta}_{1},r)$ with error $\epsilon_{0}+2(k+1)\xi$. Union bounding over a sufficiently dense grid of $\caD_{k+1}\setminus\caD_{k}$ completes the proof. 
\end{proof}

	\subsection{Proof of Proposition~\ref{pro:zignew}}
	\label{sec:zig}
	In what follows we will prove both the average local law and the isotropic local law  for matrices with some Gaussian component. 
	For simplicity we present first the proof of with $B=I$ and the isotropic law, then explain the short adjustment for general $B$. Fix $0<\gamma<\epsilon/10$, an arbitrary small $\xi>0$, $\eta_0>0$, $\caA_{0}\in\C^{N\times N}$,  and unit vectors in $\bsu,\bsy\in\C^{2N}$. Define the stopping time
 \begin{equation}\label{eq:X+A zig stopping time}
		\begin{aligned}
			\tau\deq T^*(\caA_{0},\eta_{0})
			&\wedge \inf\left\{t>0:\absv{\brkt{G_t(\ii \eta_t)-M_{ \caA_t}(\ii\eta_t)}}\geq C_* \frac{N^{\gamma+\xi}}{N\eta_t}\right\}	\\
			&\wedge \inf\left\{s>0:\absv{\bsu\adj (G_t(\ii\eta_t)-M_{ \caA_t}(\ii\eta_t))\bsy}\geq C_* N^{\gamma+\xi}\sqrt{\frac{\rho^{\caA_t}(\ii\eta_t)}{N\eta_t}}\right\},
		\end{aligned}
	\end{equation}
	for some large constant $C_*>0$. To simplify notation we will write $G_t\equiv G_t(\ii\eta_t)$, $M_t\equiv M_{ \caA_t}(\ii\eta_t)$, and $T\equiv T^*(\caA_{0},\eta_0)$ throughout this section.  By It\^o's formula
	\begin{equation}\label{eq:Ito set up}
		\begin{aligned}
			\dd \brkt{G_t-M_t}&=\left(\frac{1}{2}\brkt{G_t-M_t}+\brkt{G_t-M_t}\brkt{G_t^2}+ \frac{\brkt{G_t^3}}{N}\lone(\beta=1) \right)\dd t \\
			&\quad-\frac{1}{N^{3/2}}\sum_{i,j\in \bbrktt{2N}}\left[\left(G_t^2 \right)_{ij}+\left(G_t^2 \right)_{ji}\lone(\beta=1)\right]\dd(\mathfrak{B}_t)_{ij}\\
			&=:\caE_1(t)\brkt{G_t-M_t}\dd t+\frac{\brkt{G_t^3}}{N}\lone(\beta=1)\dd t+\dd\caE_2(t)
		\end{aligned}
	\end{equation} where $\beta=1$ if the entries of $X$ are real and $\beta=2$ if the entries of $X$ are complex, and  \begin{equation*}
		\mathfrak{B}_t=\begin{pmatrix}
			0& B_t\\
			B_t^*&0
		\end{pmatrix}.
	\end{equation*} Similarly, we have 
	\beq\begin{aligned}\label{eq:Ito_x}
		\dd\bsu\adj (G_t-M_t))\bsy=&\left(\frac{1}{2}\bsu\adj (G_t-M_t)\bsy+\brkt{G_t-M_t}\bsu\adj G_t^{2}\bsy+\frac{1}{N}\bsu\adj G_t^{3}\bsy\lone(\beta=1)\right)\dd t\\
		&\qquad-\frac{1}{N^{1/2}}\sum_{i,j}\left[(G_t\bsy\bsu\adj G_t)_{ji}+(G_t\bsy\bsu\adj G_t)_{ij}\lone(\beta=1)\right]\dd (\mathfrak{B}_t)_{ij}	\\
		=&:\frac{1}{2}\bsu\adj(G_t-M_t)\bsy\dd t+\caE_{1}(t,\bsu,\bsy)\dd t+\dd\caE_{2}(t,\bsu,\bsy).
	\end{aligned}\eeq
	
	Recall the short--hand notation $\rho_t=\rho^{\caA_t}(\ii\eta_t)$.  We start with the estimate of the stochastic terms in \eqref{eq:Ito set up}--\eqref{eq:Ito_x}. In the following we will often use the simple integration rule
	\begin{equation}
	\label{eq:intrul}
	\int_0^t\frac{\rho_s}{\eta_s^p}\,\dd s\lesssim \frac{1}{\eta_t^{p-1}},
	\end{equation}
	with $p>1$, even if not stated explicitly. This follows from simple computations using \eqref{eq:local law:eta dynamics}.
	
	\begin{lem}\label{lem:E_2 bounds}
		Uniformly over $t$ in $[0,T]$, and uniformly in $\bsu,\bsy$ deterministic unit vectors, we have
		\begin{equation}\label{eq:local:E2 average bound}
			\left|\int_{0}^{t\wedge \tau}\dd\caE_{2}(s)\right| \prec \frac{1}{N\eta_{t\wedge \tau}},  
		\end{equation} and \begin{equation}\label{eq:local:E2 isotropic bound}
			\left|\int_{0}^{t\wedge\tau} \dd\caE_{2}(s,\bsu,\bsy)\right|\prec \sqrt{\frac{\rho_{t\wedge\tau}}{N\eta_{t\wedge \tau}}}.
		\end{equation}
	\end{lem}
	
	\begin{proof}
		From the definition of $\tau$ we have that \begin{equation}
		\label{eq:imgb}
			\brkt{\im G_t(\ii {\eta}_{t\wedge\tau})}\leq\left(1+N^{-\epsilon/2} \right)\rho_{t\wedge \tau}.
		\end{equation} The quadratic variation in \eqref{eq:Ito set up} satisfies \begin{equation}\label{eq:local law:average quad var}
				\brkt{\dd\caE_{2},\dd\caE_{2} }
				=\frac{1}{N^2}\brkt{G_t^2(G_t^*)^2} \dd t\leq \frac{1}{N^2\eta_t^2}\frac{\brkt{\im G_t}}{\eta_t}\dd t.
		\end{equation} Fix $t\in[0,T]$, applying the Burkholder-Davis-Gundy (BDG) inequality to \eqref{eq:local law:average quad var}, by \eqref{eq:imgb}
		\begin{equation}
		\label{eq:local: single time E2 bound}
			\max_{\wt t\in [0,t]}\left|\int_{0}^{\wt t\wedge \tau}\dd\caE_{2}(s) \right|\prec \frac{1}{N} \left(\int_{0}^{t\wedge\tau}\frac{\rho_s}{\eta_s^3}\dd s \right)\lesssim \frac{1}{N^2\eta_{t\wedge\tau}^2},
		\end{equation} 
		 where in the last inequality we used \eqref{eq:intrul}.  This bound can then be extended uniformly over $[0,T]$ by taking a polynomial sized net of times in $[0,T]$ and union bounding over the very high probability events on which \eqref{eq:local: single time E2 bound} holds.

		We now move on to the proof of \eqref{eq:local:E2 isotropic bound}. For fixed $t\in[0,\tau]$ we have \begin{equation}\label{eq:local:iso E2 quad var}
			\begin{aligned}
				\brkt{\dd\caE_2(t,\bsu,\bsy),\dd\caE_{2}(t,\bsu,\bsy)} &=\frac{\left(\bsu^*\im G_t\bsy \right)^2}{N\eta_t^2}\dd s&\leq\frac{1}{N\eta_t^2}\left(\rho_t+N^{-\epsilon/2}\sqrt{\frac{\rho_t}{N\eta_t}} \right)^2 \dd t\lesssim \frac{\rho_t}{N\eta_t^2} \dd t,
			\end{aligned} 
		\end{equation} where we use \eqref{eq:M boundedness} in the second inequality. The proof then follows completely analogously to $\caE_2(s)$ using the BDG inequality. 
	\end{proof}
	
	We now move on to controlling the size of the drift terms $\caE_{1}(t)$, $N^{-1}\brkt{G_t^3}$, and $\frac{1}{2}\bsu^*(G_t-M_t)\bsy+\caE_1(t,\bsu,\bsy)$. 
	
Using a Schwarz inequality
\[
\big|\brkt{G_s^3}\big|\le \brkt{G_sG_s^*}^{1/2}\brkt{G_s^2(G_s^2)^*}^{1/2}\le \lVert G_s\rVert\brkt{G_sG_s^*}\le \frac{1}{\eta_s}\brkt{G_sG_s^*}
\]
followed by the Ward identity $G_sG_s^*=\Im G_s/\eta_s$, we estimate
\begin{equation}
	\label{eq:gcubeb}
	\int_{0}^{t\wedge \tau}\Absv{\frac{\brkt{G_s^3}}{N}} \dd s\le \int_0^s  \frac{\brkt{\im G_s}}{N\eta_s^2}\,\dd s\lesssim \frac{1}{N\eta_t}.
\end{equation}
We also note that
\beq\begin{aligned}
\label{eq:bounde1}
		\absv{\caE_{1}(t)}&\leq \frac{1}{2}+\frac{\brkt{\im G_t} }{{\eta}_t}\leq \frac{1}{2}+\left(1+N^{-\epsilon/2} \right)\frac{\rho_t}{{\eta}_t}= -\frac{1+O(N^{-\epsilon/2})}{{\eta}_t}\frac{\dd{\eta}_t}{\dd t}=:\wt{\caE}_{1}(t).
	\end{aligned}\eeq Let $\Omega$ be the very high probability event  on which \eqref{eq:incond0} and \eqref{eq:local:E2 average bound}--\eqref{eq:local:E2 isotropic bound} hold.
	Integrating \eqref{eq:Ito set up} with respect to time, and using \eqref{eq:incond0} to estimate the initial condition, we have that on $\Omega$
	\begin{equation}
	\label{eq:befgronw}
		\absv{\brkt{G_{t\wedge\tau }(\ii{\eta}_{t\wedge\tau}) -M_{t\wedge\tau} } }\leq\int_{0}^{t\wedge\tau}\wt{\caE}_{1}(s)\absv{\brkt{G_s(\ii{\eta}_s) -M_s}}\dd s+\frac{N^{\gamma}(\log N)^2}{N\eta_{t\wedge \tau}}.
	\end{equation} Applying a standard integral Gr\"onwall's inequality we get
	\begin{equation}\label{eq:local law: zig Gronwall}
		\absv{\brkt{G_{t\wedge\tau }(\ii{\eta}_{t\wedge\tau}) -M_{t\wedge\tau} } }\leq \int_{0}^{t\wedge \tau}\left(\frac{N^{\gamma}(\log N)^\alpha}{N\eta_{s}}\right) \wt{\caE}_1(s) \exp\left[\int_{s}^{t\wedge\tau}\wt{\caE}_{1}(u)\dd u \right]\dd s.
	\end{equation} Notice from the definition of $\wt{\caE}_1(s)$ that
	\begin{equation}
	\label{eq:asympgro}
		\int_{s}^{t\wedge\tau}\wt{\caE}_{1}(u)\dd u=(1+O(N^{-\epsilon/2}))\log\frac{{\eta}_{s}}{{\eta}_{t\wedge\tau}},
	\end{equation} and hence combining this with \eqref{eq:local law: zig Gronwall} we obtain \begin{equation}\label{eq:local:improv av zig}
		\absv{\brkt{G_{t\wedge\tau }(\ii{\eta}_{t\wedge\tau}) -M_{t\wedge\tau} } }\lesssim \frac{N^{\gamma+\xi}}{N\eta_{t\wedge \tau}},
	\end{equation}
    for any arbitrary small $\xi>0$, giving the desired result.
	
	On the event $\Omega$ we use \eqref{eq:local:iso E2 quad var} and \eqref{eq:local:improv av zig} to estimate $\caE_{1}(s,\bsu,\bsy)$ uniformly over $s\in[0,\tau]$ as \begin{equation}
		\absv{\caE_1(t,\bsu,\bsy)}\lesssim \left|\brkt{G_t-M_t}\bsu^*G_t^2\bsy \right|\lesssim\frac{N^{\gamma+\xi}}{N\eta_t^2}\rho_t.
	\end{equation} We then integrate \eqref{eq:Ito_x} to conclude that \begin{equation}\label{eq:local:iso improv}
			\absv{\bsu^*(G_{t\wedge\tau}-M_{t\wedge\tau})\bsy }
			\lesssim N^{\gamma+\xi}\sqrt{\frac{\rho_{t\wedge\tau}}{N\eta_{t\wedge\tau}}},
	\end{equation}
	where we used Schwarz inequality and Ward identity to estimate $|\bsu^*G_t^2\bsy|\le \rho_t/\eta_t$.
	 Combining \eqref{eq:local:improv av zig}--\eqref{eq:local:iso improv} we thus conclude that  $\tau=T$ holds on $\Omega$.
		
	For general $B$ (with $\beta=2$ for simplicity of presentation) \eqref{eq:Ito set up} becomes \begin{equation}\label{eq:Ito for gen B}
	\begin{aligned}
		\dd \brkt{B(G_t-M_t)}&=\left(\frac{1}{2}\brkt{B(G_t-M_t)}+\brkt{G_t-M_t}\brkt{BG_t^2} \right)\dd t-\frac{1}{N^{3/2}}\sum_{i,j\in \bbrktt{2N}}\left(G_tBG_t \right)_{ij}\dd(\mathfrak{B}_t)_{ij}.\\
	\end{aligned}
\end{equation}	Comparing the form of \eqref{eq:Ito for gen B} with \eqref{eq:Ito set up} one can see that the proof is analogous to the proof used to obtain \eqref{eq:local law: zig Gronwall}, \nc with \eqref{eq:local:improv av zig} used as an input to control the term $\brkt{G_t-M_t}\brkt{BG_t^2}$.
	\qed
	
	\medskip

	We conclude this section with the proof of \eqref{eq:imprllaw}.
	
	\begin{proof}[Proof  of \eqref{eq:imprllaw}]
We first recall some basic facts about the density $\rho^{A-z}$ of the Hermitization of $A+X-z$ (see \eqref{eq:rho_def} for its definition). If $z$ is outside of the support of $\rho_{A+\bsx}$, then the density of the Hermitization, $\rho^{A-z}(\cdot+\ii 0)$, has a gap around zero, which we denote by $\Delta=\Delta^z:=[\mathfrak{e}_-^z,\mathfrak{e}_+^z]$. We point out that the gap is symmetric around zero, i.e. $\mathfrak{e}_+^z=-\mathfrak{e}_-^z$; we keep using  both $\mathfrak{e}_-^z,\mathfrak{e}_+^z$ to keep a notation similar to \cite{Alt-Erdos-Kruger2020}. To keep the presentation simpler, in the remainder of the proof $A,z$ are fixed hence we will often omit the $A$-- and $z$--dependence of various quantities. We recall that the limiting density of the eigenvalues for $\re w\in\Delta$ is so that \cite[Theorem 7.2]{Alt-Erdos-Kruger2020} (in the regime $\kappa(w)\lesssim |\Delta|$):
\begin{equation}
\rho^{A-z}(w)\sim \frac{\im w}{\kappa(w)^{1/2}|\Delta|^{1/6}},
\end{equation}
where $\kappa(w)=\kappa^{A-z}(w):=\mathrm{dist}(w,\supp(\rho^{A-z}))$. Throughout the proof we assume $|\Delta|\gg N^{-3/4}$ and $\kappa(w)\sim \kappa(\re w)$ (and so that $\kappa(w)\lesssim |\Delta|$),
since in the opposite regime, when $|\Delta|\lesssim N^{-3/4}$ or $\eta=\im w\gg \kappa(\re w)\equiv\mathrm{dist}(\re w,\supp(\rho^{A-z}))$, then a bound as in \eqref{eq:imprllaw} for all $w$ with $\re w\in\Delta$ follows by a proof similar to the one performed in Proposition~\ref{pro:zignew} plus a comparison argument as in Section~\ref{sec:comparison}.

Even if our goal is to prove \eqref{eq:imprllaw} only on the imaginary axis, we need to control the resolvent for all $\re w\in \Delta$ to make sure that the Hermitization of $A+X-z$ does not have any outliers “well inside" $\Delta$. For this reason, fix a small $\epsilon_0>0$ and define
\begin{equation}
\label{eq:defwellinside}
\Delta^{(\epsilon_0)}:=\left[\mathfrak{e}_-+\frac{N^{\epsilon_0}|\Delta|^{1/9}}{N^{2/3}}, \mathfrak{e}_+-\frac{N^{\epsilon_0}|\Delta|^{1/9}}{N^{2/3}}\right].
\end{equation}
Then, for fixed $z,A$, we will prove
\begin{equation}
\label{eq:imprllaw1}
\big|\langle G^z(w)-M_{A-z}(w)\rangle\big|\prec \frac{1}{N\sqrt{\kappa(w)|\im w|}},
\end{equation}
uniformly on the domain
\begin{equation}
\label{eq:deffinaldomain}
\mathcal{D}=\mathcal{D}^{A,z}:=\{\re w\in \Delta^{(\epsilon_0)} : N\sqrt{\kappa(w)|\im w|}\rho^{A-z}(w)\ge N^\epsilon\},
\end{equation}
for a fixed small $\epsilon>0$ and for matrices with an (almost) order one Gaussian component. Note that in \eqref{eq:imprllaw1}--\eqref{eq:deffinaldomain} the spectral parameter $w$ and $z,A$ play a different role. In fact $z$ and $A$ are fixed while $w$ varies over the domain $\mathcal{D}$. This is the core part of the proof. We then briefly explain how to remove this Gaussian component, i.e. how to prove \eqref{eq:imprllaw1} for any i.i.d. matrix. The purpose of this weaker bound is  to exclude outlier eigenvalues away from the support of $\rho^{A-z}$ by a distance much larger than the fluctuation scale $\eta_{\mathrm{f}}$ from \eqref{etaf}. Once outliers are excluded, a simple contour integral representation will give \eqref{eq:imprllaw} in its original form with $\kappa$ and general $B$. We point out that the error term in the rhs. of \eqref{eq:imprllaw1} is what we naturally get in the estimates below along the flow \eqref{eq:Ito set up}. For our purpose any estimate $\ll 1/(N|\im w|)$ would do the job of excluding the presence of outliers.

To keep the presentation shorter, in the following we denote $\eta:=\im w$ and assume $\eta>0$. When $\eta/\rho^{A-z}(w)\ge N^{-\delta}$, for some very small fixed $\delta>0$, then \eqref{eq:imprllaw1} immediately follows from \cite[Theorem 2.1]{Erdos-Kruger-Schroder2019} for any i.i.d. matrix. We now show that \eqref{eq:imprllaw1} can be propagated to smaller $\eta$'s using the flow \eqref{eq:Ito set up} at the expense of adding a Gaussian component. This proof is similar to the proof of Proposition~\ref{pro:zignew} (as well as to \cite[Section 4]{Adhikari-Huang2020}), so we only present the main differences. 

From now on we fix an order one time $0\le T<1$, the terminal time of the flow \eqref{eq:Ito set up}, and we will parametrize the flow by $s$, i.e.  we will consider times $0\le s\le T$. Recall that $z,A$ are fixed above \eqref{eq:imprllaw1}, and consider any $w\in \mathcal{D}$.
By standard ODE theory there exist  $\caA_0$ and $w_0$ such that $\wh\kappa_0\deq \dist(w_{0},\supp\rho_{\absv{\caA_{0}+\bsx}})
\gtrsim N^{-\delta}$ and that the characteristics $\partial_s \caA_s=-\caA_s/2$ and
\begin{equation}
\label{eq:chargap}
\partial_s w_s=-\brkt{M_{\caA_{s}}(w_{s})}
	-\frac{w_s}{2},
\end{equation}
with initial condition $\caA_0$ and $w_0$,
 end up at  $\caA_T=A-z$ and $w_T=w$
after time $T$. To keep the presentation shorter, we define the time--dependent versions of $\kappa$ and $\rho^{A-z}$ naturally as $\kappa_{s}\deq \kappa^{\caA_{s}}$ and $\rho_{s}\deq \rho^{\caA_{s}}$. Then, writing $\Delta_{s}:=[\mathfrak{e}_-(s),\mathfrak{e}_+(s)]$ for the gap of $\rho^{A_{s}-z_{s}}$ around the origin and $\wh{\kappa}_{s}\deq\kappa_{s}(w_{s})$,
we have $\wh\kappa_T\lesssim |\Delta_T|$ and so by the properties of the characteristics we have $\wh\kappa_s\lesssim |\Delta_s|$, for any $0\le s\le T$. We point out that the flow $(\caA_s, w_s)$ depends on $A,z,w$, however, to keep the presentation simple, we omit this from the notation. In particular, throughout the proof, all the quantities denoted by hat are functions of $w$, even if not stated explicitly, as a consequence of the fact that $w_s$ depends on $w$ via the condition $w_T=w$. All the following estimates holds uniformly $w=w_T\in\mathcal{D}_T$.

Recall $\epsilon_0>0$ from \eqref{eq:defwellinside}, we introduce the function 
\begin{equation}
\label{eq:deff}
f(s)= f_{\epsilon_0}(s):=\frac{1}{|\Delta_s|^{1/3}}\left[\mathfrak{c}\left(T-s\right)\vee \frac{|\Delta_s|^{2/9} N^{\epsilon_0}}{N^{1/3}}\right]^2,
\end{equation}
for some small $N$--independent $\mathfrak{c}>0$ we will choose shortly.
It is well known \cite[Lemma 4.2]{Cipolloni-Erdos-Kruger-Schroder2019} that the  length of $\Delta_s$ is $|\Delta_s|\sim (t_*-s)^{3/2}$, where $t_*$ denotes the time at which the Hermitized density $\rho_{t_{*}}$
has an exact cusp. Note that $|\Delta_T|\gg N^{-3/4}$ and choosing $\mathfrak{c}$ sufficiently small ensures $f(s)<|\Delta_s|/2$. From now on, with a slight abuse of notation, by $\Delta_s$ we denote the gap as well as its length. In the current proof we always consider
terminal times $T$ such that
$T< t_*$, i.e. the flow terminates before the gap closes. We also introduce the time--dependent domains
\begin{equation}\label{eq:caDs_def}
\mathcal{D}_s:=\{\re w\in \Delta_{s}^{(\epsilon_0)} :  N\sqrt{\kappa_s(w)\eta}\rho_s(w)\ge N^\epsilon\},
\end{equation}
where $\Delta_s^{(\epsilon_0)}:=[\mathfrak{e}_-(s)+f_{\epsilon_0}(s),\mathfrak{e}_+(s)-f_{\epsilon_0}(s)]\subset \Delta_s$. Note that by the definition of $f(s)$ and $\Delta_s^{(\epsilon_0)}$ it follows that $\Delta_T^{(\epsilon_0)}=\Delta^{(\epsilon_0)}$, and so that $\mathcal{D}_T=\mathcal{D}$. Furthermore, an easy calculation using the monotonicity of the characteristics shows that $\mathfrak{f}_{s, T}^{-1}(\mathcal{D}_T)\subset \mathcal{D}_s$, where  $\mathfrak{f}_{s, T}:\C_+\to \C_+$, $s\le T$, denotes the  propagator of \eqref{eq:chargap}, e.g. $\mathfrak{f}_{s, T}(w_s)=w_T$. 

We now state some properties of $f(s)$ and of $\widehat{\kappa}_s=\kappa_s(w_s)$ for $w_s\in \mathfrak{f}_{s, T}^{-1}(\mathcal{D}_T)$. From $\Delta_s\sim (t_*-s)^{3/2}$ and \cite[Eq. (4.15e)]{Cipolloni-Erdos-Kruger-Schroder2019} it directly follows that for $\alpha \in [0,2/3]$ we have
\begin{equation}
\label{eq:changegap}
\Delta_s^\alpha-\Delta_T^\alpha\sim \frac{T-s}{\Delta_s^{2/3-\alpha}}, \qquad 0\le s\le T\le t_*.
\end{equation}
Then, by \eqref{eq:deff} and \eqref{eq:changegap} used for $\alpha=1/9$, it follows that
\begin{equation}
\label{eq:ineqf}
\sqrt{f(s)}\le \sqrt{f(T)}+\mathfrak{c}\frac{(T-s)}{\Delta_s^{1/6}}, \qquad 0\le s\le T\le t_*.
\end{equation}
Furthermore, by explicit, but fairly tedious, computations from \eqref{eq:chargap} we find that for $w_s\in \mathfrak{f}_{s, T}^{-1}(\mathcal{D}_T)$ there exists $C>0$ such that
\begin{equation}
\label{eq:kappaineq}
\sqrt{\wh{\kappa}_s}\ge \sqrt{\wh\kappa_T}+C\frac{(T-s)}{\Delta_s^{1/6}}, \qquad 0\le s\le T\le t_*.
\end{equation}
Next, combining \eqref{eq:ineqf} and\eqref{eq:kappaineq} and  choosing $\mathfrak{c}\le C/2$, we readily have
\begin{equation}
\label{eq:implowerb}
\wh\kappa_s-f(s)\gtrsim \sqrt{\wh\kappa_s}\frac{T-s}{\Delta_s^{1/6}}, \qquad 0\le s\le T\le t_*.
\end{equation}

We now consider the evolution of the resolvent $G_s(w)=(H_{\caA_{s}+X_{s}}-w)^{-1}$ (recall the definition from \eqref{eq:G_t}) along the flow \eqref{eq:Ito set up}.
Define the stopping time (here $M_s(w):=M_{A_{s}-z_s}(w)$):
\begin{equation}
\label{eq:defstoptime}
\tau:=\inf \left\{s\ge 0 : \, \exists \, w\in \mathcal{D}_s \, \mathrm{s.t.}\, \big|\langle G_s(w)-M_s(w)\rangle\big|=\frac{N^\gamma (\log N)^2}{N\sqrt{\kappa_{s}(w)\eta}}\right\}\wedge T,
\end{equation}
for some $\gamma\le \epsilon/10$, where $\epsilon$ is from the definition of $\caD_{s}$ in \eqref{eq:caDs_def}. Our goal is to show that $\tau=T$. Given the definition of $\tau$, we claim that for any $s\le\tau$ there are no eigenvalues in $\Delta_{s}^{(\epsilon_0)}$. In fact, assume by contradiction that there is an eigenvalue $\lambda\in \Delta_s^{(\epsilon_0)}$, and let  $w_*:=\lambda+\ii\eta$, with $\eta:=N^{\epsilon_0/10}\Delta_s^{1/6}/N^{1/3}$. Then, $\langle \Im G_s(w_*)\rangle\gtrsim 1/(N\eta)$ (note that  $\re w_*\in\Delta_s^{(\epsilon_0)}$). On the other hand, by \eqref{eq:defstoptime}, it follows that for any $\widetilde{w}\in\mathcal{D}_s$, $s\le \tau$,  we have
\begin{equation}
\label{eq:imgbound}
\langle \im G_s(\widetilde{w})\rangle \le\rho_s(\widetilde{w})+\frac{N^\gamma(\log N)^{2}}{N\sqrt{\kappa_s(\widetilde{w})\im \widetilde{w}}}\lesssim \rho_s(\widetilde{w})\sim \frac{\eta}{\kappa_s(\widetilde{w})^{1/2}\Delta_s^{1/6}}.
\end{equation}
Using this for $\widetilde{w}:=w_*\in \mathcal{D}_s$, we have
\begin{equation}
\langle \im G_s (w)\rangle\lesssim \frac{\eta}{\kappa_s(w)^{1/2}\Delta_s^{1/6}}\le \frac{N^{-\epsilon_0/10}}{N\eta}.
\end{equation}
This leads to a contradiction.

We introduce the function $d_s(w):=\mathrm{dist}(w,\mathrm{Spec}(\absv{A_{s}+X_{s}-z_{s}}))$, and denote $\widehat{d}_s:=d_s(w_s)$ to measure the distance of the spectral parameter $w_s$ from the spectrum of the Hermitization of $A+X-z$. We point out that $\widehat{d}_s$ is a function on $\caD_{T}\equiv\caD$, since $w_{s}$ depends on the final point $w_{T}\in\caD_{T}$ (see the discussion below \eqref{eq:chargap}).
Note that with very high probability we have
\begin{equation}
\label{eq:lbdist}
\widehat{d}_s \gtrsim \sqrt{\widehat{\kappa}_s}(T-s)/\Delta_s^{1/6}+\eta_s,  \qquad s\le \tau,
\end{equation}
which follows from the fact that there are no eigenvalues in $\Delta_s^{(\epsilon_0)}$ and from \eqref{eq:implowerb}. Here we defined $\eta_s:=\im w_s$. Even though $\eta_s$ is also a function of $w_T$, we do not use the hat notation for it since it depends on it in a simple way and to keep the notation simpler.

We now explain how to estimate all the terms in the rhs. of \eqref{eq:Ito set up} by
$N^\epsilon /(N\sqrt{\widehat{\kappa}_T\eta_T})$. Using \eqref{eq:imgbound} for $\widetilde{w}=w_s$ and \eqref{eq:lbdist}, we can estimate the quadratic variation of the stochastic term in \eqref{eq:Ito set up} by
\begin{equation}
\begin{split}
\int_0^{T\wedge\tau} \frac{\langle \im G_s \im G_s\rangle}{N^2\eta_s^2}\, \rmd s \lesssim \int_0^{T\wedge \tau} \frac{\langle \im G_s\rangle}{N^2\eta_s \widehat{d}_s^2}\, \rmd s &\lesssim \int_0^{T\wedge \tau} \frac{1}{N^2\sqrt{\widehat{\kappa}_s}\Delta_s^{1/6}\widehat{d}_s^2}\, \rmd s \\
&\lesssim \int_0^{T\wedge \tau} \frac{1}{N^2\widehat{\kappa}_s^{1/2}\Delta_s^{1/6}[\widehat{\kappa}_s^{1/2}(T-s)/\Delta_s^{1/6}+\eta_s]^2}\,\rmd s\\
&\lesssim \frac{1}{N^2\widehat{\kappa}_T\eta_T}.
\end{split}
\end{equation}
In the first inequality  we used that
\[
\langle \im G_s \im G_s\rangle\le \lVert \im G_s\rVert \langle \im G_s\rangle \le \frac{\eta_s}{\widehat{d}_s^2} \langle \im G_s\rangle,
\]
while in the last inequality we used that $\eta_s\gtrsim \eta_T$, $\widehat{\kappa}_s\gtrsim\widehat{\kappa}_T$, $\Delta_s\gtrsim\Delta_T$, $\widehat{\kappa}_s^{1/2}/\Delta_s^{1/6}\gtrsim \widehat{\kappa}_T^{1/2}/\Delta_T^{1/6}$, for any $s\le T$. 
Finally,  by the BDG inequality, we obtain the desired bound for the stochastic term in \eqref{eq:Ito set up}. Proceeding in a similar way we obtain the same bound for the term $\langle G_s^3\rangle/N$ in \eqref{eq:Ito set up} as well; we omit this for brevity.
	
We thus get (cf. \eqref{eq:befgronw})
\begin{equation}\label{eq:lastb-1}
\absv{\brkt{G_{T\wedge\tau }(w_{T\wedge\tau}) -M_{T\wedge\tau}(w_{T\wedge\tau})} }\leq\int_{0}^{T\wedge\tau}\wt{\caE}_{1}(s)\absv{\brkt{G_s(w_s) -M_s(w_{s})}}\dd s+\frac{N^{\gamma}(\log N)^2}{N \sqrt{\widehat{\kappa}_{T\wedge \tau}\eta_{T\wedge \tau}}},
\end{equation}
with $\wt{\caE}_{1}(s)$ from \eqref{eq:bounde1}. Finally, using Gronwall inequality and \eqref{eq:asympgro}, the desired result follows from (here we use the notation $\widehat{\rho}_s:=\rho_s(w_s)$)
\begin{equation}
\label{eq:lastb}
\int_0^{T\wedge\tau} \frac{\widehat{\rho}_s}{\eta_s} \cdot \frac{\eta_s}{\eta_{T\wedge\tau}}\cdot \frac{1}{N\sqrt{\widehat{\kappa}_s\eta_s}}\,\rmd s \lesssim \frac{1}{N\widehat{\kappa}_{T\wedge\tau}^{1/2}\Delta_{T\wedge\tau}^{1/6}}\int_0^{T\wedge\tau} \frac{1}{\sqrt{\widehat{\kappa}_s\eta_s}}\,\rmd s\lesssim \frac{\log N}{N\sqrt{\widehat{\kappa}_{T\wedge\tau}\eta_{T\wedge\tau}}},
\end{equation}
where each factor in the integrand on the leftmost side of \eqref{eq:lastb} controls $\wt{\caE}_{1}(s)$, $\exp(\int_{s}^{T\wedge\tau}\wt{\caE}_{1}(s))$, and $\absv{\brkt{G_{s}-M_{s}}}$, in that order. For the proof of the second inequality above, we used the following side
calculation. Recall 
 \eqref{eq:kappaineq}, that without loss of generality we assume $\widehat{\kappa}_s\lesssim \Delta_s$, and recall that for any $0\le s\le t'$ we have
 \[
\widehat{\rho}_s\sim \widehat{\rho}_{t'}, \qquad\quad \eta_s\sim\eta_{t'}+\wh\rho_s(t'-s). 
\]
Using these relations, for any $t'\le T$ we have
\begin{equation}
\begin{split}
\label{eq:tedcomp}
\int_0^{t'}\frac{1}{\sqrt{\widehat{\kappa}_s\eta_s}} \,\rmd s&\lesssim \int_0^{t'}\frac{1}{[\sqrt{\widehat{\kappa}_{t'}}+(t'-s)/\Delta_s^{1/6}][\eta_{t'}+\widehat{\rho}_s(t'-s)]^{1/2}}\,\rmd s \\
&\lesssim \Delta_{t'}^{1/6}\int_0^t\frac{1}{[\Delta_{t'}^{1/6}\sqrt{\widehat{\kappa}_{t'}}+(t'-s)][\eta_{t'}+\widehat{\rho}_s(t'-s)]^{1/2}}\,\rmd s \\
&\quad+ \int_0^{t'}\frac{t'-s}{\Delta_s^{1/2}[\Delta_{t'}^{1/6}\sqrt{\widehat{\kappa}_{t'}}+(t'-s)][\eta_{t'}+\widehat{\rho}_s(t'-s)]^{1/2}}\,\rmd s \\
&\lesssim \frac{\Delta_{t'}^{1/6}}{\sqrt{\eta_{t'}}}+\int_0^{t'}\frac{1}{\Delta_s^{1/2}[\eta_{t'}+\widehat{\rho}_s(t'-s)]^{1/2}}\,\rmd s\lesssim \frac{\Delta_{t'}^{1/6}\log N}{\sqrt{\eta_{t'}}}.
\end{split}
\end{equation}
In the second inequality we also used \eqref{eq:changegap} for $\alpha=1/6$, and in the last inequality we used that
\[
\int_0^{t'}\frac{1}{\Delta_s^{1/2}[\eta_{t'}+\widehat{\rho}_s(t'-s)]^{1/2}}\,\rmd s\le \left(\int_0^{t'}\frac{1}{\widehat{\rho}_s\Delta_s}\,\rmd s\right)^{1/2}\left(\int_0^{t'}\frac{\widehat{\rho}_s}{\eta_{t'}+\widehat{\rho}_s(t'-s)}\, \rmd s\right)^{1/2}\lesssim  \frac{\log N}{\sqrt{\widehat{\rho}_{t'}}\Delta_{t'}^{1/6}},
\]
which follows from $\Delta_s\sim (t_*-s)^{3/2}$ and $\Delta_{t'}\sim(t_*-t')^{3/2}$. Plugging \eqref{eq:tedcomp}, with $t'=T\wedge \tau$, into the lhs. of the second inequality of \eqref{eq:lastb} gives the bound in the rhs. of \eqref{eq:lastb}. 

This proves \eqref{eq:imprllaw1}  for matrices with an order one Gaussian component. We then remove this Gaussian component with a fairly standard moment matching argument (see e.g. \cite{Cipolloni-Erdos-Henheik2023, Erdos-Yau2017, Tao-Vu2011}). More precisely, given an i.i.d. matrix $X$, we can match the first three moments with a matrix having an order one Gaussian component \cite[Lemma 16.2]{Erdos-Yau2017}. Then, proceeding in a similar fashion to the argument in \cite[Section 5]{Cipolloni-Erdos-Henheik2023}\footnote{ That is, given the usual local law $\absv{\brkt{G-M}}\prec 1/(N\eta)$ for deformed i.i.d. matrices and the improved bound $1/(N\sqrt{\kappa\eta})$ for the Gaussian divisible model, we use a four-moment-matching argument to prove $\absv{\brkt{G-M}}\prec 1/(N\sqrt{\kappa\eta})$ for deformed i.i.d. matrices. We also point out that with some additional effort this comparison argument could be replaced by a two--moment matching GFT similar to the one presented in Section~\ref{sec:comparison}.}, we  get that \eqref{eq:imprllaw1} holds for general i.i.d. matrices. 
Given \eqref{eq:imprllaw1}, by a very similar proof, on $\mathcal{D}$, we also get
\begin{equation}
\label{eq:imprllaw2}
\big|\langle (G^z(w)-M_{A-z}(w))B\rangle\big|\prec \frac{1}{N\sqrt{\kappa(w)|\im w|}}
\end{equation}
for general $B$'s with $\lVert B\rVert\lesssim 1$.

We now prove that the error term in the rhs. of \eqref{eq:imprllaw2} can be improved to $1/(N\kappa(w)|\im w|)$, which implies the original \eqref{eq:imprllaw}.
First, from \eqref{eq:imprllaw2}, following the same argument as the one around \eqref{eq:imgbound}, we know that there are no outliers in $\Delta^{(\epsilon_0)}=[\mathfrak{e}_-+\Delta^{1/9}N^{\epsilon_0}/N^{1/3}, \mathfrak{e}_+-\Delta^{1/9}N^{\epsilon_0}/N^{1/3}]$, for any small $\epsilon_0>0$. Then, we point out that \eqref{eq:imprllaw2} easily extends (see e.g. \cite[Appendix A]{Cipolloni-Erdos-Schroder2021}) the usual local law
\begin{equation}
\label{eq:smalllaw}
\big|\langle (G^z(w)-M_{A-z}(w))B\rangle\big|\prec \frac{1}{N\im w}
\end{equation}
uniformly over all $w$ with $\absv{\re w}\leq \fre_{+}-\Delta^{1/9}N^{\epsilon_0/2}/N^{1/3}$, including those with $\im w\ll 1/N$. Lastly, for any fixed $w$ such that $\Re w\in \Delta^{(\epsilon_0)}$, we write (cf. \cite[Remarks 4.4--4.10]{Adhikari-Huang2020})
\begin{equation}
\langle (G^z(w)-M_{A-z}(w))B\rangle=\oint_\Gamma \frac{\langle (G^z(u)-M_{A-z}(u))B\rangle}{u-w}\, \rmd u,
\end{equation}
where $\Gamma$ is a rectangle of height $4\kappa$, with $\kappa=\kappa(w)$, and width $2\kappa-\Delta^{1/9}N^{\epsilon_0}/(2N^{1/3})$ centered at $\Re w$. Note that for any $u\in \Gamma$ we have $|u-w|\gtrsim \kappa$. Thus, using \eqref{eq:smalllaw} for the points $u\in \Gamma$, we get
\begin{equation}
\big|\langle (G^z(w)-M_{A-z}(w))B\rangle\big|=\oint_\Gamma \frac{|\langle (G^z(u)-M_{A-z}(u))B\rangle|}{|u-w|}\, \rmd u\lesssim \frac{1}{\kappa}\int_\Gamma |\langle (G^z(u)-M_{A-z}(u))B\rangle|\,\rmd u\prec \frac{\log N}{N\kappa},
\end{equation}
concluding the proof. We point out that in the last inequality to estimate the integral on
the vertical lines of $\Gamma$ we used that the regime $|\im u|\le N^{-100}$ can easily be removed as a consequence of the fact that there are no outliers in $\Delta^{(\epsilon_0)}$.
	\end{proof}

	\subsection{Proof of Proposition~\ref{prop:new zag}}
	\label{sec:comparison}
	
	Fix $r\geq 1$ and let $\epsilon>0$ be as in \eqref{eq:finscale}. We will prove Proposition \ref{prop:new zag} by first proving the high moment bounds  \begin{equation}\label{eq:local law: zag avaerage}
		\E\left|\left\langle B\left(G^0(\ii\eta)-M_{r\caA}(\ii \eta)\right) \right\rangle\right|^p\lesssim  \E\left|\left\langle B\left(G^{t}(\ii\eta)-M_{r\caA}(\ii \eta)\right) \right\rangle\right|^p+\left(\frac{1}{N\eta} \right)^p,
	\end{equation} and \begin{equation}\label{eq:local law: zag isotropic}
		\E\left| \bsu^*\left(G^0(\ii\eta)-M_{r\caA}(\ii \eta) \right)\bsy \right|^p\lesssim \E\left| \bsu^*\left(G^{t}(\ii\eta)-M_{r\caA}(\ii \eta) \right)\bsy \right|^p+\left( \sqrt{\frac{\rho^{r\caA}(\ii\eta)}{N\eta} } \right)^{p},    
	\end{equation} on the scale $\widetilde{\mathfrak{s}}_{1}$, i.e. for $\eta=\wt{\eta}_{1}$, for $p\in\N$.
	The proof of these two bounds will be presented shortly.
	
	\begin{proof}[Proof of Proposition~\ref{prop:new zag}]

	By assumption, \eqref{eq:llawzag} holds for $G^t-M_{r\mathcal{A}}$. Using this fact and \eqref{eq:local law: zag avaerage}--\eqref{eq:local law: zag isotropic} we immediately obtain the desired bound \eqref{eq:llawzag} for $t=0$ as well.
	\end{proof}

	 Recall that $\wt{\eta}_j$ satisfies $\widetilde{\mathfrak{s}}_{j}=\mathfrak{s}(\widetilde{\eta}_{j},r)$ for $j=0,1$. For the remainder of the section we drop the tilde and write $\eta_{j}\equiv\wt{\eta}_{j}$. Define
	\begin{equation}\label{eq:local law Psi}
		\Psi(z,\eta)\deq\sqrt{\frac{\rho^{r\caA}(\ii\eta)}{N\eta}}.
	\end{equation} We begin with the following lemma which uses monotonicity to control the entries of $G^s(\ii {\eta})$ at ${\eta_1}$ 
	 given a local law at $\eta_0$.   We point out that $t$ will be fixed throughout this section.

	\begin{lem}\label{lem:X+A monotone}
		Fix $s\in [0,t]$.  Then, under the assumptions of Proposition \ref{prop:new zag}, we have
		\begin{equation}\label{eq:X+A monotone}
			\left|\bsu^* G^{s}(\ii\eta_1)\bsu \right| \lesssim\frac{ \eta_0}{\eta_1}\rho^{r\caA}(\ii{\eta_0})\quad\text{and}\quad \left|\bsu^* G^{s}(\ii\eta_1)\bsy \right| \lesssim\frac{ \eta_0}{\eta_1},
		\end{equation} with very high probability uniformly in $\bsu,\bsy \in \C^{2N}$. 
	\end{lem}
	
	\begin{proof}
		We have
		\beq\label{eq:X+A monotone_1}
			\absv{\bsu\adj(G^{s}(\ii \eta_1)-G^{s}(\ii{\eta_0}))\bsy}\leq \int_{\eta_1}^{{\eta_0}}\absv{\bsu\adj G^{s}(\ii \wh{\eta})^{2}\bsy}\dd \wh{\eta}
			\leq \frac{{\eta_0}}{2\eta_1}(\bsu\adj \im G^{s}(\ii{\eta_0})\bsu+\bsy\adj \im G^{s}(\ii{\eta_0})\bsy).
\eeq
		where in the last inequality we used\footnote{We point out that in the remainder of the paper we will often use this inequality even if not stated explicitly.}
		\begin{equation}
		\absv{\bsu\adj G^{s}(\ii \wh{\eta})^{2}\bsy}\le \big(\bsu\adj G^{s}(\ii \wh{\eta})G^{s}(\ii \wh{\eta})^*\bsu\big)^{1/2} \big(\bsy\adj G^{s}(\ii \wh{\eta})G^{s}(\ii \wh{\eta})^*\bsy\big)^{1/2}
		\end{equation}
		which follows by the Cauchy-Schwarz inequality and the Ward identity, then followed by an AM--GM inequality, and that $\wh{\eta}\mapsto \wh{\eta}\im G(\ii\wh{\eta})$ is monotone increasing with respect to the Loewner order. As, by assumption, $\mathfrak{s}_{\mathrm{fin}}<\widetilde{\mathfrak{s}}_0$ and \eqref{eq:llawzag} holds on a scale $\widetilde{\mathfrak{s}}_0$ we obtain from \eqref{eq:finscale} that $\Psi(z,\eta_{0})\lesssim \rho^{r\caA} N^{-\epsilon/2}$ and
		\begin{equation}\label{eq:X+A monotone 2}
			\absv{\bsu\adj(G^{s}(\ii \eta_1)-G^{s}(\ii{\eta_{0}}))\bsy}\leq \frac{{\eta_0}}{\eta_1}\rho^{r\caA}(\ii{\eta_0})\left(1+O\left(N^{-\epsilon/5}\right) \right)
		\end{equation} with very high probability. 
	Similarly, by \eqref{eq:llawzag}, we obtain \begin{equation}\label{eq:X+A monotone 3}
			\left| \bsu^*G^s(\ii{\eta_0})\bsy \right|\leq \frac{{\eta_0}}{\eta_1}N^{\gamma}\left(\log N\right)^{\alpha}\Psi({z},{\eta_0})+\left|\bsu^*M_{r\caA}(\ii{\eta_0})\bsy \right|.
		\end{equation} Noting again by \eqref{eq:finscale} that $\Psi(z,\eta_{0})\lesssim \rho^{r\caA} N^{-\epsilon/2}$, and combining \eqref{eq:X+A monotone 2} with \eqref{eq:X+A monotone 3} we complete the proof. 
	\end{proof}

	 For simplicity we will use the short--hand notation $\wt \rho\equiv\rho^{r\caA}(\ii\eta_0)$. Lemma \ref{lem:X+A monotone} will be our main tool for estimating the entries of $G$ below. At times we will keep track of the $\wt\rho$ term, however often we will apply \eqref{eq:M boundedness} and simply use that $\wt\rho\lesssim 1$.  To complete the proofs of \eqref{eq:local law: zag avaerage} and \eqref{eq:local law: zag isotropic}, we can apply Gr\"onwall's inequality to the following propositions. For simplicity of presentation in the remainder of the section we use the short--hand notation $\delta=\frac{\epsilon_0}{100}$ to denote the step size in our scale. We note that Lemma \ref{lem:X+A monotone} is presented as an inequality holding with high probability. This implies that the inequality also holds in the sense of stochastic domination, which will be more convenient to use in the proofs of Propositions \ref{prop:X+A derivative estimate average} and \ref{prop:X+A derivative estimate isotropic} below. 
	
	\begin{prop}\label{prop:X+A derivative estimate average}
		Define \begin{equation}
			R_s\deq \left\langle B\left(G^s(\ii\eta)-M_{r\caA}(\ii \eta) \right)\right\rangle.
		\end{equation} Let $p> 2$  and $\delta:=\epsilon_0/100$,  then
				\begin{equation}
				\label{eq:pregroav}
			\left|\E\frac{\dd}{\dd s}\absv{R_s}^p \right|\lesssim  N^{\frac{1}{2}+ 10 \delta}\Psi(z,\eta_0)\left[\absv{R_s}^p+\left(\frac{1}{N\eta_1}\right)^p \right],
		\end{equation} for all $0\leq s\leq t$. 
	\end{prop}

	\begin{proof}[Proof of \eqref{eq:local law: zag avaerage}]
	
	By \eqref{eq:pregroav} with $\eta_1=\eta$, using Gronwall inequality we obtain
	\begin{equation}
	\begin{split}
	\E|R_0|^p&\lesssim \exp\left(tN^{\frac{1}{2}+ 10\delta}\Psi(z,\eta_0)\right) \left[\E|R_t|^p+\left(\frac{1}{N\eta}\right)^p\right] \\
	&\le\exp\left(N^{10\delta}\left[\sqrt{\frac{\eta_0}{\rho^{r\mathcal{A}}(\ii\eta_0)}}+\frac{\eta_0}{\sqrt{N}\rho^{r\mathcal{A}}(\ii\eta_0)}\right]\right) \left[\E|R_t|^p+\left(\frac{1}{N\eta}\right)^p\right] \\
	&\le \exp\big(N^{10\delta-\epsilon_0/2}\big) \left[\E|R_t|^p+\left(\frac{1}{N\eta}\right)^p\right],
	\end{split}
	\end{equation}
	where in the first inequality we used that $t\le \widetilde{\mathfrak{s}}_0=\eta_0/\rho^{r\mathcal{A}}(\ii\eta_0)$, and in the second inequality that $\eta_0/\rho^{r\mathcal{A}}(\ii\eta_0)\le N^{-\epsilon_0}$, by assumption. Finally, using $10\delta-\epsilon_0/2=-2\epsilon_0/5$ we obtain $\exp\big(N^{10\delta-\epsilon_0/2}\big)\lesssim 1$, which concludes the proof.
	
		\end{proof}

	\begin{proof}
		 Let $H_s$ be the Hermitization of $r\mathcal{A}+X_s$. Differentiating $R_s=R_s(H_s)$,  by It\^o's lemma we have
		\begin{equation}\label{eq:local law: average zag Ito}
			\E\frac{\dd}{\dd s}|R_s|^{p}=\E\left[-\frac{1}{2}\sum_{a,b}w_{ab}(s)\partial_{ab} |R_s|^{p}+\frac{1}{2}\sum_{a,b,c,d}\kappa_s(ab,cd)\partial_{ab}\partial_{cd} |R_s|^{p} \right],
		\end{equation} where $a,b,c,d\in[2N]$, $w_{ab}(s)$ is the $ab$-th entry of $H_s$, $\kappa_{s}(ab,a_1b_1,\dots)$ are the joint cumulants of $w_{ab}(s),w_{a_1b_1}(s),\dots$, and $\partial_{ab}=\partial_{w_{ab}(s)}$  denotes the directional derivative in the direction $w_{ab}(s)$. For a smooth function $f$ of the entries of $H_{s}$ one can expand expectations in the cumulants:\begin{equation}
		\E w_{ab}f(H)=\sum_{k=0}^L\sum_{a_{1},b_{1},\dots,1_{k}b_{k}}\frac{\kappa_{s}(ab,a_{1}b_{1},\dots,a_{k}b_{k} )}{k!}\E\partial_{a_{1}b_{1}}\cdots\partial_{a_kb_k}f(H)+\Omega_{f}(L),
	\end{equation} for some error $\Omega_{f}(L)$ going to zero as $L\rightarrow\infty$.  Using this cumulant expansion for the first sum in \eqref{eq:local law: average zag Ito} we see that \begin{equation}\label{eq:local law average cumulant expansion}
			\left|\E\frac{\dd}{\dd s}|R_s|^{p} \right|\leq\left|\sum_{k=2}^{L-1}\sum_{a,a_1,\dots a_k,b,b_1,\dots,b_k}\frac{\kappa_s(ab,a_1b_1,\dots,a_kb_k)}{k!}\E\partial_{a_kb_k}\cdots\partial_{ab}|R_s|^{p} \right|+\Omega(L),
		\end{equation} 	where $\Omega(L)\le\left(\frac{1}{N\eta_1} \right)^p$ for  some sufficiently large $L$ independent of $N$. 
		
	 We now first present the estimates for the third order cumulants, which are the most delicate, and then present the ones for higher order cumulants.  For the third order cumulants we aim to control the term \begin{equation*}
			\left|N^{-3/2}\sum_{a,b}\E \partial_{ab}^3\absv{R_s}^p\right|,
		\end{equation*}which after applying the chain rule can be broken up into terms involving derivatives of the form $(\partial_{ab}R_s)^3\absv{R_s}^{p-3}$, $(\partial_{ab}R_s)(\partial_{ab}^2R_s)\absv{R_s}^{p-2}$, and $(\partial_{ab}^3R_s)\absv{R_s}^{p-1}$.  For simplicity of notation, here we neglected the difference between $R_s$ and $\overline{R_s}$ since terms involving one or the other will be estimated in a completely analogous way. 
		
		We begin with the term involving $(\partial_{ab}^3R_s)\absv{R_s}^{p-1}$. First we note \begin{equation}
			\partial_{ab}^3 R_s=\frac{1}{2N}\sum_{x} \partial_{ab}^3(BG)_{xx}, 
		\end{equation} and $\partial_{ab}^3 G_{xx}$ is a sum of terms of the form $(BG)_{xa}G_{bb}G_{aa}G_{bx}$, $(BG)_{xa}G_{ba}G_{bb}G_{ax}$, $(BG)_{xa}G_{ab}^2G_{bx}$, or similar terms with some transposition of $a$ and $b$. We will see that there are essentially two cases, when the product of $G$ entries involves an off diagonal term, and when it does not. We now illustrate how to stochastically dominate each of these representative examples. In all three examples,  as well as in the case when the derivatives hit more $R_s$ or for higher cumulants,  the terms involving $(BG)_{x\cdot}$ will be dealt with by eventually bounding the norm of a column of $BG$ by $\norm{B}$ times the identically indexed column of $G$. To avoid repetition we assume $B=I$ for the remainder of the proof.
		Applying Cauchy--Schwarz to the sum of $x$, followed by the Ward identity and Lemma \ref{lem:X+A monotone}, we obtain (recalling $\eta_0/\eta_1\leq N^{3\delta/2}$)  
		\begin{equation}\label{eq:local:third der abab} 
			\begin{aligned}
				N^{-5/2}\left| \sum_{a,b,x}G_{xa}G_{ab}G_{ab}G_{bx} \right|&\lesssim N^{-5/2}\frac{1}{\eta_1}\sum_{a,b}|G_{ab}|^2\sqrt{\im G_{aa}}\sqrt{\im G_{bb}}\\
				&\prec N^{-3/2}\frac{1}{\eta_1^2}\left(\frac{\eta_0}{\eta_1}\wt \rho\right)^2
				\lesssim \frac{N^{10\delta}\wt \rho^2}{\eta_1 \sqrt{N}}\frac{1}{N\eta_1}.
			\end{aligned}
		\end{equation}
	Similarly, 
	we bound\begin{equation}\label{eq:local:third der babb}
							N^{-5/2}\left| \sum_{a,b,x}G_{xa}G_{ba}G_{bb}G_{ax} \right|
				\prec \frac{N^{10\delta}\wt \rho^{5/2}}{\sqrt{\eta_1}}\frac{1}{N\eta_1}.
		\end{equation}

		For the remaining term, we need to handle the sum over the trace index $x$ more carefully, since we do not have other off-diagonal terms to apply the Ward identity to. Applying Lemma \ref{lem:X+A monotone} and the Ward identity 
		 \begin{equation}\label{eq:local:third der aabb}
			\begin{aligned}
				N^{-5/2}\left| \sum_{a,b,x}G_{xa}G_{aa}G_{bb}G_{bx} \right|&\lesssim N^{-5/2}\sum_{a,b}\absv{G_{aa}}\absv{G_{bb}}\absv{(G^2)_{ba}}\\
				&\prec N^{-5/2}\left(\frac{\eta_0}{\eta_1}\wt \rho \right)^2\sum_{b}\left(N\sum_{a}\absv{(G^2)_{ab}}^2 \right)^{1/2}\\
				&\lesssim \frac{1}{\sqrt{\eta_1}}\left(\frac{\eta_0}{\eta_1}\wt \rho \right)^{5/2}\frac{1}{N\eta_1}.
			\end{aligned} 
		\end{equation}
%
		Applying Young's inequality, from \eqref{eq:local:third der abab}, \eqref{eq:local:third der aabb}, and \eqref{eq:local:third der babb}  we conclude that \begin{equation}
		\label{eq:almthere}
			\left|N^{-3/2}\sum_{a,b} (\partial_{ab}^3R_s)\absv{R_s}^{p-1} \right| \prec  \frac{N^{5\delta}\wt \rho^2}{\sqrt{\eta_1}}\left[\absv{R_s}^{p}+\left(\frac{1}{N\eta_1}\right)^{p} \right] \le N^{1/2+6\delta}\Psi(z,\eta_0)\left[\absv{R_s}^{p}+\left(\frac{1}{N\eta_1}\right)^{p}\right]. 
		\end{equation} 
		
		Moving on to terms involving $(\partial_{ab}R_s)(\partial_{ab}^2R_s)\absv{R_s}^{p-2}$, we again use Ward's identity, Cauchy--Schwarz inequality, and Lemma \ref{lem:X+A monotone} to estimate
		\begin{equation}
			\begin{aligned}
				N^{-3/2}\left|\sum_{a,b}(\partial_{ab}R_s)(\partial_{ab}^2R_s)\absv{R_s}^{p-2} \right|&\lesssim N^{-3/2}\left|\sum_{a,b}\left[\frac{1}{2N}\sum_{x} G_{xa}G_{bx}\right]\left[\frac{1}{2N}\sum_{x} G_{xb}G_{aa}G_{bx}+G_{xa}G_{ba}G_{bx} \right]\absv{R_s}^{p-2} \right|\\
				&\prec \left(\frac{\eta_0}{\eta_1}\wt \rho \right)^{3/2}\left(\frac{1}{N\eta_1}\right)^{2}N^{-3/2}\sum_{a,b}\absv{(G^2)_{ab}}\absv{R_s}^{p-2}\\
				&\lesssim \frac{1}{\sqrt{\eta_1}}\left(\frac{\eta_0}{\eta_1}\wt \rho \right)^{2}\left(\frac{1}{N\eta_1}\right)^{2}\absv{R_s}^{p-2}.
			\end{aligned}
		\end{equation} 
		
		Finally, for the terms involving $(\partial_{ab}R_s)^3\absv{R_s}^{p-3}$ we apply the same approaches to see
		\begin{equation}
				\left| N^{-3/2}\sum_{a,b}\left[\partial_{ab}R_s\right]^{3}\absv{R_s}^{p-3}  \right|
				\prec \left(\frac{\eta_0}{\eta_1}\wt \rho \right)^{2}\left(\frac{1}{N\eta_1} \right)^{3}\frac{1}{\eta_1\sqrt{N}}\absv{R_s}^{p-3}.
		\end{equation} 
		 By a simple Young's inequality, together with $\eta_0/\eta_1\le N^{3\delta/2}$,
		 we then obtain a bound as in the rhs. of \eqref{eq:almthere}.

		For the the higher order terms we aim to estimate
		\beq\begin{aligned}
			\frac{\kappa_{n}}{n!}\sum_{a,b}\E\partial_{ab}^{n}\absv{R_s}^{p},\qquad n\geq 4.
		\end{aligned}\eeq
		 Here we chose the assignment $\kappa_n:=\kappa_n(ab,ab,\dots)$ for simplicity of notation, all other cases being analogous. 
		We expand the derivative as 
		\beq\label{eq:local law high order der expansion}
		\frac{\kappa_{n}}{n!}\sum_{a,b}\partial_{ab}^{n}\absv{R_s}^{p}=\frac{\kappa_{n}}{n!}\sum_{a,b}\sum_{k=1}^{n}\sum_{\substack{\bsn=(n_{1},\ldots, n_{k})\in\N^{k} \\ n_{1}+\cdots+n_{k}=n}}c_{\bsn,p}\prod_{i=1}^{k}\partial_{ab}^{n_{i}}[R_s]\absv{R_s}^{p-k},
		\eeq
		for some combinatorial numbers $c_{\bsn,p}$ depending on $\bsn\in\N^{k}$ and $p\in\N$.  For any $n_i \geq 1$ the partial derivative $\partial_{ab}R_{s}$ is the sum over a trace index $x$ over products over $n_{i}+1$ entries of $G$ of the form $G_{x\cdot}\cdots G_{\cdot x}$ where the middle terms are entries $G_{aa},G_{ba},G_{bb}$, or $G_{bb}$. We can apply Lemma \ref{lem:X+A monotone} to the middle terms, and Cauchy--Schwarz followed by two applications of the Ward identity and Lemma \ref{lem:X+A monotone} again to estimate \begin{equation}\label{eq:local law partial der estimate}
			\Absv{\partial_{ab}^{n_i}R_s }\prec \frac{1}{N\eta_1}\left(\frac{\eta_0}{\eta_1}\right)^{n_i-1}\left(\frac{\eta_0\wt\rho}{\eta_1}\right).
		\end{equation} Hence, \begin{equation}\label{eq:local law higher order bound}
			N^{-n/2}\sum_{a,b}\prod_{i=1}^{k}\Absv{\partial_{ab}^{n_{i}}[R_s]}\prec N^{2-n/2} \left( \frac{\eta_0}{\eta_1} \right)^{k-1}\left(\frac{\eta_0\wt\rho}{\eta_1}\right)\left(\frac{1}{N\eta_1}\right)^{k}.
		\end{equation}  Using $\kappa_n\lesssim N^{-n/2}$,  after applying Young's inequality to \eqref{eq:local law high order der expansion} and \eqref{eq:local law higher order bound} we see that  for $n\ge 4$ we have 
		\begin{equation}
			\frac{\kappa_{n}}{n!}\sum_{a,b}\E\partial_{ab}^{n}\absv{R_s}^{p}\prec N^{2-n/2}\left(\frac{ \eta_0}{\eta_1} \right)^{n}\wt\rho\left[\absv{R_s}^{p}+\left(\frac{1}{N\eta_1} \right)^{p} \right] \le N^{1/2+10\delta}\Psi(z,\eta_0)\left[\absv{R_s}^{p}+\left(\frac{1}{N\eta_1} \right)^{p} \right].
		\end{equation}
	\end{proof}
	
	\begin{prop}\label{prop:X+A derivative estimate isotropic}
		Let $\bsu,\bsy\in\C^{2N}$ and define \begin{equation}
			S_s(\bsu,\bsy)\deq\bsu^*\left(G^s(\ii\eta)-M_{A-z}(\ii \eta) \right)\bsy. 
		\end{equation} Let $p> 2$ be even  and let $\delta:=\epsilon_0/100$,  then
		 \begin{equation}
		\label{eq:isobGFT}
			\left|\E \frac{\dd}{\dd s}|S_s(\bsu,\bsy)|^{p} \right|\lesssim N^{\frac{1}{2}+ 10 \delta}\Psi(z,\eta_0)\left[\absv{S_s(\bsu,\bsy)}^p+\left(\Psi (z,\eta_0) \right)^p \right],
		\end{equation} for all $0\leq s\leq t$. 
	\end{prop}

	\begin{proof}[Proof of \eqref{eq:local law: zag isotropic}]
	
	Given \eqref{eq:isobGFT}, the proof of  \eqref{eq:local law: zag isotropic} is completely analogous to the one of \eqref{eq:local law: zag avaerage}, and so omitted.
	
	\end{proof}

	\begin{proof}[Proof of Proposition~\ref{prop:X+A derivative estimate isotropic}]
		Similar to Proposition \ref{prop:X+A derivative estimate isotropic} after applying It\^o's lemma and a cumulant expansion we see that  \begin{equation}
			\left|\E\frac{\dd}{\dd s}|S_s(\bsu,\bsy)|^{p} \right|\leq\left|\sum_{k=2}^{L-1}\sum_{a,a_1,\dots a_k,b,b_1,\dots,b_k}\frac{\kappa_s(ab,a_1b_1,\dots,a_kb_k)}{k!}\E\partial_{a_kb_k}\cdots\partial_{ab}|S_s(\bsu,\bsy)|^{p} \right|+\Omega(L). 
		\end{equation} 
		
		For simplicity we will write $G_{\bsu \bsy}=\bsu^*(G^{s}({\ii\eta}))\bsy$ and $S_s=S_s(\bsu,\bsy)$.
		  Here $\Omega(L)$ is a negligible error if $L$ is sufficiently large. 
We start with the estimate of the third order cumulants, i.e. we present the estimate of  $\frac{\kappa_3}{6}\sum_{a,b}\E \partial_{ab}^3 |S_s|^p$, with $\kappa_3:=\kappa_3(ab,ab,ab)$. The estimate of all the other terms coming from third order cumulants is completely analogous and so omitted.  Expanding the derivative, we obtain
\begin{equation}\begin{aligned}
				\frac{\kappa_3}{6}\sum_{a,b}\E \partial_{ab}^3 |S_s|^p=	\frac{\kappa_3}{6}\sum_{a,b}&\E p(\partial_{ab}^3 G_{\bsu \bsy}) |S_s|^{p-1}+\E c_p(\partial_{ab}^2 G_{\bsu\bsy})(\partial_{ab} G_{\bsu\bsy}) |S_s|^{p-2}\\&\quad+\E c'_p(\partial_{ab} G_{\bsu\bsy})^3 |S_s|^{p-3}.
			\end{aligned}
		\end{equation}   Note that $\kappa_3\lesssim N^{-3/2}$, and we begin with the term $\E p(\partial_{ab}^3 G_{\bsu \bsy}) |S_s(\bsu,\bsy)|^{p-1}$. $\partial_{ab}^3 G_{\bsu \bsy}$ is a linear combination of terms of the form $G_{\bsu a}G_{bb}G_{aa}G_{b\bsy}$, $G_{\bsu a}G_{ab}G_{ba}G_{b\bsy}$, $G_{\bsu a}G_{ab}^2G_{b\bsy}$, or similar terms up to a transposition of $a$ and $b$. The main observation is that we always have two off-diagonal terms at the beginning and end of the product. All middle terms can be bounded using Lemma \ref{lem:X+A monotone}. For the term with two diagonal elements \beq
		N^{-3/2}\sum_{a,b}G_{\bsu a}G_{bb}G_{aa}G_{b\bsy}\cdot \absv{S_s}^{p-1},
		\eeq
		we estimate the diagonal factors with Lemma \ref{lem:X+A monotone}:
		\beq\begin{split}
			N^{-3/2}\sum_{a,b}\absv{G_{\bsu a}G_{b\bsy}}\absv{G_{aa}G_{bb}}&\prec N^{-3/2+2\delta}
			\left(N\sum_{a}\absv{G_{\bsu a}}^{2}\right)^{1/2}\left(N\sum_{b}\absv{G_{b\bsy}}^{2}\right)^{1/2}	\\
			&\prec N^{-1/2+2\delta}\frac{{\eta_0}^{2}}{\eta_1^{2}}\frac{\wt \rho}{{\eta_0}}\leq N^{1/2+8\delta}\Psi({z},{\eta_0})^{2},
		\end{split}\eeq
		and hence we arrive at
		\beq
		N^{-3/2}\sum_{a,b}G_{\bsu a}G_{bb}G_{aa}G_{b\bsy}\cdot \absv{S_{s}}^{p-1}\prec N^{1/2+8\delta}\Psi( z,\eta_0)^{2}\absv{S_{s}}^{p-1}.
		\eeq
		As we can always apply Lemma \ref{lem:X+A monotone} twice to the middle two terms in the product of four $G_{\cdot\cdot}$ we conclude \begin{equation}\label{eq:third cumulant 3}
			N^{-3/2}\sum_{a,b}(\partial_{ab}^3G_{\bsu \bsy})\absv{S_{s}}^{p-1}\prec N^{1/2+8\delta}\Psi( z,\eta_0)^2 |S_{s}|^{p-1}.
		\end{equation} Similarly terms of the form $(\partial_{ab} G_{\bsu\bsy})^3$ can also be bounded. \begin{align*}
			N^{-3/2}\sum_{a,b}\absv{G_{\bsu a}G_{b\bsy}}^{3}\absv{S_{s}}^{p-3}\prec N^{-\frac{3}{2}+4\delta}\sum_{a,b}\absv{G_{\bsu a}}^2\absv{G_{b\bsy}}^2\absv{S_{s}}^{p-3}
			\prec N^{\frac{1}{2}+8\delta}\Psi( z_0,{\eta_0})^{4}\absv{S_{s}}^{p-3},
		\end{align*} and \begin{equation}\label{eq:third cumulant 1 cubed}
			N^{-3/2}\sum_{a,b}\left(\partial_{ab}G_{\bsu\bsy}\right)^3|S_{s}|^{p-3}\prec N^{\frac{1}{2}+8\delta}\Psi(  z,\eta_0)^{4}\absv{S_{s}}^{p-3}
		\end{equation} Again,  nearly identical computations yields \begin{equation}\label{eq:third cumulant 2+1}
			N^{-3/2}\sum_{a,b}(\partial_{ab} G_{\bsu\bsy} )(\partial_{ab}^2 G_{\bsu\bsy})\absv{S_{s}}^{p-2}\prec N^{\frac{1}{2}+10\delta}\Psi(  z,\eta_0)^{3}\absv{S_{s}}^{p-2}.
		\end{equation} Taking expected values of \eqref{eq:third cumulant 3}, \eqref{eq:third cumulant 1 cubed}, and \eqref{eq:third cumulant 2+1} we conclude 
		\begin{equation}
			\left|\frac{\kappa_3}{6}\sum_{a,b}\E \partial_{ab}^3 |S_s|^p \right|\lesssim  N^{1/2+10\delta}\Psi(z,\eta_0)\big[\E |S_{{s}}|^{p}+\left(\Psi (z,\eta_0)\right)^p \big] .
		\end{equation}		
		This concludes the estimate of the third order cumulants.


		 We now present the bound for cumulants of order four and higher, i.e.  our goal is to estimate $\frac{\kappa_{n}}{n!}\sum_{a,b}\E\partial_{ab}^{n}|S_{s}|^{p}$, for $n\geq 4$, where we again chose the assignment $\kappa_n:=\kappa_n(ab,ab,\dots)$ for simplicity of notation. 
		We expand the derivative as 
		\beq
		\frac{\kappa_{n}}{n!}\sum_{a,b}\partial_{ab}^{n}(G_{\bsu\bsy})^{p}=\frac{\kappa_{n}}{n!}\sum_{a,b}\sum_{k=1}^{n}\sum_{\substack{\bsn=(n_{1},\ldots, n_{k})\in\N^{k} \\ n_{1}+\cdots+n_{k}=n}}c_{\bsn,p}\prod_{i=1}^{k}\partial_{ab}^{n_{i}}[G_{\bsu\bsy}](S_{s})^{p-k},
		\eeq
		for some combinatorial numbers $c_{\bsn,p}$ depending on $\bsn\in\N^{k}$ and $p\in\N$.
		Next, we claim for each $k\in\llbra 1,n\rrbra$ and $\bsn\in\N^{k}$, with $n_1+\dots+n_k=n$, that
		\beq\label{eq:higher_goal}
		N^{-n/2}\sum_{a,b}\prod_{i=1}^{k}\Absv{\partial_{ab}^{n_{i}}[G_{xy}]}\prec \big[N^{1/2+10\delta}\Psi(z,\eta_0)\big]\cdot (\Psi(z,\eta_0))^{k}.
		\eeq
		To this end, without loss of generality we assume $n_{1}\leq n_{2}\leq \cdots\leq n_{k}$ and divide the proof in two cases, $n_{2}=1$ and $n_{2}>1$.
		
		Firstly when $n_{2}=1$, notice that 
		\beq
		\absv{\partial_{ab}^{n_{1}}[G_{\bsu\bsy}]\partial_{ab}^{n_{2}}[G_{\bsu\bsy}]}= \absv{G_{\bsu a}G_{b\bsy}+G_{\bsu b}G_{a\bsy}}^{2}\leq 2(\absv{G_{\bsu a}}^{2}\absv{G_{b\bsy}}^{2}+\absv{G_{\bsu b}}^{2}\absv{G_{a\bsy}}^{2}).
		\eeq
		For the remaining derivatives of order $n_{3},\ldots,n_{k}$ we apply Lemma \ref{lem:X+A monotone} so that
		\beq\label{eq:higher_rough}
		\absv{\partial_{ab}^{n}[G_{\bsu\bsy}]}\prec N^{(n+1)\delta},\qquad n\in\N.
		\eeq
		In conclusion we have
		\beq\begin{aligned}
			N^{-n/2}\sum_{a,b}\prod_{i=1}^{k}\Absv{\partial_{ab}^{n_{i}}[G_{\bsu\bsy}]}
			\prec&N^{-n/2}N^{(n+k-4)\delta}\frac{\im G_{\bsu\bsu}\im G_{\bsy\bsy}}{\eta_1^{2}}\prec N^{-n/2}N^{(n+k-4)\delta}\frac{{\eta_0}^{2}}{\eta_1^{4}}\wt\rho^{2}\\
			&=N^{-n/2+2}N^{(n+k)\delta}\Psi({z},{\eta_0})^{4}.
		\end{aligned}\eeq
		Finally, we use $N^{-1/2+\epsilon_0/2}\le N^{-1/2}\widetilde{\mathfrak{s}}_0^{-1/2}=\Psi({z},{\eta_0})$ so that 
		\beq\begin{aligned}
			N^{-n/2+2}N^{(n+k)\delta}\Psi({z},{\eta_0})^{4}\leq \sqrt{N}N^{2n\delta- (n-4)\epsilon_0/2}\Psi( z,\eta_0)\cdot (\Psi( z,\eta_0))^{n}\le \big[N^{1/2+8\delta}\Psi(z,\eta_0)\big] (\Psi(z,\eta_0))^k,
		\end{aligned}\eeq
		where we also used that $n\ge k$ and $\Psi(z,\eta_0)\le 1$. This concludes the proof of \eqref{eq:higher_goal} for $n_{2}=1$.
		
		We next prove \eqref{eq:higher_goal} when $n_{2}>1$. Here we will only use the first factor corresponding to $n_{1}$, and applying using Lemma \ref{lem:X+A monotone} we have that
		\beq
		\absv{\partial_{ab}^{n_{1}}[G_{\bsu\bsy}]}\prec N^{(n_{1}-1)\delta}(\absv{G_{\bsu a}G_{b\bsy}}+\absv{G_{\bsu b}G_{a\bsy}})\leq \frac{1}{2}N^{(n_{1}-1)\delta}\sum_{i\in\{\bsu,\bsy\},j\in\{a,b\}}\absv{G_{ij}}^{2}.
		\eeq
		Applying \eqref{eq:higher_rough} to the remaining derivatives, using again $N^{-1/2+\epsilon_0/2}\le \Psi(z,\eta_0)$, we find
		\beq
		\begin{split}
			N^{-n/2}\sum_{a,b}\prod_{i=1}^{k}\Absv{\partial_{ab}^{n_{i}}[G_{\bsu,\bsy}]}
			&\prec N^{-n/2+2}N^{(n+k)\delta-(n-4)\epsilon_0/2}\Psi(z,\eta_0)^{2}\leq N^{1/2+8\delta}\Psi( z,\eta_0)^{n-1} \\
			&\le \big[N^{1/2+8\delta} \Psi(z,\eta_0)\big] (\Psi(z,\eta_0))^k.
		\end{split}
\eeq
	We point out that in the last inequality we used that $n_{2}>1$ implies $n\geq 2(k-1)+n_{1}\geq 2k-1$, which together with $n\geq 4$ gives $n\geq (2k+3)/2$ so that $n\geq k+2$. This concludes the proof of \eqref{eq:higher_goal}.
	\end{proof}

\appendix
\section{Preliminaries on $A+\bsx$}\label{sec:meta}
\subsection{Freeness of $A$ and $\bsx$: Proof of Lemma \ref{lem:meta}}
\begin{proof}[Proof of Lemma \ref{lem:meta}]
	We need to prove two statements, that (i) $\bsx$ is a circular element, and that (ii) $A$ and $\bsx$ are $*$-free. 
	
	To prove the first claim, we approximate (in $*$-moments) $\bsx_{ij}$ with independent complex Ginibre matrices $X^{(ij)}$ of size $N'\gg 1$, so that $X'=(X^{(ij)})_{ij}\in\C^{(NN'\times NN')}$ is a complex Ginibre matrix of size $NN'$. Notice that $\sqrt{2}\re X'$ and $\sqrt{2}\im X'$ are independent Gaussian unitary ensembles of size $NN'$, so that by their asymptotic freeness for any $*$-polynomial $p$ of two variables we have (see \cite[Theorem 4, Chapter 4]{Mingo-Speicher2017})
	\begin{equation}\label{eq:meta}
		\lim_{N'\to\infty} \frac{1}{NN'}\Tr p(\re X',\im X')=\brkt{p(\sqrt{2}\bss_{1},\sqrt{2}\bss_{2})}_{\caM},
	\end{equation}
	where $\bss_{1},\bss_{2}$ are free semicircular elements in $(\caM,\brkt{\cdot}_{\caM})$. On the other hand, we can expand the left-hand side of \eqref{eq:meta} as another $*$-polynomial of the $(N'\times N')$ GUE's $\{\re X^{(i,j)},\im X^{(i,j)}:1\leq i,j\leq N\}$. More precisely, we write
	\begin{equation}
		\wt{p}_{i}:\caM^{N\times N}\times\caM^{N\times N}\to\caM,\qquad \wt{p}_{i}(\bsa,\bsb)=(p(\bsa,\bsb))_{ii}\in\caM,
	\end{equation}
	where $(p(\bsa,\bsb))_{ii}$ stands for the $(i,i)$-th entry of $p(\bsa,\bsb)\in\caM^{N\times N}$ (which obviously is a $*$-polynomial in the entries of $\bsa,\bsb$). Then, by the asymptotic freeness of $\{\re X^{(i,j)},\im X^{(i,j)}:1\leq i,j\leq N\}$ we get
	\begin{equation}\begin{aligned}\label{eq:meta1}
			&\lim_{N'\to\infty}\frac{1}{NN'}\Tr p(\re X',\im X')=\frac{1}{N}\sum_{i=1}^{N}\lim_{N'\to\infty}\frac{1}{N'}\Tr \wt{p}_{i}(\re X',\im X')	\\
			=\,&\frac{1}{N}\sum_{i=1}^{N}\brkt{\wt{p}_{i}(\re\bsx,\im\bsx)}_{\caM}=\brkt{\brkt{p(\re \bsx,\im\bsx)}_{\caM}}.
	\end{aligned}\end{equation}
	Combining \eqref{eq:meta} and \eqref{eq:meta1} proves that $\sqrt{2}\re \bsx$ and $\sqrt{2}\im\bsx$ are free semi-circular elements, so that $\bsx$ is a circular element.
	
	For the second claim, we use the fact that the two $(NN'\times NN')$ matrices $A\otimes I_{N'}$ and $X'$ are asymptotically $*$-free as $N'\to\infty$, and that their joint $*$-moments converge to those of $A$ and $\bsx$. We omit further details for brevity.
\end{proof}

\subsection{The density of $\rho_{A+\bsx}$ around an edge point $z_{0}$}
Here we collect preliminary results on the Brown measure density $\rho_{A+\bsx}$. The content of this section applies to generic operators $A,\bsx$ in a von Neumann algebra $\wt{\caM}$, and we apply them to our choice of $A,\bsx\in\caM^{N\times N}$ in Definition \ref{defn:vN}.
Firstly, we have the following description of the support:
\begin{lem}\label{lem:supp}
	Let $A$ and $\bsx$ be $*$-free operators in a $W\adj$-probability space $(\wt{\caM},\brkt{\cdot})$ such that $\bsx$ is a circular element. Then the following hold true:
	\begin{itemize}
		\item[(i)] $\supp\rho_{A+\bsx}=\ol{\{z\in\C:\brkt{\absv{A-z}^{-2}}>1\}}$.
		\item[(ii)] For $z\in\C$ with $\norm{(A-z)^{-1}}<\infty$, $z\in\partial\supp\rho_{A+\bsx}$ if and only if $\brkt{\absv{A-z}^{-2}}=1$.
	\end{itemize}
\end{lem}
\begin{proof}
	The first statement is due to \cite[Theorem 4.6]{Zhong2021}. The second statement follows from the first statement and the fact the function $z\mapsto\brkt{\absv{A-z}^{-2}}$ is a continuous, strictly subharmonic function (since $\partial_{z}\partial_{\ol{z}}\brkt{\absv{A-z}^{-2}}=\brkt{(A-z)^{-2}((A-z)^{-2})\adj}>0$) on the resolvent set $\{z\in\C:\norm{(A-z)^{-1}}<\infty\}$.
\end{proof}

Next, we record the main result of \cite{Erdos-Ji2023circ}: Below we copied \cite[Theorem 2.10, (iii) and (iv)]{Erdos-Ji2023circ} verbatim, except that we took $t=1$, made the theorem local and quantitative (see \cite[Remark 2.12]{Erdos-Ji2023circ}), wrote $A$ instead of $\bsa$, and spelled out the function $f_{A}(z)\deq \brkt{\absv{A-z}^{-2}}$ and its derivatives; e.g. we have used the trivial algebraic identity
\begin{equation}\begin{aligned}
		\absv{\nabla f_{A}(z)}^{2}
		&=	\absv{(\partial_{z}+\partial_{\ol{z}})f_{A}(z)}^{2}+\absv{(\partial_{z}-\partial_{\ol{z}})f_{A}(z)}^{2}	\\
		&=2(\absv{\partial_{z}f_{A}(z)}^{2}+\absv{\partial_{\ol{z}}f_{A}(z)}^{2})=4\absv{\brkt{(A-z)\absv{A-z}^{-4}}}^{2}.
\end{aligned}\end{equation}
\begin{lem}[Theorem 2.10 of \cite{Erdos-Ji2023circ}]\label{lem:sharp}
	Let $A$ and $\bsx$ be two operators in a $W\adj$-probability space $(\caM,\brkt{\cdot})$ such that $A$ and $\bsx$ are $*$-free and $\bsx$ is a circular element. Suppose that $\norm{A}\leq C_{0}$ for a constant $C_{0}>0$. Define
	\begin{equation}
		\caD\deq \{z\in\C:\brkt{\absv{A-z}^{-2}}>1\},
	\end{equation}
	\beq\label{eq:def_C}
	\caC\deq\left\{z\in\partial\caD:\brkt{(A-z)\absv{A-z}^{-4}}=0\right\}.
	\eeq
	Let $z_{0}\in\partial\caD$ satisfy $\norm{(A-z_{0})^{-1}}\leq C$ for a constant $C>0$.
	Then there exists a density $\rho$ for $\rho_{\bsa+\sqrt{t}\bsx}$ that satisfies the following:
	\begin{itemize}
		\item[(iii)] {\rm(Sharp edge)} Suppose that $\absv{\brkt{(A-z_{0})\absv{A-z_{0}}^{-4}}}>1/C$. Then the density $\rho(z)$ satisfies the following asymptotics as $z\to z_{0}$;
		\beq\label{eq:edge_sharp}
		\rho(z)=\lone_{\caD}(z)\left(\frac{1}{\pi}\brkt{\absv{A-z_{0}}^{-4}}^{-1}\absv{\brkt{(A-z_{0})\absv{A-z_{0}}^{-4}}}^{2}+O(\absv{z-z_{0}})\right),
		\eeq
		where the implicit constant depends only on $C$ and $C_{0}$.
		
		\item[(iv)] {\rm(Quadratic edge)} Define the $(2\times 2)$ Hessian matrix
		\begin{equation}
			\caH_{z_{0}}=H_{z}[\brkt{\absv{A-z}^{-2}}]_{z=z_{0}}
		\end{equation} 
		where $z\mapsto\brkt{\absv{A-z}^{-2}}$ is viewed as a function in $\R^{2}\cong\C$. Then $\Tr \caH_{z_{0}}>c_{1}$ for some constant $c_{1}>0$ depending only on $C_{0}$ and $C$, and the density $\rho(z)$ satisfies the following quadratic asymptotics as $z\to z_{0}$;
		\beq\label{eq:edge_quad}
		\rho(z)= \lone_{\caD}(z)\left(Q_{z_{0}}[z-z_{0}]+O(\absv{z-z_{0}}^{3})\right),
		\eeq
		where the implicit constant depends only on $C$ and $C_{0}$, and $u\mapsto Q_{z_{0}}[u]$ is the real-valued quadratic form on $\C$ defined by\footnote{In \eqref{eq:quad_def}, we identified $u\in\C$ with $(\re u,\im u)^{\tp}\in\R^{2}$ and $\caH_{z_{0}}$ with the corresponding $2\times 2$ real matrix. This choice is purely cosmetic as $\brkt{u,\caH_{z_{0}}u}$ and $\norm{\caH_{z_{0}}u}$ do not depend on the choice of an orthonormal basis of $\C\cong\R^{2}$.}
		\beq\label{eq:quad_def}
		Q_{z}[u]\deq \frac{1}{2\pi}\frac{\brkt{(A-z)^{-2}((A-z)^{-2})\adj}}{\brkt{\absv{A-z}^{-4}}}\brkt{u,\caH_{z_{0}}u}
		+\frac{1}{4\pi}\frac{1}{\brkt{\absv{A-z}^{-4}}}\norm{\caH_{z_{0}}u}^{2}.
		\eeq
	\end{itemize}
\end{lem}
\subsection{Proof of Lemma \ref{lem:fc}}
Recall $f_{A}(z)=\brkt{\absv{A-z}^{-2}}$. Before proceeding with the proof of Lemma \ref{lem:fc}, we present the following estimate for the fluctuation scale $\sigma_{\rmf}(z)$ for an edge point $z$ in terms of $\absv{\nabla f_{A}(z)}$, that is, we generalize the naive estimates $\sigma_{\rmf}(z)\sim N^{-1/2}$ at a sharp edge and $\sigma_{\rmf}(z)\sim N^{-1/4}$ at a critical edge to the whole edge:
\begin{lem}\label{lem:rmf}
	Under Assumption \ref{ass:newassout}, the following holds uniformly over $z\in\partial\supp\rho_{A+\bsx}$:
	\begin{equation}\label{eq:rmf}
		\sigma_{\rmf}(z)\sim \frac{N^{-1/2}}{N^{-1/4}+\absv{\nabla f_{A}(z)}}.
	\end{equation}
	In particular $\sigma_{\rmf}(z)\lesssim N^{-1/4}$ holds uniformly over $z\in\partial\supp\rho_{A+\bsx}$. 
\end{lem}
\begin{proof}
	Note that \eqref{eq:rmf} would follow immediately once we have
	\begin{equation}\label{eq:rmf_e}
		\frac{1}{N}=\int_{\absv{w-z}\leq\sigma_{\rmf}(z)}\rho_{A+\bsx}(w)\dd^{2}w\sim \absv{\nabla f_{A}(z)}^{2}\sigma_{\rmf}(z)^{2}+\sigma_{\rmf}(z)^{4},
	\end{equation}
	by solving the quadratic equation in $\sigma_{\rmf}(z)^{2}$. In what follows we now prove \eqref{eq:rmf_e}. 
	
	First, we recall the Taylor expansion of the density $\rho_{A+\bsx}(w)$ in \cite[Eq. (5.18)]{Erdos-Ji2023circ} which reads as, for $v\equiv v(w,0)=\brkt{M_{A-w}(+\ii 0)}$ from Lemma \ref{lem:MDE_asymp},
	\begin{equation}\begin{aligned}\label{eq:rho_e}
			\pi &\rho_{A+\bsx}(w)=v^{2}\Delta f_{A}(w)	\\
			&+(\brkt{\absv{A-w}^{-4}}-2v^{2}\brkt{\absv{A-w}^{-6}})^{-1}\Absv{\nabla f_{A}(w)-2v^{2}\nabla \brkt{\absv{A-w}^{-4}}}^{2}+O(v^{4}).
	\end{aligned}\end{equation}
	Notice that around an edge point $z$ we have $v\ll 1$, so that the expansion is valid.  Next, we estimate the integral of \eqref{eq:rho_e} over the disk $\absv{w-z}\leq\sigma_{\rmf}(z)$ to prove \eqref{eq:rmf_e}. For the first term of \eqref{eq:rho_e}, expanding $w$ around $z$ and plugging in $v^{2}\sim (f_{A}(w)-1)_{+}$ from \eqref{lem:MDE_asymp} and $\Delta f_{A}(w)=4\brkt{\absv{(A-w)^{2}}^{-2}}\sim 1$, we have
	\begin{equation}\begin{aligned}
			v^{2}\Delta f_{A}(w)
			&\sim (f_{A}(w)-1)_{+}	=(f_{A}(w)-f_{A}(z))_{+}	\\
			&=\left(\brkt{\nabla f_{A}(z),(w-z)}+\brkt{\caH_{z}(w-z),(w-z)}\right)_{+}+O(\absv{w-z}^{3}).
	\end{aligned}\end{equation}
	This in turn shows that as $\sigma\to0^{+}$
	\begin{equation}\label{eq:rho_e_1}
		\int_{\absv{w-z}\leq\sigma}v^{2}\Delta f_{A}(w)\dd^{2}w\sim \absv{\nabla f_{A}(z)}\sigma^{3}+\sigma^{4},
	\end{equation}
	where the lower bound follows from integrating over a sector of $\absv{w-z}\leq \sigma$ on which
	\begin{equation}\label{eq:sec}
		\frac{\brkt{\nabla f_{A}(z),(w-z)}}{\absv{w-z}}\geq c\absv{\nabla f_{A}(z)}, \qquad 
		\frac{\brkt{\caH_{z}(w-z),(w-z)}}{\absv{w-z}^{2}}\geq c
	\end{equation}
	for some small constant $c\sim 1$. We can choose the angle of the sector to be $\gtrsim 1$, since the first inequality in \eqref{eq:sec} is satisfied in a sector of angle $\pi-O(c)$ and the second one is true in a double cone of angle $\sim 1$ by $\Delta f_{A}(z)\sim 1$, so that the intersection of the two sets must contain a sector of angle $\gtrsim 1$. 
	Likewise for the second term of \eqref{eq:rho_e}, notice that for $\absv{w-z}\ll1$
	\begin{equation}
		v^{2}\sim (f_{A}(w)-1)_{+}\leq \absv{\nabla f_{A}(z)}\absv{w-z}+O(\absv{w-z}^{2}),
	\end{equation}
	and similarly we expand
	\begin{equation}
		\nabla f_{A}(w)=\nabla f_{A}(z)+O(\absv{w-z}).
	\end{equation}
	Thus we find a constant $C>0$ for which the second term of \eqref{eq:rho_e} is comparable to $\absv{\nabla f_{A}(w)}^{2}$ (i.e. $v^{2}$ is negligible) provided $C\absv{w-z}<\absv{\nabla f_{A}(z)}$, using that the first factor therein is trivially of order $1$.
	Then we divide the proof into two cases, according to $\absv{\nabla f_{A}(z)}>C\sigma_{\rmf}(z)$ or not. In the former case the integral of the second term of \eqref{eq:rho_e} is bounded from below by
	\begin{equation}
		\int_{\absv{w-z}\leq\sigma_{\rmf}(z)}\absv{\nabla f_{A}(z)}^{2}\dd^{2}w\sim \absv{\nabla f_{A}(z)}^{2}\sigma_{\rmf}(z)^{2}.
	\end{equation}
	In the latter case when $\absv{\nabla f_{A}(z)}\leq C\sigma_{\rmf}(z)$ we only need an upper bound of $O(\absv{\sigma_{\rmf}(z)}^{4})$ for the integral, since \eqref{eq:rho_e_1} already establishes the lower bound required in \eqref{eq:rho_e} and the second term of \eqref{eq:rho_e} is positive. This upper bound follows directly from Taylor expansions, concluding the proof of \eqref{eq:rmf}.
\end{proof}
\begin{proof}[Proof of Lemma \ref{lem:fc}]
	Notice that $f_{A}$ is a strictly subharmonic function on $\C\setminus\mathrm{Spec}(A)$ since $\Delta f_{A}(z)=4\brkt{\absv{(A-z)^{2}}^{-2}}\sim 1$. Therefore by the maximum principle 
	applied to each component of the complement of $\mathrm{Spec}_\epsilon(A)$ (each of them is included in 
	$\C\setminus\mathrm{Spec}(A))$,	
	we only need to show $f_{A}(z)\leq 1-N^{-1/2+c_{*}\epsilon}$ for $z\in\partial\mathrm{Spec}_{\epsilon}(A+\bsx)$.
	
	Now we fix a point $z\in\partial\mathrm{Spec}_{\epsilon}(A+\bsx)$, and take $\wt{z}\in\partial\supp \rho_{A+\bsx}$ 
	to be (one of) the closest point to $z$, i.e. such that
	\begin{equation}
		\absv{z-\wt{z}}=\min\{\absv{z-w}:w\in\partial\supp\rho_{A+\bsx}\}.
	\end{equation}
	Then it follows that $N^{-1/4+\epsilon}\geq \absv{z-\wt{z}}\geq N^{\epsilon}\sigma_{\rmf}(\wt{z})$, 
	and that the disc $\{z+w:\absv{w}\leq \absv{z-\wt{z}}\}$ is tangent to the set  $\supp\rho_{A+\bsx}$ at $\wt{z}$. This already shows that $\wt{z}$ cannot be a critical edge, since no circle can be tangent to the level set of a smooth function at a saddle point, and a critical edge is indeed a saddle point as $f_{A}$ is subharmonic and an edge of $\rho_{A+\bsx}$ cannot be a local minimum. 	In fact, we claim the following concrete lower bound for $\absv{\nabla f (\wt{z})}$ for some constant $c>0$;
	\begin{equation}\label{eq:not_crit}
		\absv{\nabla f(\wt{z})}\geq c\absv{z-\wt{z}}.
	\end{equation}
	To see \eqref{eq:not_crit}, we Taylor-expand $f_{A}$ around $\wt{z}$ so that
	\begin{equation}\label{eq:Taylor}
		f_{A}(\wt{z}+w)=1+\brkt{\nabla f_{A}(\wt{z}),w}+ \brkt{\caH_{\wt{z}}w,w}+O(\absv{w}^{3}).
	\end{equation}
	Then notice that the Hessian $\caH_{\wt{z}}$ of $f_A$ at $\wt{z}$
	is guaranteed to have one positive eigenvalue of size 
	at least $\Delta f_{A}(\wt{z})/2\sim 1$, so that there exist a constant $c'>0$ and a double cone $\caC$ of angle $\theta\sim 1$ based at the origin around the leading eigenvector of $\caH_{\wt{z}}$ such that 
	\begin{equation}
		\brkt{\caH_{\wt{z}}w,w}\geq c'\absv{w}^{2},\qquad w\in \caC.
	\end{equation}
	Plugging this into \eqref{eq:Taylor} shows $f_{A}(\wt{z}+w)>1$ 
	holds whenever $w$ is in the set
	\begin{equation}
		\caC\cap \left\{w:\frac{2\absv{\nabla f_{A}(\wt{z})}}{c'}\leq \absv{w}\leq \frac{c'}{2}\right\}.
	\end{equation}
	Since the set above 
	should not intersect with the disk 
	$\{w:\absv{w}<\absv{z-\wt{z}}\}$ 
	and the angle $\theta$ of the cone is bounded from below, we conclude
	\begin{equation}\label{eq:nab}
		\absv{\nabla f_{A}(\wt{z})}\geq \frac{c'\theta}{100} \absv{z-\wt{z}}\geq c\absv{z-\wt{z}}.
	\end{equation}
	As a consequence of \eqref{eq:not_crit}, we combine \eqref{eq:nab}, \eqref{eq:rmf}, and the assumption $\absv{\wt{z}-z}\geq N^{\epsilon}\sigma_{\rmf}(\wt{z})$ to deduce
	\begin{equation}
		\absv{\nabla f_{A}(\wt{z})}\gtrsim \absv{z-\wt{z}}\geq N^{\epsilon}\sigma_{\rmf}(\wt{z})\gtrsim N^{\epsilon}\frac{N^{-1/2}}{N^{-1/4}+\absv{\nabla f_{A}(\wt{z})}},
	\end{equation}
	which proves $\absv{\nabla f_{A}(\wt{z})}\gtrsim N^{-1/4+\epsilon/2}$.
	
	We then divide the proof into two cases, depending on whether $\absv{\nabla f_{A}(z)}\leq N^{\delta} \absv{z-\wt{z}}$ or not, where we set $\delta=\epsilon/100$. In the former case, we use the assumption $\absv{\det \caH_{\wt{z}}}\sim 1$ from \eqref{eq:nondegenerate} to find a point $\wh{z}$ with $\absv{z-\wh{z}}\lesssim \absv{\nabla f_{A}(z)}$ such that $\nabla f_{A}(\wh{z})=0$. Then the function $f_{A}(\wh{z}+\cdot)$ is approximately a non-degenerate quadratic form, since $\absv{\det\caH_{\wt{z}}}\sim 1$ and $\absv{\wt{z}-\wh{z}}\lesssim N^{\delta}\absv{z-\wt{z}}\lesssim N^{-1/4+\epsilon+\delta}$:
	\begin{equation}
		f_{A}(\wh{z}+w)=f_{A}(\wh{z})+\brkt{\caH_{\wh{z}}w,w} +O(\absv{w}^{3}),\qquad w\in\C,\absv{w}\ll 1.
	\end{equation}
	Notice that the point $\wt{z}-\wh{z}$ is on the level set of $f_{A}(\wh{z}+\cdot)$, and $z-\wh{z}$ is  
	away from this level set at a distance $O(N^{\delta}\absv{z-\wt{z}})$.
	Since $\caH_{\wh{z}}$ is  non-degenerate the cubic error term is negligible. An explicit calculation with 
	the non-degenerate quadratic function (both in the hyperbolic and the (easier) elliptic cases) then shows that 
	\begin{equation}
		1-f_{A}(z)=f_{A}(\wt{z})-f_{A}(z)\gtrsim \absv{z-\wt{z}}^{2}\geq N^{-2\delta}\absv{\nabla f_{A}(\wt{z})}^{2}\gtrsim N^{-1/2+\epsilon-2\delta}.
	\end{equation}
	
	In the other case when $\absv{\nabla f_{A}(z)}\geq N^{\delta}\absv{z-\wt{z}}$, we have $\absv{\nabla f_{A}(\wt{z})}\gtrsim N^{\delta}\absv{z-\wt{z}}$ since the Hessian is bounded. Then we expand around $\wt{z}$ to conclude\begin{equation}\begin{aligned}
			1-f_{A}(z)=f_{A}(\wt{z})-f_{A}(z)&=\absv{\nabla f_{A}(\wt{z})}\absv{z-\wt{z}}+O(\absv{z-\wt{z}}^{2})	\\
			&\gtrsim \absv{\nabla f_{A}(\wt{z})}\frac{N^{-1/2+\epsilon}}{N^{-1/4}+\absv{\nabla f_{A}(\wt{z})}}\\
			&\gtrsim N^{-1/2+\epsilon},
	\end{aligned}\end{equation}
	where we used Lagrange multiplier theorem in the first, linear term to show that the vectors (in $\R^{2}\cong \C$) $\nabla f_{A}(\wt{z})$ and $z-\wt{z}$ are in the opposite directions.
\end{proof}

\end{document}